\documentclass[12pt]{report}
\usepackage{amstex}
\title{Nahm transform \\ of doubly-periodic instantons}
\author{Marcos Benevenuto Jardim \\ St. Anne's College \\ 
        University of Oxford}
\date{\vfill Thesis submitted for the degree of Doctor of Philosophy \\ Trinity Term 1999}

\newcommand{\pf}{{\em Proof: }}
\newcommand{\pfend}{\hfill $\Box$ \linebreak}
\newcommand{\eq}{\begin{equation}}
\newcommand{\matriz}{\left( \begin{array}{cc}}
\newcommand{\column}{\left( \begin{array}{c}}
\newcommand{\matfim}{\end{array} \right)}
\newcommand{\seta}{\rightarrow}
\newcommand{\torus}{T\times\cpx} \newcommand{\tproj}{T\times\proj}
 
\newcommand{\dual}{\hat{T}} \newcommand{\dproj}{\dual\times\proj}
\newcommand{\del}{\overline{\partial}}

\newcommand{\see}{\Leftrightarrow}
\newcommand{\imply}{\Rightarrow} \newcommand{\hproj}{\hat{{\Bbb P}}^1}

\newcommand{\ksi}{\xi} 
\newcommand{\real}{{\Bbb R}} \newcommand{\cpx}{{\Bbb C}}
\newcommand{\zed}{{\Bbb Z}} 
\newcommand{\proj}{{\Bbb P}^1} \newcommand{\vv}{{\cal V}}
\newcommand{\ee}{{\cal E}}  \newcommand{\pp}{{\cal P}}
  \newcommand{\call}{{\cal L}}
\newcommand{\cala}{{\cal A}} \newcommand{\calg}{{\cal G}} \newcommand{\jj}{{\cal J}}
\newcommand{\ind}{{\rm index}}  
\newcommand{\modspc}{{\cal M}_{(k,\ksi_0)}} \newcommand{\as}{\pm\ksi_0}
\newcommand{\oo}{{\cal O}}  \newcommand{\rr}{{\cal R}}
 
 \newcommand{\cfg}{{\cal A}_{(k,\ksi_0)}}

\newtheorem{mthm}{Main Theorem}
\newtheorem{thm}{Theorem}[chapter]
\newtheorem{cor}[thm]{Corollary}
\newtheorem{lem}[thm]{Lemma}
\newtheorem{prop}[thm]{Proposition}
\newtheorem{conj}[thm]{Conjecture}

\begin{document}
\maketitle \newpage
\renewcommand{\thepage}{\roman{page}}
\baselineskip18pt

\centerline{{\large \bf Nahm transform of doubly-periodic instantons}}
\vskip10pt
\centerline{Marcos Jardim}
\centerline{St. Anne's College}
\vskip10pt
\centerline{Thesis submitted for the degree of Doctor of Philosophy}
\centerline{Trinity Term 1999}
\vskip10pt
\centerline{{\bf ABSTRACT}}
\vskip10pt

\baselineskip15pt
This work concerns the study of certain finite-energy
solutions of the anti-self-dual Yang-Mills equations on
Euclidean 4-dimensional space which are periodic in two directions,
so-called {\em doubly-periodic instantons}. We establish a
circle of ideas involving equivalent analytical and
algebraic-geometric descriptions of these objects.

In the first introductory chapter we provide an overview of the
problem and state the main results to be proven in the thesis.

In chapter 2, we study the asymptotic behaviour of the connections
we are concerned with, and show that the coupled Dirac operator is
Fredholm.

After laying these foundations, we are ready to address the main
topic of the thesis, the construction of a {\em Nahm transform} of
doubly-periodic instantons. By combining differential-geometric and
holomorphic methods, we show in chapters 3 through 5 that
doubly-periodic instantons correspond bijectively to certain
singular {\em Higgs pairs}, i.e. meromorphic solutions of Hitchin's
equations defined over an elliptic curve.

The circle of ideas is finally closed in chapter 7. We start by
presenting a construction due to Friedman, Morgan \& Witten that
associates to each doubly-periodic instanton a {\em spectral pair}
consisting of a Riemann surface plus a line bundle over it. On the
other hand, it was shown by Hitchin that Higgs pairs are equivalent
to a similar set of data. We show that the Friedman, Morgan
\& Witten spectral pair associated with a doubly-periodic instanton
coincides with the Hitchin spectral pair associated with its Nahm
transform.

\newpage 
\baselineskip18pt
\vskip50pt

\centerline{{\bf ACKNOWLEDGEMENTS}}

\vskip50pt

First of all, I must thank my supervisors, Simon Donaldson and
Nigel Hitchin, for their constant support and guidance. Not only
they have taught me a lot of mathematics during the last three
years, but they also showed me what it is to be a mathematician.

\vskip20pt

I also benefited from conversations with various other
mathematicians; Alexei Kovalev, Brian Steer, Olivier Biquard and
Richard Thomas are just some of the people whose comments where
incorporated in the text. I also thank Jos\'e Nat\'ario for his
valuable friendship throughout my PhD years.

\vskip20pt

I was supported by CNPq, an agency from the Brazilian  Ministry of Education; many thanks for all the helpful  people that work there. 
I also thank the Mathematical Institute and St. Anne's College for 
travel grants at various stages of my research.

\vskip20pt

I owe much of what I have achieved to my parents, who have put me
and kept me always in the right track.

\vskip20pt

But this work is wholly dedicated to Anelyssa. For every single
moment she has spent with me. For every smile and spark of light
coming out her eyes. For making me undeservingly happy. \`A minha
Princesa Cora\c c\~ao, Nani Pequenina cara menina, com todo o meu
amor!

\newpage \tableofcontents \newpage
\setcounter{page}{1}
\renewcommand{\thepage}{\arabic{page}}
\pagestyle{plain}

\chapter{Overview and statement of the results} \label{intro} 

Since the appearance of the Yang-Mills equation on the mathematical
scene in the late 70's, its anti-self-dual (ASD) solutions have
been intensively studied. The first major result in the field was the
ADHM construction of instantons on $\real^4$ \cite{ADHM}. Soon
after that, W. Nahm adapted the ADHM construction to obtain the
{\em time-invariant} ASD solutions of the Yang-Mills equations, the
so-called monopoles \cite{N}. It turns out that these constructions
are two examples of a much more general framework.

The {\em Nahm transform} can be defined in general for ASD
connections on $\real^4$, which are invariant under some sub-group
of translations $\Lambda\subset\real^4$. In these generalised
situations, the Nahm transform gives rise to {\em dual instantons}
on $(\real^4)^*$, which are invariant under:
$$ \Lambda^*=\{\alpha\in(\real^4)^*\ |\ \alpha(\lambda)\in\zed\ \forall\lambda\in\Lambda\} $$
There are plenty of examples of such constructions available
in the literature, namely:

\begin{itemize}
\item The trivial case $\Lambda=\{0\}$ is closely related to the
celebrated ADHM construction of instantons, as described by Donaldson
\& Kronheimer \cite{DK}; in this case, $\Lambda^*=(\real^4)^*$ and an
instanton on $\real^4$ corresponds to some algebraic data.

\item If $\Lambda=\zed^4$, this is the Nahm transform of Braam \& van
Baal \cite{BVB} and Donaldson \& Kronheimer \cite{DK}, defining a
hyperk\"ahler isometry of the moduli space of instantons over two
dual 4-tori.

\item $\Lambda=\real$ gives rise to monopoles, extensively studied by
Hitchin \cite{H3}, Donaldson \cite{D3} and Hurtubise \& Murray
\cite{HM}, among several others; here, $\Lambda^*=\real^3$, and the
transformed object is, for SU(2) monopoles, an analytic solution of
certain matrix-valued ODE's (the so-called Nahm's equations), defined
over the open interval $(0,2)$ and with simple poles at the end-points.

\item $\Lambda=\zed$ correspond to the so-called calorons, studied
by Nahm \cite{N}, Garland \& Murray \cite{GM} and others; the
transformed object is the solution of certain nonlinear Nahm-type
equations on a circle.

\end{itemize}

The purpose of this work fits well into this larger mathematical
programme. We study instantons, i.e. finite energy solutions of the
Yang-Mills anti-self-dual equations, on $SU(2)$ bundles
$E\seta\torus$, which can be seen as solutions over $\real^4$
invariant under a two-dimensional lattice. More precisely, we
search for a definition of a Nahm transform in this situation.

According to the general scheme outlined above, the dual object
should be an instanton over $(\real^4)^*$ invariant under
$\Lambda^*=\zed^2\times\real^2$. This is the same as a solution of
the so-called Hitchin's equations \cite{H} over a two-dimensional
torus $\dual$, which we call the {\em dual torus}. Indeed, our
first main result, theorem \ref{nahmthm} below, addresses such a
correspondence. As in the case of monopoles, some singularities
appear \cite{H3}, essentially due to the non-compactness of
$\torus$.

Although the moduli space of singular solutions of Hitchin's
equations is relatively well studied \cite{Ko} \cite{K}, nothing
has been said about the moduli space of doubly-periodic instantons.
This is actually one of the main advantages of our approach, since we
can then use known information about the moduli space of Higgs pair to
probe the structure of the instanton moduli. In particular, existence
of Higgs pairs will imply existence of doubly-periodic instantons.

We then move to a more traditional approach and study this moduli
space within the usual framework of gauge theory, and the second
main result in this work is a characterisation of some if its basic
properties.

Another recurrent theme on the study of instantons on Euclidean space
is the equivalence with certain algebraic curves. They appear as
{\em jumping lines} in the original ADHM construction and as
{\em spectral curves} in Hitchin's construction of monopoles and on
the study of instantons invariant under $\real^2$. One might then
expect that some suitable algebraic curves will also play a
significant role. This turns out to be indeed the case, as we shall
see in theorem \ref{specthm}. Again, useful information about the
instanton moduli space is gained from this point of view.


\section{Instantons and Hitchin's equations}

Before we state the results to be proven in this thesis, it is
convenient to gather some relevant definitions here. More precisely,
we set up our configuration space of connections on a vector bundle
$E\seta\torus$ in order to make clear what we mean by an {\em instanton}.
Due to the non-compactness of our base manifold $\torus$, this really
requires some extra work. We then proceed to briefly recall the
definition of the Hitchin's equations over an elliptic curve.

\paragraph{On the choice of metric and complex structure.}
The surface we want to consider has at least three reasonable models:
$$ \torus \simeq T\times(\proj\setminus\{\infty\}) \simeq T\times S^1\times[0,\infty] $$
which we respectively call the {\em plane}, {\em round} and {\em
cylindrical} models. Of course, these surfaces are all
diffeomorphic, but each one has its own natural choice of a
riemannian metric, namely the product one.

Moreover, the respective product metrics are not conformal to one another.
This leads to three different concepts of anti-self-duality and finite energy,
so that instantons in one model are not instantons on the others.

There is a good amount of literature studying the round and cylindrical
cases (see \cite{B2} and \cite{K}, respectively). In this work, however,
we are interested only on the plane model, since we want to think of
$\torus$ as the quotient of $\real^4$ by a two-dimensional lattice
$\zed^2$. Hence, $\torus$ will always be equipped with its product
riemannian metric; a complex structure $I$ coming from the product of
a complex structure on the torus with a complex structure on the complex
line is assumed to be fixed and we denote by $\kappa$ the associated
K\"ahler form. Moreover, the compactified version $\tproj$ will always
be equipped with its product riemannian metric and a complex
structure compatible with $I$ is chosen; we denote the associated K\"ahler form by
$\overline{\kappa}$.

Actually, note that $\torus$ inherits a hyperk\"ahler structure from
$\real^4$; the two other complex structures arise if we regard $\torus$
as the product of two cylinders $(S^1\times\real)\times(S^1\times\real)$.

On the other hand, we also want to think of the dual torus as a
quotient of $(\real^4)^*$ by the dual group of translations
$\Lambda^*$. Thus, $\dual$ is given the flat, Euclidean metric.
Moreover, the choice of a complex structure of $\torus$ also fixes a
complex structure on $\dual$, since this is seen as a lattice in
$(\real^4)^*$.


\subsection{Instantons over $\torus$.}
An {\em instanton} is a smooth, anti-self-dual connection $A$
on an $SU(2)$ bundle $E\seta\torus$ with a system of transitions
functions lying in $L^2_3({\rm Aut}E)$. As we mentioned above,
anti-self-duality is taken with respect to the product metric
$\kappa$ on the base.

Alternatively, $\torus$ can be thought as a quotient of $\real^4$
by a two-dimensional lattice $\zed^2$. In this way, $A$ is regarded
as a $SU(2)$ connection on a bundle over $\real^4$ which is
invariant under the action of $\zed^2$ by translations, i.e. $A$ is
periodic in two directions of the 4-plane, fitting therefore in the
framework described at the introduction.

Given a function $f:\cpx\seta\real$, we say that $f\sim O(|w|^n)$ if:
$$ \lim_{w\seta\infty} \frac{|f(w)|}{|w|^n} < \infty $$

In this work, to avoid deeper analytical problems, we will consider only
anti-self-dual connections $A$ on $E\seta\torus$ satisfying the following
conditions:
\begin{enumerate}
\item $|F_A|\sim O(r^{-2})$;
\item there is a holomorphic vector bundle $\ee\seta\tproj$ with trivial
determinant such that $\ee|_{T\times(\proj\setminus\{\infty\})}\simeq(E,\del_A)$,
where $\del_A$ is the holomorphic structure on $E$ induced by the instanton
connection $A$;
\end{enumerate}
Such connections are said to be {\em extensible}. Moreover, we assume
the restriction of the extended bundle to the added divisor splits as a sum
of flat line bundles, i.e.:
$$ \ee|_{T_\infty}=L_{\ksi_0}\oplus L_{-\ksi_0} $$
and $\pm\ksi_0$ can be seen as points in the Jacobian torus
$\jj(T)=\dual$. We say $\ksi_0$ is the {\em asymptotic state}
of the connection $A$. We also fix the topological type of the extended
bundle $\ee$ by making $c_2(\ee)=k>0$; the integer $k$ is the {\em instanton
number} of the connection $A$.

Finally, we also assume that $A$ is {\em irreducible} as an $SU(2)$ connection.
In particular, this implies that $E$ admits no square-integrable covariantly
constant sections, i.e.:
\begin{equation} \label{irred}
||\nabla_As||_{L^2}>0
\end{equation}
for all $s\in L^2(E)$ not constant.

\paragraph{Spectral curve.}
The holomorphic extension of $E\seta\torus$ to $\ee\seta\tproj$ we
mentioned above leads us to look at a construction due to Friedman,
Morgan \& Witten \cite{FMW}. These authors have shown how one can
associate a pair of {\em spectral data}, consisting of a complex curve
$S$ plus a line bundle $\call\seta S$, to holomorphic vector bundles
over elliptic surfaces. We shall pursue this point of view in section
\ref{instspec}.


\subsection{Hitchin's equations.}

If $\Lambda=\zed^2$ then $\Lambda^*=\zed^2\times\real^2$. According
to the scheme outlined in the introduction, we must look at ASD
connections on a suitable $(\real^4)^*$ which do not depend on two
coordinates and are periodic on the other two. These objects were
studied by Hitchin \cite{H} and correspond to solutions of the
so-called {\em Hitchin's equations} over the two-dimensional torus
$\dual=(\real^4)^*/\Lambda^*$; these can be obtained via
dimensional reduction of the usual ASD equations from four to two
dimensions.

More precisely, let $V\seta\real^4$ be a rank $k$ vector bundle with a
connection $\tilde{B}$ which does not depend on two coordinates. Pick
up a global trivialisation of $V$ and write down $\tilde{B}$ as
a 1-form:
$$ \tilde{B}=B_1(x,y)dx+B_2(x,y)dy+\phi_1(x,y)dz+\phi_2(x,y)dw $$
Hitchin then defined a {\em Higgs field} $\Phi=(\phi_1+i\phi_2)d\xi$,
where $d\xi=dx+idy$. So $\Phi$ is a section of $\Lambda^{1,0}{\rm End}V$,
where $V$ is now seen as a bundle over $\real^2$ with a connection
$B=B_1dx+B_2dy$.

The ASD equations for $\tilde{B}$ over $\real^4$ can then be
rewritten as a pair of equations on $(B,\Phi)$ over $\real^2$:
\begin{equation} \label{hiteqs} \left\{ \begin{array}{l}
F_B+[\Phi,\Phi^*]=0 \\
\del_B\Phi=0
\end{array} \right. \end{equation}
These equations are also conformally invariant, so they depend only on
the conformal class of the Euclidean metric on $\dual$. Solutions
$(B,\Phi)$ are often called {\em Higgs pairs}.

As we mentioned above, the Nahm transform will produce singular
solutions of (\ref{hiteqs}); in fact, there are very few smooth
solutions for bundles over elliptic curves (see \cite{H}). The
particular class of singular solutions that will appear was
studied by several authors \cite{S} \cite{Ko} \cite{K} \cite{B3},
and are related to the parabolic vector bundles of Mehta \& Seshadri
\cite{MS}. The presence of singularities in the dual object
is not at all surprising. In fact, we shall see that they encode the
asymptotic behaviour of the original connections over $\torus$, just
as in the case of monopoles \cite{H3}.

Therefore, we will study solutions of (\ref{hiteqs}) over $\dual$ with
the singularities removed. The Euclidean metric becomes incomplete,
and one cannot expect to have a finite dimensional moduli space of
solutions. However, since the equations depend only on the conformal
structure, we are allowed to perform conformal changes in the
metric. Indeed, we will follow Biquard \cite{B3} and consider the
so-called Poincar\'e metric (which is complete) when we study the
relevant singular Higgs pairs in section \ref{inv}.

We have one last important hypothesis. A Higgs pair $(B,\Phi)$ is
said to be {\em admissible} if the bundle $V$ has no covariantly
constant sections, i.e.:
$$ ||\nabla_Bs||_{L^2}>0 $$
for all $s\in L^2(V)$ not constant.

\paragraph{Spectral curves.}
In \cite{H2}, Hitchin has shown that smooth solutions of
(\ref{hiteqs}) are equivalent to a set of {\em spectral data},
consisting of a complex curve $C$ plus a line bundle ${\cal N}\seta
C$. This was later generalised to singular solutions by various
people. We will review this construction more carefully in section
\ref{hitspec}.


\section{Statement of the main results}

We are now in position to state the first main result to be proven
here. It provides a correspondence between finite energy ASD
connections over $\torus$ and singular solutions of Hitchin's
equations over the punctured dual torus, where the Higgs field is
allowed to have simple poles with a definite residue. More
precisely, we have:

\begin{mthm} \label{nahmthm}
The Nahm transform is a bijective correspondence between the following objects:
\begin{itemize}
\item gauge equivalence classes of irreducible, extensible $SU(2)$
instanton connections on $E\seta\torus$ with fixed instanton number
$k$ and asymptotic state $\ksi_0$; and
\item admissible $U(k)$ solutions of the Hitchin's equations over
the dual torus $\dual$, such that the Higgs field has at most
simple poles at $\as\in\dual$; moreover, its residues are semi-simple
and have rank $\leq2$, if $\ksi_0$ is an element of order 2 in the
Jacobian of $T$, and rank $\leq1$ otherwise.
\end{itemize} \end{mthm}

It is interesting to note here that the behaviour of the Higgs
field $\Phi$ near the singularities $\as$ is determined by the
behaviour at infinity of the original instanton, and vice-versa.
This is analogous to what happens in the monopole case \cite{H3}.

The proof of this theorem will be carried out in chapters
\ref{nahm} to \ref{conc}. There are two possible approaches: the
gauge-theoretical construction of sections \ref{diffgeo} and \ref{inv}
and the purely holomorphic approach of sections \ref{holo} and
\ref{invholo}. These actually complement each other, and the whole
proof uses a mixture of both.

The above result has a physical interpretation in terms of certain
supersymmetric theories, given by Kapustin \& Sethi \cite{KS}
\cite{Ka}. The four dimensional theory containing the instanton is
regarded as the low energy limit of a type IIA string theory
containing NS5- and D4-branes wrapped around a torus $T$. A version
of mirror symmetry (T-duality) plays the role of the Nahm's transform,
mapping this theory to another one containing only D-branes wrapped
around the dual torus $\dual$. Simultaneously, the Coulomb branch of
the original 4-dimensional theory (i.e. the moduli space of doubly-periodic
instantons) is mapped onto the Higgs branch of a 5-dimensional
impurity theory (i.e. the moduli space of Higgs pairs on $\dual$).

In appendix B we will indicate how the above result could be
modified to assume a more general condition on the instanton
connection. More precisely, one can expect to exchange the
extensibility hypothesis for a pointwise estimate for the
asymptotic behaviour of the curvature $F_A$.

\smallskip

In chapter \ref{spec} we turn to the study of the spectral
curves associated to each side of the Nahm transform. We start
by reviewing the construction of the spectral data associated
to holomorphic vector bundles over elliptic surfaces \cite{FMW}
and to singular Higgs pairs \cite{H2}. After establishing various
facts about them, we show that:

\begin{mthm} \label{specthm}
If $(V,B,\Phi)$ is the holomorphic Nahm transform of $(E,A)$, then
the instanton spectral data $(S,\call)$ associated to $(E,A)$
coincide with the Higgs spectral data $(C,{\cal N})$ associated to
$(V,B,\Phi)$, in the sense that the curves $S$ and $C$ coincide
pointwise and there is a natural line bundle isomorphism
$\call\seta{\cal N}$.
\end{mthm}

One of the consequences this last result is a nice picture of the
moduli space of doubly-periodic instantons: it has the structure of
a fibration over the space of spectral curves (of complex dimension
$2k+1$), with fibres given by the Jacobian of the given curve (of
genus $2k-1$). Thus, we conclude that the moduli space of
extensible doubly-periodic instanton connections is a smooth,
complex manifold of dimension $4k$. Moreover, the latter and the
moduli space of singular Higgs pair are explicitly seen to be
diffeomorphic, with the Nahm transform as a diffeomorphism.

Finally, theorem \ref{specthm} closes a circle of ideas analogous
to the one considered by Hitchin in the case of monopoles
\cite{H3}, giving a correspondence between doubly-periodic
instantons, the Nahm transformed singular Higgs pair and the
associated spectral data.


\chapter{Analytical background} \label{back}

The first stage towards the proof of our main theorem is to
sort out a few analytical problems caused by the non-compactness
of $\torus$. Clearly, the extensibility hypothesis saves us
some hard work (see however appendix B). In this chapter we
will look at the Dirac operator coupled to an extensible
connection, proving that it is Fredholm in section \ref{fred}.

First, let us recall the conditions for extensibility; an
instanton connection $A$ is {\em extensible} if it satisfies:
\begin{enumerate}
\item $|F_A|\sim O(r^{-2})$;
\item there is a holomorphic rank two vector bundle $\ee\seta\tproj$ with trivial
determinant such that $\ee|_{T\times(\proj\setminus\{\infty\})}\simeq(E,\del_A)$,
where $\del_A$ is the holomorphic structure on $E$ induced by the instanton
connection $A$.
\end{enumerate}


\section{Instanton number and asymptotic states}

We now use the extensibility hypothesis to study the compatibility
between the instanton connection $A$ and the extended bundle
$\ee\seta\tproj$. More precisely, we first want to show that the
holomorphic type of the restriction of the extended bundle to the
added divisor  $T_\infty=T\times\{\infty\}$ is indeed directly
determined by the asymptotic behaviour of the instanton connection
$A$. Then we argue that the topology of $\ee$ is fixed by the
energy ($L^2$-norm) of $A$.

Before that, we must fix an appropriate trivialisation at infinity.


\subsection{Good gauge at infinity}

Let $B_R$ denote a closed ball in $\cpx$ of radius R, and
let $V_R$ be its complement. Also, consider the obvious
projection $p:T\times V_R\seta T$. We shall need the
following technical result, whose proof we postpone to appendix B.

\begin{prop} \label{goodgauge}
If $|F_A|\sim O(r^{-2})$, then, for $R$ sufficiently large, there
is a gauge over $T\times V_R$ and a constant flat connection $\Gamma$
on a topologically trivial rank two bundle over the elliptic curve
such that:
$$ A-p^*\Gamma \ = \ \alpha \ \sim O(r^{-1}\cdot\log r) $$
\end{prop}

\paragraph{Asymptotic states.}
By general theory, a constant flat connection
on a bundle $S\seta T$ determines uniquely a holomorphic structure on this
bundle. Moreover, $S$ must split, holomorphically, as the sum of two line
bundles, i.e. $S=L_{\ksi_0}\oplus L_{-\ksi_0}$, uniquely up to $\pm1$. Here,
$\as$ are seen as points in $\dual$, the Jacobian of the elliptic curve $T$.

Therefore, by proposition \ref{goodgauge}, to each extensible instanton
connection we can associate an unique pair of opposite points $\as\in\dual$.
Such points are called the {\em asymptotic states} of $A$.

\begin{lem}
If an extensible instanton connection $A$ has asymptotic states $\as$, then
$\ee|_{T_\infty}=L_{\ksi_0}\oplus L_{-\ksi_0}$.
\end{lem}

\pf
Let $V_\infty\subset\proj$ be a small neighbourhood centred at
$\infty\in\proj$; let $w$ be a coordinate there. We can regard
$\ee|_{T\times V_\infty}$ as a family of rank 2 bundles over $T$,
parametrised by $w$, Furthermore,
If $\del$ denotes the holomorphic structure on $\ee$, let
$\del_w$ be the holomorphic structure on the restriction $\ee|_{T_w}$.
Clearly, as operators:
$$ \lim_{w\seta\infty}\del_w=\del_\infty $$
However, from condition (2) in the definition of extensibility,
we know that $\del_w=\del_{A|_{T_w}}$ away from $\infty$. But
proposition \ref{goodgauge} tells us that $\del_{A|_{T_w}}$
approaches $\del_{\Gamma}$ as $w\seta\infty$. Therefore,
$\del_\infty=\del_{\Gamma}$, and the lemma follows. \pfend


\subsection{The instanton number}
Moreover, as we mentioned before, the topological type of $\ee$
is determined by the energy of the instanton connection:

\begin{lem} \label{c2=i}
$c_2(\ee)=\frac{1}{8\pi^2}\int_{\torus}|F_A|^2$
\end{lem}

\pf Again, let $V$ be a small neighbourhood of $\infty\in\proj$. Let
$\Gamma_{\as}$ be the canonical connection on the bundle
$L_{\ksi_0}\oplus L_{-\ksi_0}$ over an elliptic curve and consider the
projection $p:T\times V\seta T$.

Now consider a connection $A'$ on the extended bundle $\ee$ that coincides
with $p^*\Gamma_{\as}$ on $T\times V$. Therefore
\begin{eqnarray} \label{c2}
c_2(\ee) & = & \frac{1}{8\pi^2} \int_{\tproj} {\rm Tr}(F_{A'}\wedge F_{A'}) \ = \
\frac{1}{8\pi^2} \int_{T\times(\proj\setminus\{\infty\})} {\rm Tr}(F_{A'}\wedge F_{A'}) \nonumber \\
& = & \frac{1}{8\pi^2} \lim_{R\seta\infty}\int_{T\times B_R} {\rm Tr}(F_{A'}\wedge F_{A'})
\end{eqnarray}

On the other hand, we have from lemma \ref{goodgauge} that $A-A'=\alpha$ is a
1-form in $O(r^{-1}\cdot\log(r))$. Define the 1-parameter family of connections
$A_t=A'+t\cdot\alpha$, so that the corresponding curvatures:
\begin{eqnarray}
& F_{A_t}=t\cdot F_{A} + (1-t)\cdot F_{A'} -
\left( t-\frac{t^2}{2}\right)\cdot\alpha\wedge\alpha & \nonumber \\
& \Longrightarrow\ \ \ |F_{A_t}|\sim O(r^{-2}\cdot \log^2r)\ \ 
\forall t\in[0,1] \label{famcurvs}
\end{eqnarray}

So let:
\begin{equation} \label{i}
i(A) = \frac{1}{8\pi^2} \int_{\torus} {\rm Tr}(F_{A}\wedge F_{A}) \ = \
\frac{1}{8\pi^2} \lim_{R\seta\infty}\int_{T\times B_R} {\rm Tr}(F_{A}\wedge F_{A})
\end{equation}
Usual Chern-Weil theory tells us that:
\begin{eqnarray*}
c_2(\ee)-i(A) & = & \frac{1}{8\pi^2} \lim_{R\seta\infty} \left\{ \int_{T\times B_R} \left(
{\rm Tr}(F_{A'}\wedge F_{A'})-{\rm Tr}(F_{A'}\wedge F_{A'}) \right) \right\} \ = \\
& = & \frac{1}{4\pi^2} \lim_{R\seta\infty} \left\{ \int_{T\times B_R} d \left(
\int_0^1 {\rm Tr}(\alpha\wedge F_{A_t}) \right) \right\} \ = \\
& = & \frac{1}{4\pi^2} \lim_{R\seta\infty} \left\{ \int_{T\times S_R^1} \left(
\int_0^1 {\rm Tr}(\alpha\wedge F_{A_t}) \right) \right\} \ = 0
\end{eqnarray*}
by our estimates in proposition \ref{goodgauge} and in equation
({\ref{famcurvs}). This completes the proof. \pfend

In particular, the integral in the right hand side of the equation in
lemma \ref{c2=i} has to equal an integer number $k>0$, which we call
the {\em instanton number} of $A$.

Finally, we say that an extensible connection $A$ on the bundle
$E\seta\torus$ belongs to $\cala_{(k,\ksi_0)}$ if it has instanton
number $k$ and asymptotic state $\ksi_0$.


\subsection{Estimating the Dolbeault operator}
Finally, we need a final lemma that will be useful in the following
section section, where we develop a Fredholm theory for the Dirac
operator coupled to an instanton connection $A\in\cfg$.

First, note that the bundle $L_{\ksi_0}\oplus L_{-\ksi_0}\seta T$ admits a
flat connection with constant coefficients, which we denote by
$\Gamma_{\ksi_0}$. Use the projection $T\times V_R\stackrel{p_1}{\seta}T$ to
pull it back to $T\times V_R$. We show that:

\begin{lem} \label{model}
Let $A\in\cfg$ be any extensible instanton connection. Given
$\epsilon>0$, there is $R$ sufficiently large such that:
$$ ||\del_A-\del_{\Gamma_{\ksi_0}}||_{L^2(T\times V_R)}<\epsilon $$
\end{lem}

\pf Since $\del_A-\del_{\Gamma_{\ksi_0}}$ is just the $(0,1)$-part
of the 1-form $\alpha=A-\Gamma_{\ksi_0}$, the statement is a simple
consequence of the gauge-fixing proposition \ref{goodgauge}. \pfend






\section{The Poincar\'e line bundle} \label{poin}

We now quickly review some facts regarding holomorphic vector bundles
over elliptic curves and surfaces that will be useful later on. We are
particularly interested in the definition of the Poincar\'e line bundle
and on Atiyah's classification result \cite{A3}.

Recall that an elliptic curve is a two-dimensional torus $T$ with a
complex structure, plus the choice of a point $e\in T$ which plays
the role of the identity element of the torus as an abelian group.
For simplicity, we denote an elliptic curve only by $T$, letting
the choice of the identity element always implicit.

The {\em Jacobian} $\jj(T)=\hat{T}$ of $T$ is defined as the set of
flat holomorphic line bundles over $T$. Such bundles can be
parametrised by $T$ itself in the following way: to each $z\in T$,
we associate the bundle $\call_z=\oo_T(e)\otimes\oo_T(z)^{-1}$.
Hence $T$ and $\dual$ are isomorphic as elliptic curves, and the
identity element $\hat{e}\in\dual$ corresponds to the
holomorphically trivial line bundle $\underline{\cpx}\seta T$.
Moreover, the set of flat holomorphic line bundles over $\dual$ is
again $T$. Throughout the thesis, points in $T$ are denoted by $z$
and points in $\dual$ are denoted by $\ksi$.

An element $\ksi$ of $\dual$ has order 2 if $L_\ksi\otimes L_\ksi=\underline{\cpx}$.
The are four such elements, one of them being the identity $\hat{e}$.

Moreover, the line bundles $L_\ksi\seta T$ and $L_z\seta\dual$ can be
given a natural constant connection compatible with the holomorphic
structure. This follows from the differential-geometric definition of $\dual$:
$$ \dual=\{\ksi\in(\real^4)^*\ |\ \ksi(z)\in\zed,\forall z\in\Lambda^2 \} $$
where $\Lambda^2\subset\real^4$ is the two-dimensional lattice
generating $\torus$. Hence each $\ksi\in\dual$ can be regarded as a
constant, real 1-form over T, so that $\omega_\ksi=i\ksi$ is a connection on a
topologically trivial line bundle $L\seta T$. Each such connection
defines a different holomorphic structure on $L$, which we denote
by $L_\ksi$. The holomorphic line bundles $L_z\seta\dual$ are
defined on the same way.

Note that, in the notation of lemma \ref{model},
$\Gamma_{\ksi_0}=\omega_{\ksi_0}\oplus\omega_{-\ksi_0}$.

\paragraph{The Poincar\'e bundle.}
The {\em Poincar\'e line bundle} {\bf P}$\seta T\times\dual$ is the
unique holomorphic line bundle satisfying:
$$ \begin{array}{ccc} {\mathbf P}|_{T\times\{\ksi\}}\simeq L_\ksi & \ \ &
{\mathbf P}|_{\{z\}\times\dual}\simeq L_{-z} \end{array} $$
It can be constructed as follows. Identifying $T$ and $\dual$ as before,
let $\Delta$ be the diagonal inside $T\times\dual$, and consider the
divisor $D=\Delta-T\times{\hat{e}}-{e}\times\dual$. Then
${\mathbf P}=\oo_{T\times\dual}(D)$; it is easy to see that the sheaf so
defined restricts as wanted.

Note that although the two restrictions above are flat line bundles over $T$ and
$\dual$ respectively, the Poincar\'e bundle itself is not topologically trivial;
in fact, $c_1({\mathbf P})\in H^1(T)\otimes H^1(\dual)\subset H^2(T\times\dual)$.
More precisely, the unitary connection and its corresponding curvature are given by:
$$ \begin{array}{ccc}
\omega(z,\ksi)=i\sum_{\mu=1}^{2}\ksi_\mu dz_\mu-z_\mu d\ksi_\mu & \ \ &
\Omega(z,\ksi)=i\sum_{\mu=1}^{2}d\ksi_\mu\wedge dz_\mu
\end{array} $$
Restricted to $T\times\{\ksi\}$, these give the bundles $L_\ksi\seta T$
flat connections $\omega_\ksi=i\sum_{\mu=1}^{2}\ksi_\mu dz_\mu$, with
constant coefficients. Similarly, the bundles $L_z\seta\dual$ also have
canonical flat connections $\omega_z=-i\sum_{\mu=1}^{2}z_\mu d\ksi_\mu$.

Finally, note that $c_1({\bf P})^2=2\cdot t\wedge\hat{t}$, where $t$ and
$\hat{t}$ are the generators of $H^2(T)$ and $H^2(\dual)$,
respectively.

\paragraph{Atiyah's classification result.}
Holomorphic vector bundles $\vv\seta T$ are classified by the following
result due to Atiyah \cite{A3}. The building blocks for Atiyah's classification are the
holomorphic vector bundles constructed as follows. Start by defining
${\mathbf F}_1=\underline{\cpx}$; then ${\mathbf F}_n$ is defined recursively as the
unique non-trivial extension of ${\mathbf F}_{n-1}$ by $\underline{\cpx}$:
$$ 0\seta\underline{\cpx}\seta {\mathbf F}_{n} \seta {\mathbf F}_{n-1} \seta0 $$

\begin{thm} \label{atiyah}
Let $\vv\seta T$ be an indecomposable rank $r$ holomorphic vector bundle such that
${\rm deg}\vv=0$. Then $\vv={\mathbf F}_r\otimes L$, for some flat line bundle
$L\seta T$.
\end{thm}

In particular, for the case of rank $2$ bundles we have:

\begin{thm} \label{ranktwo}
Let $\vv\seta T$ be a semi-stable rank 2 holomorphic vector bundle such that
${\rm det}\vv=\underline{\cpx}$. Then either of two possibilities must hold:
\begin{itemize}
\item $\vv$ is decomposable, and $\vv=L\oplus L^{-1}$, where
$L\in\dual$ is uniquely determined up to $\pm1$;
\item $\vv$ is indecomposable, and $\vv={\mathbf F}_2\otimes L$, where
$L\in\dual$ is an uniquely determined element of order $2$.
\end{itemize} \end{thm}

Note that semi-stability excludes only decomposable bundles looking like
$Q\oplus Q^{-1}$, where $Q\seta T$ has degree $n>0$. Moreover,
semi-stability is a Zariski open condition.


\paragraph{Elliptic surfaces.}

Recall that an elliptic surface is a complex surface $S$ with a map $\pi:S\seta B$
to a compact curve $B$ such that $\pi^{-1}(b)$ is an elliptic curve for generic
$b\in B$; there might be points $b\in B$ such that $\pi^{-1}(b)$ is singular or
multiple. This is a vast class of complex surfaces and there is a large theory
about them, but we are interested here only in a quite simple case: $S=\tproj$ and
$\pi$ the usual projection onto the second factor (hence $B=\proj$).

The Jacobian surface ${\cal J}(S)$ of $S$ is defined to be the elliptic surface
obtained, roughly speaking, in the following manner. For each $b\in B$, we replace
the elliptic curve $T=\pi^{-1}(b)$ by its Jacobian curve, so that they fit together
to form a new elliptic surface. In our case of interest, ${\cal J}(S)=\dual\times\proj$.

It is also possible to define a Poincar\'e bundle ${\bf P}_{S}$ over an elliptic
surface. For the case we are interested in, ${\bf P}_{S}=p_{13}^*{\bf P}$,
where $p_{13}:\tproj\times\dual\seta T\times\dual$ is the natural
projection on to the first and third factors. For the most general
definition, see \cite{FMW}, p. 688.

\section{Fredholm theory of the Dirac operator} \label{fred} 

We begin by recalling that the dual torus $\dual$ parametrises the set of flat
holomorphic line bundles $L\seta T$. Moreover, such bundles have a
natural choice of connection, denoted $i\ksi$, which is consistent
with the holomorphic structure.

In fact, $\dual$ also parametrises the set of flat holomorphic line
bundles over $\torus$. Using the projection $p_1:\torus\seta T$,
one obtains the holomorphic line bundle $p_1^*(L_{\ksi})$ over
$\torus$, which we shall also denote by $L_{\ksi}$ for simplicity;
let $\omega_{\ksi}$ be the pullback of the flat constant connection
on $L_{\ksi}\seta T$ described above; clearly, such connection is
also flat.

As usual, let $E\seta\torus$ be a rank 2 bundle provided with an
instanton connection $A\in\cala_{(k,\ksi_0)}$. Form the bundle
$E\otimes L_{\ksi}$ with the corresponding connection
$A_{\ksi}=A\otimes I + I\otimes\omega_{\ksi}$; since all we have done
was to add a flat term to our original instanton, $A_{\ksi}$ is
still an instanton on the twisted bundle. We also require $A$ to be
irreducible; clearly, its twisted version $A_{\ksi}$ is also
irreducible.

Consider now the Dirac operator acting on the bundle $E(\ksi)=E\otimes L_{\ksi}$,
coupled to the connection $A_{\ksi}$, and its adjoint:
\eq \left\{ \begin{array}{c}
D_{A_{\ksi}}:\Gamma(E(\ksi)\otimes S^+)\seta\Gamma(E(\ksi)\otimes S^-) \\
D_{A_{\ksi}}^*:\Gamma(E(\ksi)\otimes S^-)\seta\Gamma(E(\ksi)\otimes S^+)
\end{array} \right. \end{equation}
where the spaces of sections are provided with norms suitably defined.
Since the base manifold is flat and the connection is anti-self-dual,
the Weitzenb\"ock formula on $E(\ksi)\otimes S^+\seta\torus$ is simply:
\begin{eqnarray}
D_{A_{\ksi}}^* D_{A_{\ksi}} & = & \nabla_{A_{\ksi}}^*\nabla_{A_{\ksi}}  \label{weit} \\
\Rightarrow \ \ ||D_{A_{\ksi}}s||^2 & = & ||\nabla_{A_{\ksi}}s||^2 \nonumber
\end{eqnarray}
Hence, if $A_{\ksi}$ is irreducible, there are no covariantly constant
sections of \linebreak $E(\ksi)\otimes S^+$. This means that the
kernel of  $D_{A_{\ksi}}$ is trivial. Now, if $D_{A_{\ksi}}$ is a
Fredholm operator, then  ${\rm ker}D_{A_{\ksi}}^*$ (which coincides
with ${\rm coker}D_{A_{\ksi}}$)  is a finite dimensional subspace of
$\Gamma(E(\ksi)\otimes S^-)$.

In this rather technical but fundamental section, we prove that this is indeed the case:

\begin{thm} \label{dirac-fred}
Given any extensible instanton connection $A\in\cfg$, the Dirac operators:
\begin{equation} \label{oopp}
D_{A_{\ksi}}^*:L^2_1(E(\ksi)\otimes S^-)\seta L^2(E(\ksi)\otimes S^+)
\end{equation}
form a smooth family of Fredholm operators parametrised by $\dual\setminus\{\as\}$.
Moreover, $\ind D_{A_\ksi}^*=k$, for all $\ksi\in \dual\setminus\{\as\}$.
\end{thm}

\smallskip

The Sobolev norm in the left hand side of (\ref{oopp}) is defined
as follows. Let $D_\ksi^*$ be the Dirac operator $L_\ksi\otimes
S^-\seta L_\ksi\otimes S^+$. Then $L^2_1(E(\ksi)\otimes S^-)$ is
the completion of $\Gamma(E(\ksi)\otimes S^-)$ with respect to the
norm:
\begin{equation} \label{n.thm}
||s||_{L^2_1}=||s||_{L^2}+||D_\ksi^*s||_{L^2}
\end{equation}

The proof consists of three steps, which we now outline. We first
prove that the operator $D_\ksi^*:L^2_1(L_\ksi\otimes S^-)\seta L^2(L_\ksi\otimes S^+)$
is invertible for nontrivial $\ksi\in\dual$. A gluing argument then
shows that the Dirac operator coupled to a twisted instanton $A_\ksi$ is
Fredholm if $\ksi\neq\ksi_0$, after using the fact that the set of
Fredholm operators is open. To compute the index, we use an argument
based on the Gromov-Lawson Relative Index Theorem \cite{GL}; the
details are left to the appendix.


\paragraph{The flat model.}
Let $L_\ksi\seta\torus$ be the flat line bundle described above, provided
with the constant connection $\omega_\ksi$. Our starting point to prove
theorem \ref{dirac-fred} is the following proposition.

\begin{prop} \label{flatmodel}
For non-trivial $\ksi\in\dual$, the coupled Dirac operator
$$ D_\ksi^*:L^2_1(L_\ksi\otimes S^-)\seta L^2(L_\ksi\otimes S^+) $$
is invertible. Its inverse is denoted by $Q^\infty_\ksi$.
\end{prop}

\pf
Let $L_{\ksi}\seta\torus$ be the pull-back via $p_1:\torus\seta T$ of
a flat line bundle over the 2-torus, provided with the constant
connection $\omega_{\ksi}=p^*(-i\ksi)$, as described in section
\ref{poin}. Consider the twisted Dirac operator:
$$ D_{\ksi}:\Gamma(L_{\ksi}\otimes S^+)\seta\Gamma(L_{\ksi}\otimes S^-) $$
and its adjoint $D_{\ksi}^*$.

Since $M=\torus$ is a K\"ahler surface, we have the
following decompositions:
\eq \label{decomp} \left\{ \begin{array}{l}
S^+=\Lambda^{(0,0)}_M L_{\ksi} \oplus \Lambda^{(0,2)}_M L_{\ksi} \\
S^-=\Lambda^{(0,1)}_M L_{\ksi} =\Lambda^{(0,1)}_{T} L_{\ksi} \oplus
    \Lambda^{(0,1)}_{\cpx}
\end{array} \right. \end{equation}
With respect to these decompositions, the Dirac operator and its adjoint are given by:
\eq \label{dirac} \begin{array}{cc}
D_{\ksi}=\matriz \del^{(z)}_{\ksi} & -\del^{(w),*}_{\ksi} \\
\del^{(w)}_{\ksi} & -\del^{(z),*}_{\ksi}  \matfim &
D^*_{\ksi}=\matriz -\del^{(z),*}_{\ksi} & -\del^{(w),*}_{\ksi} \\
\del^{(w)}_{\ksi} & \del^{(z)}_{\ksi}  \matfim
\end{array} \end{equation}
where $\del^{(z,w)}_{\ksi}$ denotes the Dolbeault operator twisted by $\omega_{\ksi}$ along the
toroidal ($z$) and plane ($w$) complex coordinates, i.e. the components of the covariant
derivative. Hence, the coupled Dirac laplacian $\triangle_{\ksi}=D_{\ksi}^*D_{\ksi}$ mapping
$\Lambda^{(0,0)}_M L_{\ksi}\oplus \Lambda^{(0,2)}_M L_{\ksi}$ to itself is just:
\eq \label{lap} \matriz
\triangle_{\ksi}^{(z)}+\triangle_{\ksi}^{(w)} & 0 \\
0 & \triangle_{\ksi}^{(z)}+\triangle_{\ksi}^{(w)}
\matfim \end{equation}
The off-diagonal terms are cancelled, for they are proportional to
the curvature, which was supposed to vanish. Moreover, the flat
connection $\omega_\ksi$ is a pull back from the torus, so that
$\triangle^{(w)}_{\ksi}$ is just the usual plane laplacian $\triangle^{(w)}$.
Let us concentrate on a single component, say $\Lambda^{(0,0)}_M L_{\ksi}$.

First, we want to solve the homogeneous equation $\triangle_{\ksi}f=0$ for \linebreak
$f\in\Lambda^{(0,0)}_M(L_{\ksi})$ and a fixed $\ksi\in\dual$.
Now, separate variables, supposing that $f(z,w)=\varphi(z)g(w)$:
$$ \triangle_{\ksi}f=0\ \see \ g\triangle_{\ksi}^{(z)}\varphi+\varphi\triangle^{(w)}g=0 $$
Therefore:
\eq \label{sep} \left\{ \begin{array}{l}
(i)\ \triangle_{\ksi}^{(z)}\varphi=\lambda^2\varphi \\
(ii)\ \triangle^{(w)}g=-\lambda^2 g\ \rightarrow\  (\triangle^{(w)}+\lambda^2)g=0
\end{array} \right. \end{equation}
where $\lambda^2$ are the eigenvalues of the $\ksi$-twisted laplacian over
the torus. They form a discrete, unbounded set $\{\lambda_n(\ksi)\}$
of $\real^+$, each being a function of the parameter $\ksi$. Note that
since $H^0(T,L_\ksi)=0$ for nontrivial $\ksi\in\dual$, we can indeed
guarantee that $\lambda_n(\ksi)>0$ for all nontrivial $\ksi$. On the
other hand, for $L_\ksi=\underline{\cpx}$, the laplacian has a
1-dimensional kernel, i.e. one zero eigenvalue.

As usual, we can decompose $f$ on the eigenstates of $\triangle_{\ksi}^{(z)}$,
i.e.:
\begin{equation} \label{f.decomp}
f=\sum_{n}g_n(w)\varphi_{n}(z)
\end{equation}
where $\{\varphi_n\}$ is an orthonormal basis for the $L^2$ norm on
$\Lambda^{(0,0)}_M(L_{\ksi})$ of eigenstates with eigenvalues
$\{\lambda^2_n\}$; so, $||f||_{L^2(\torus)}^2=\sum_{n}||g_n||^2_{L^2(\cpx)}$.
Moreover:
\eq \label{evalue}
\triangle_{\ksi}f=\sum_{n}[(\triangle^{(w)}+\lambda_n^2)g_n]\varphi_n
\end{equation}

\begin{prop} \label{est1}
Let $\rho\in L^2(L_{\ksi}\otimes S^+)$ be compactly supported and suppose
that $\ksi$ is nontrivial. Then there is $f\in L^2(L_{\ksi}\otimes S^+)$
and a constant $k<\infty$ such that $\Delta_{\ksi}f=\rho$ and
$||f||_{L^2}\leq k||\rho||_{L^2}$.
\end{prop}

\pf Given (\ref{evalue}), solving the equation $\triangle_{\ksi}f=\rho$ amounts to solve
\linebreak $(\triangle^{(w)}+\lambda^2_n)g_n=\rho_n$ for each $n$, where $g_n,\rho_n$ are the
components of $g,\rho$ along the eigenspaces of $\lambda^2_n$, respectively.

Fix some integer $n$ and denote by $F_n$ the fundamental solution of
$(\Delta^{(w)}+\lambda^2_n)F_n(w)=0$. Rescale the plane coordinate
$w^{\prime}=\lambda_n w$, which transforms the previous equation to
$(\triangle^{(w^{\prime})}+1)F_n(\frac{w^\prime}{\lambda_n})=0$. The
unique integrable solution for this equation is the Bessel function
$K_0$ (see below), so that $F_n(w)=K_0(\lambda_n w)$. Solutions to
the non-homogeneous equations will then be given by the convolution:
\begin{equation} \label{conv}
g_n(w)=\int_{\real^{2}}F_n(w-x)\rho_n(x)dxd\overline{x}
\end{equation}
and recall that $||g_n||_{L^2}\leq||F_n||_{L^1}||\rho_n||_{L^2}$.
So, all we need is an estimate for $||F_n||_{L^1}$ which is independent
of $n$.

From the expression above, one sees that each $F_n$ is integrable if
the Bessel function $K_0$ is, so that $||F_n||_{L^1}=\lambda_n^{-2}||K_0||_{L^1}$.
So, let $\lambda={\rm min}\{\lambda_n\}_{n\in{\Bbb N}}$; therefore,
$||F_n||_{L^1}\leq\lambda^{-2}||K_0||_{L^1}$; putting $k=\lambda^{-2}||K_0||_{L^1}$
we have \linebreak $||g_n||_{L^2}\leq k||\rho_n||_{L^2}$ for each $n$.
This completes the proof.
\pfend

Consider the Hilbert space $L^2_2(L_\ksi\otimes S^\pm)$ obtained by the completion of
$\Gamma(L_{\ksi}\otimes S^\pm)$ with respect to the norm:
\begin{equation} \label{norm}
||s||_{L^2_2}=||s||_{L^2}+||\triangle_{\ksi}s||_{L^2}
\end{equation}
The map $\triangle_{\ksi}:L^2_2(L_{\ksi}\otimes S^-)\seta L^2(L_{\ksi}\otimes S^-)$
is then bounded, for clearly $||\Delta_\ksi s||_{L^2}\leq||s||_{L^2_2}$.
Let $G_\ksi:L^2(L_{\ksi}\otimes S^-)\seta L^2_2(L_{\ksi}\otimes S^-)$
be the inverse of $\triangle_\ksi$ given by proposition
\ref{est1}. Using the inequality of the proposition, one shows that
$G_\ksi$ is also bounded, if $\ksi$ is nontrivial:
\begin{eqnarray*}
||G_\ksi s||_{L^2_2}& = & ||G_\ksi s||_{L^2}+||\triangle_\ksi G_\ksi s||_{L^2} =
||G_\ksi s||_{L^2}+||s||_{L^2}\leq \\
& \leq & k ||s||_{L^2}+||s||_{L^2} \leq (k+1)\cdot||s||_{L^2}
\end{eqnarray*}
Moreover, we also conclude that:
\begin{equation} \label{gr.est}
||G_\ksi||< 1+\frac{C}{\lambda^2}
\end{equation}

Hence, $G_\ksi$ is an invertible operator when acting between the above
Hilbert spaces, if $\ksi$ is non-trivial.

\vskip18pt

\noindent {\bf Remark:} We emphasise the necessity of assuming that
$\ksi$ is nontrivial. If $\ksi=\hat{e}$, then the equation
(\ref{sep}$(i)$) admits one zero eigenvalue; on the other hand, the
fundamental solution of $\triangle^{(w)}g=0$ is essentially $\log
r$, which is not integrable. It is then impossible to get the
estimate of proposition \ref{est1}, in other words, the operator
$\triangle_{(\ksi=\hat{e})}$ fails to be invertible. In addition,
the parameter $k$ also depends on $\ksi$, and $k\seta\infty$ (i.e.
$\lambda\seta0$) as $\ksi\seta0$.

\vskip18pt

Now, define the norms:
\begin{equation} \label{metrics} \left\{ \begin{array}{l}
||s||_{L^2_1}=||s||_{L^2}+||D^*_{\ksi}s||_{L^2}\ {\rm if}\ s\in \Gamma(L_{\ksi}\otimes S^-)\\
||s||_{L^2_{l+1}}=||s||_{L^2_l}+||D_{\ksi}s||_{L^2_l}\ {\rm if}\ s\in \Gamma(L_{\ksi}\otimes S^+)
\end{array} \right. \end{equation}
and consider the Dirac operators as maps between the following Hilbert
spaces, obtained by the completion of $\Gamma(L_{\ksi}\otimes S^\pm)$
with respect to the above norms:
\begin{equation} \label{mapsops} \left\{ \begin{array}{l}
D^*_\ksi: L^2_1(L_{\ksi}\otimes S^-)\seta L^2(L_{\ksi}\otimes S^+)\\
D_\ksi: L^2_{l+1}(L_{\ksi}\otimes S^+)\seta L^2_l(L_{\ksi}\otimes S^-)
\end{array} \right. \end{equation}
Then $D^*_\ksi$ is clearly bounded. Furthermore, it has an inverse given by \linebreak
$(D^*_\ksi)^{-1}=D_\ksi G_\ksi:L^2(L_{\ksi}\otimes S^+)\seta L^2_1(L_{\ksi}\otimes S^-)$,
which is also bounded:
\begin{eqnarray*}
||(D^*_\ksi)^{-1}s||_{L^2_1}&=&||(D^*_\ksi)^{-1}s||_{L^2}+||D^*_\ksi(D^*_\ksi)^{-1}s||_{L^2}= \\
& = &||D_\ksi G_\ksi s||_{L^2}+||s||_{L^2} \  = \
||D_\ksi G_\ksi s||_{L^2_1} \leq \\
& \leq & ||G_\ksi s||_{L^2_2} \leq (k+1)\cdot||s||_{L^2}
\end{eqnarray*}

So, $D^*_\ksi$ is also Fredholm when acting as in (\ref{mapsops}),
and our proof is complete. To further reference, we shall denote
$Q_\ksi^\infty=(D^*_\ksi)^{-1}$; note, moreover, that this is a
bounded, elliptic, pseudo-differential operator of order $-1$.
\pfend

We are left with one point to establish: the integrability of the
fundamental solution of $(\triangle+1)F=0$ in the plane. Indeed, first note
that since the operator $\triangle+1$ has polar symmetry, then the
fundamental solution $F$ also has. After imposing this symmetry,
we obtain the following ODE, for $r>0$:
$$ (\triangle+1)F(r)=0 \imply F^{\prime\prime}+\frac{1}{r}F^\prime-F=0 $$
This is a Bessel equation with parameter $\nu=0$. Its
solutions are linear combinations of the Bessel functions of
imaginary argument $I_0$ and $K_0$ (see \cite{A}, chapter 11).
Below are possible integral representations for these functions:
\begin{eqnarray*}
K_0(r)=\int_{1}^{\infty}e^{-rt}(t^2-1)^{-\frac{1}{2}}dt
& \cite{GR} & 8.432.3 \\
I_0(r)=\int_{-1}^{1}\cosh(rt)(t^2-1)^{-\frac{1}{2}}dt
& \cite{GR} & 8.431.2
\end{eqnarray*}
It is easy to see that $I_0(r)$ increases exponentially with $r$; it is also
finite for $r=0$. For the purpose of finding a Green's function for the operator
$\triangle+1$, this solution can be eliminated.

With the help of a table of integrals, one finds out that $K_0$ is integrable; indeed:
$$ \int_{\real^2}K_0(r)d^2vol=\int_{0}^{\infty}\int^{2\pi}_{0}K_0(r)rdrd\theta=
2\pi\int_{0}^{\infty}rK_0(r)dr=2\pi $$
by \cite{GR} 6.561.16 (choosing $\mu=1$, $\nu=0$, $a=1$). This means that
\linebreak $||K_0||_{L^1}=2\pi$.

\begin{prop} \label{est2}
The solution $f$ of the flat laplacian problem $\Delta_{\ksi}f=\rho$
of proposition (\ref{est1}) decays exponentially if $\ksi$ is
nontrivial, in the sense that there is a real constant $\lambda>0$
such that:
$$ \lim_{r\seta\infty}e^{\lambda r}|f|<\infty $$
\end{prop}

\pf As $r\seta\infty$, the Bessel function $K_0$ admits the
following asymptotic expansion (\cite{Wa}, p.202):
\begin{equation} \label{k0exp}
K_0(r) \sim \left( \frac{\pi}{2} \right)^{\frac{1}{2}} \frac{e^{-r}}{\sqrt{r}}
              \left[ 1-\frac{1}{8r}+\frac{9}{128r^2}+\dots \right]
\end{equation}
Now since each $\rho_n$ has compact support, it follows from (\ref{conv})
that each $g_n$ will also decay exponentially:
$$ g_n(w) \sim \left( \frac{\pi}{2} \right)^{\frac{1}{2}}\cdot
\int_\Omega \frac{e^{-\lambda_n|w-x|}}{\sqrt{\lambda_n|w-x|}}
\left[ 1-\frac{1}{8\lambda_n|w-x|}+\dots \right]\rho_n(x)
dxd\overline{x} $$
where $\Omega$ is the support of $\rho$. As $|w|\seta\infty$, then also
$|w-x|\sim |w|$ for all $x\in\Omega$. Therefore,
$$ g_n(w) \sim \left( \frac{\pi}{2} \right)^{\frac{1}{2}}
\frac{e^{-\lambda_n|w|}}{\sqrt{\lambda_n|w|}}
\left[ 1-\frac{1}{8\lambda_n|w|}+\dots \right]
\cdot \int_\Omega \rho_n(x) dxd\overline{x},
\ \ {\rm as}\ \ |w|\seta\infty $$
Choosing $0<\lambda<{\rm min}\{\lambda_n\}_{n\in{\Bbb N}}$, the
statement follows from the eigenspace decomposition of $f$
(\ref{f.decomp}) and (\ref{evalue}).
\pfend

In particular, note that $(f/w)$ also belongs to $L^2(L_\ksi\otimes S^+)$.
Define \linebreak $\widetilde{L^2}(L_\ksi\otimes S^+)$ as the space of sections
$\psi$ such that $\psi/w$ is square-integrable. The proposition just proved
implies that the flat model laplacian acting as follows:
$$ \triangle_\ksi:\widetilde{L^2}(L_\ksi\otimes S^{\pm})\seta L^2(L_\ksi\otimes S^{\pm}) $$
is an invertible operator. Since $\triangle_\ksi=D_\ksi D^*_\ksi$, we conclude
that:
\begin{equation} \label{df2}
D_\ksi^*: \widetilde{L^2}(L_\ksi\otimes S^-) \seta L^2(L_\ksi\otimes S^+)
\end{equation}
is also invertible.


\paragraph{Completing the proof of the theorem \ref{dirac-fred}.}
To show that $D_{A_{\ksi}}^*$ is Fredholm, first note that usual
elliptic theory for compact manifolds guarantees the existence of
a parametrix for $D^*_{A_\ksi}$ inside this compact core $T\times K$;
this is a bounded, elliptic, pseudo-differential operator:
$$ Q^K_{A_\ksi}:L^2(E(\ksi)\otimes S^+|_{T\times K})\seta
   L^2_1(E(\ksi)\otimes S^-|_{T\times K}) $$
of order $-1$.

On the other hand, it follows from lemma \ref{model} that:
$$ ||D^*_{A_\ksi} - (D^*_{\ksi_0+\ksi}\oplus D^*_{-\ksi_0+\ksi})||
    _{L^2(T\times D_R)}^2<2\epsilon $$
where $\epsilon$ can be made arbitrarily small. Thus,
$D^*_{A_\ksi}|_{T\times D_R}$ is also invertible for sufficiently
large $R\gg 0$, if $\ksi\neq\as$. Denote this inverse by
$Q^\infty_{A_\ksi}$; this is also a bounded, elliptic, pseudo-differential
operator of order $-1$.

Now choose $\beta_1,\beta_2:\cpx\seta\real$ respectively supported over
$K$ and $D_R$ and satisfying $\beta_1^2+\beta_2^2=1$ everywhere.
We can patch together our two parametrix $Q^K_{A_\ksi}$ and $Q^{\infty}_{A_\ksi}$
in the following way:
\eq \label{parametrix}
P_{A_\ksi}g=\beta_1Q^K_{A_\ksi}(\beta_1g)+\beta_2Q^\infty_{A_\ksi}(\beta_2g)
\end{equation}
This is the same as restricting the section $g$ to $T\times K$
(respectively, \linebreak $T\times D_R$), apply
$Q^K_{A_\ksi}$ ($Q^\infty_{A_\ksi}$) and restricting the result again
to $T\times K$ ($T\times D_R$). Note that $P_{A_\ksi}$ acts as follows:
$$ P_{A_\ksi}:L^2(E(\ksi)\otimes S^+)\seta L^2_1(E(\ksi)\otimes S^-). $$

We want to show that this is a parametrix for $D_{A_{\ksi}}^*$.
In fact, take \linebreak $g\in L^2(E(\ksi)\otimes S^+)$; then:
\begin{eqnarray}
D_{A_{\ksi}}^*P_{A_\ksi}g &=& D_{A_{\ksi}}^*[\beta_1Q^K_{A_\ksi}(\beta_1g)]+
D_{A_{\ksi}}^*[\beta_2Q^\infty_{A_\ksi}(\beta_2g)]= \nonumber \\
& = & \{\beta_1D_{A_{\ksi}}^*Q^K_{A_\ksi}(\beta_1g)+
\beta_2D_{A_{\ksi}}^*Q^\infty_{A_\ksi}(\beta_2g)\}+ \label{xx} \\
& & +\underbrace{d\beta_1.Q^K_{A_\ksi}(\beta_1g)+
d\beta_2.Q^\infty_{A_\ksi}(\beta_2g)} \nonumber \\
& & \hspace{6.5em} S^\infty g \nonumber
\end{eqnarray}
where ``$.$'' means Clifford multiplication.

Since $Q^K_{A_\ksi}$ is a parametrix for $D_{A_{\ksi}}^*$ inside $T\times K$,
the first term (inside brackets) equals the identity plus a compact
operator $S^K$ acting on $\beta_1g$. Similarly, in the second term,
$Q^{\infty}_{A_\ksi}$ is the inverse of the Dirac operator outside $K$.
Together, the first two terms form the identity operator plus $S^K$. Hence:
$$ (D_{A_\ksi}^*P_{A_\ksi}-I)g=S^Kg+S^\infty g $$
where $S^\infty:L^2(E(\ksi)\otimes S^+)\seta L^2(E(\ksi)\otimes S^+)$
is the operator over the brackets in
(\ref{xx}). Since $Q^K_{A_\ksi}$ and $Q^\infty_{A_\ksi}$ are bounded operators,
so is $S^\infty$; we argue that this is also a compact operator.

In fact, let $\widetilde{\partial K}$ denote the closure of the the support of
$d\beta_1$ and $d\beta_2$ (which is an annulus around the boundary of $K$).
Consider the diagram:
$$ \begin{array}{ccc}
L^2(E(\ksi)\otimes S^+)&\stackrel{s}{\longrightarrow}&
L^2_1(E(\ksi)\otimes S^+|_{T^2\times\widetilde{\partial K}})\\
& &\downarrow\ i \\
& &L^2(E(\ksi)\otimes S^+|_{T^2\times\widetilde{\partial K}})\\
& & \cap \\ & &L^2(E(\ksi)\otimes S^+)
\end{array} $$
Now, let $\Upsilon\subset L^2(E(\ksi)\otimes S^+)$ be a
bounded set; since $s$ is a bounded operator, $s(\Upsilon)$ is also
bounded. By the Rellich lemma (see, for instance, \cite{BB}), the
map $i$ is a compact inclusion; note that $\widetilde{\partial K}$
is a compact subset of the plane. Hence, $i(s(\Upsilon))$ is a
relatively compact subset of $L^2(E(\ksi)\otimes S^+|_{T^2\times\widetilde{\partial K}})$,
and clearly also a relatively compact subset of $L^2(E(\ksi)\otimes S^+)$. This means
that $S^\infty=i\circ s:L^2(E(\ksi)\otimes S^+)\seta L^2(E(\ksi)\otimes S^+)$
is a compact operator, as have we claimed.

We conclude that
$$ D_{A_{\ksi}}^*P_{A_\ksi}-I=[{\rm compact\ operator}] $$
and (\ref{parametrix}) is indeed a parametrix for $D_{A_\ksi}^*$ if
$\ksi\neq\pm\ksi_0$.

Finally, to compute the index of $D_{A_\ksi}^*$ we need a relative index
theorem, which is stated and proved in the appendix \ref{apprit}.
There, we show that:

\begin{cor} \label{indmod}
If $A\in\cala_{(k,\ksi_0)}$, then ${\rm index}D^*_{A_\ksi}=k$.
\end{cor}


\paragraph{The Green's operator.} Clearly, the Dirac laplacian, with the norms as
in (\ref{norm}):
\begin{equation} \label{green} \begin{array}{c}
\Delta_{A_\ksi}:L^2_2(E\otimes L_\ksi\otimes S^+)\seta L^2(E\otimes L_\ksi\otimes S^+) \\
\Delta_{A_\ksi}=D_{A_\ksi}^*D_{A_\ksi}
\end{array} \end{equation}
is also a Fredholm operator. In particular, by general Fredholm
theory, there is a bounded operator $G_{A_\ksi}$, called the Green's
operator, such that $\Delta_{A_\ksi}G_{A_\ksi}=Id-H_\ksi$, where
$H_\ksi$ is the finite rank orthogonal projection operator
$H_\ksi:L^2_2(E\otimes L_\ksi\otimes S^+)\seta{\rm ker}(\Delta_{A_\ksi})$.


\subsection{Harmonic spinors and cohomology} \label{spch}

To conclude this chapter, we want to interpret the harmonic spinors \linebreak
$\psi\in{\rm ker}D_A^*$ as some holomorphic object defined in terms of the
compactified bundle $\ee\seta\tproj$. Indeed, we aim to establish the
following identification:

\begin{prop} \label{spin/coho}
If $A$ has nontrivial asymptotic state $\ksi_0\in\dual$ and $k>0$,
then there is an isomorphism $H^1(\tproj,\ee)\equiv{\rm ker}D_A^*$.

\end{prop}

Note that ${\rm ker}D_A^*\subset L^2_1(E\otimes S^-)$, with the
norm defined in (\ref{n.thm}). First, we must show that
$H^1(\tproj,\ee)$ has the correct dimension.

\paragraph{Vanishing theorem.}
Since $\chi(\ee)=-k$, in order to conclude that \linebreak
$h^1(\tproj,\oo(\ee))=k$, it is enough to show that the
cohomologies of orders $0$ and $2$ vanish.

Let us assume that the restriction of $\ee$ to the elliptic
curves $\ee|_{T\times\{w\}}$ is semi-stable for all $w\in\proj$.
We can regard $\ee\seta\tproj$ as a family of extensions:
$$ 0 \seta L_\ksi \seta \ee|_{T_w} \seta L_{-\ksi} \seta 0 $$
of a flat line bundle $L_\ksi$ by its dual $L_{-\ksi}$, where
$\ksi\in\dual$ depends holomorphically on $w\in\proj$; in other
words, the family is parametrised by $\proj$.

Since such extension bundles can be indecomposable if and only if
$\ksi=-\ksi$ (i.e. $\ksi$ has order 2 in $\dual$), we conclude that
$\ee|_{T_w}$ splits as a sum of flat line bundles for all but
finitely many points $w\in\proj$. Furthermore, these flat line
bundles are holomorphically nontrivial for all but finitely many
points $w\in\proj$.

This observation leads to the desired vanishing result:

\begin{lem} \label{vanish}
If $\ee$ is irreducible and $k>0$, then:
$$ h^0(\tproj,\ee)=h^2(\tproj,\ee)=0 $$
\end{lem}

Let $L_\ksi\seta T$ be a flat line bundle with its canonical
connection, as described in section \ref{poin};  denote:
$$ \ee(\ksi)=\ee\otimes p_1^*L_\ksi \ \ \ {\rm and} \ \ \
   \tilde{\ee}(\ksi)=\ee\otimes p_1^*L_\ksi\otimes p^*_2\oo_{\proj}(1) $$
Note that we can regard $p^*_2\oo_{\proj}(1)$ as the line bundle corresponding
to the divisor $T_\infty$. It follows from the lemma that:
$$ h^1(\tproj,\ee(\ksi))=h^1(\tproj,\tilde{\ee}(\ksi))=k $$
for every $\ksi\in\dual$.

\pf Take $w\in\proj$ such that $\ee(\ksi)|_{T_w}=L_{\ksi_1}\oplus L_{\ksi_2}$
for some non-trivial $\ksi_1,\ksi_2\in\dual$; the existence of such
point follows from the observations made prior to the statement of
the lemma. Let $V\subset\proj$ be an open neighbourhood of $w$ such
that every point of $V$ satisfy a the same condition.

Suppose there is a holomorphic section $s\in H^0(M,\ee(\ksi))$; it
gives rise to a holomorphic section $s_w$ of $\ee(\ksi)|_{T_w}\seta
T_w$. On the other hand, we have that \linebreak
$h^0(T,\ee(\ksi)|_{T\times\{w\}})=0$, hence $s_w\equiv0$. Moreover,
$s_w\equiv0$ for all $w\in V$, so that $s$ must vanish identically
on the open set $T\times V$, hence vanish everywhere and
$h^0(\ee(\ksi))=0$. The vanishing of $h^0(\tilde{\ee}(\ksi))$ is
proved in the very same way by noting
$\tilde{\ee}(\ksi)|_{T_w}\equiv \ee(\ksi)|_{T_w}$ since
$p^*_2\oo_{\proj}(1)|_{T_w}=\underline{\cpx}$.

The vanishing of the $h^2$'s follows from Serre duality and a
similar argument for the bundle $\ee(\ksi)\otimes K_{\tproj}$. More
precisely, Serre duality implies that:
$$ \begin{array}{rcl}
H^2(\tproj,\ee(\ksi)) & = & H^0(\tproj,\ee(\ksi)^\vee\otimes K_{\tproj})^* \\
& = & H^0(\tproj,\ee(\ksi)^\vee\otimes p_2^*\oo_{\proj}(-2))^*
\end{array} $$
On the other hand, it is easy to see that:
$$ \ee(-\ksi)|_{T_w} \equiv (\ee(\ksi)^\vee\otimes p_2^*\oo_{\proj}(-2))|_{T_w} $$
so that we can apply the same argument as above to show that
\linebreak $h^0(\tproj,\ee(\ksi)^\vee\otimes K_{\tproj})=0$.
\pfend


\paragraph{Proof of proposition \ref{spin/coho}.}

Let $\{w_i\}$ be the set of points in $\proj$ for which
$h^0(T_{w_i},\ee|_{T_{w_i}})\neq0$. As we argued above, there are
only finitely many such points; in fact, it can be shown that there
are at most $k$ such points (see lemma \ref{bp}). Suppose that
$\#\{w_i\}=p\leq k$; note also that $\infty\notin\{w_i\}$ if
$\ksi_0$ is nontrivial.

Denote by $B$ the divisor in $\tproj$ consisting of the elliptic
curves lying over these points, i.e. $B=\sum_i T\times\{w_i\}$. Also,
denote $\ee(p)=\ee\otimes\oo_{\tproj}(B)$.

Consider the exact sequence of sheaves:
$$ 0\seta\oo(\ee)\seta\oo(\ee(p))\seta\oo(\ee(p)|_B)\seta 0 $$
which induces the following sequence of cohomology:
\footnotesize
\begin{equation} \label{xy}  \begin{array}{ccccc}
0 \seta H^0(B,\ee(p)|_B) \seta & \underbrace{H^1(\tproj,\ee)}
  & \seta & \underbrace{H^1(\tproj,\ee(p))} & \seta H^1(B,\ee(p)|_B) \seta 0 \\
 & {\rm dim}=k & & {\rm dim}=k &
\end{array} \end{equation}
\normalsize \baselineskip18pt
and note that $p\leq h^0(B,\ee(p)|_B)=h^1(B,\ee(p)|_B)\leq2k$. It
follows from (\ref{xy}) that $h^0(B,\ee(p)|_B)=h^1(B,\ee(p)|_B)=k$,
so that the left map in the sequence (\ref{xy}) above
$H^0(B,\ee(p)|_B)\seta H^1(\tproj,\ee)$ is an isomorphism.

This means that each element in $H^1(\tproj,\ee)$ can be represented
by a $(0,1)$-form $\theta$ supported on tubular neighbourhoods of the
fibres $T\times\{w_i\}$. Pulling $\theta$ back to $\torus$, we obtain
a compactly supported $(0,1)$-form, which we also denote by $\theta$,
since $\ksi_0$ is nontrivial.

We want to fashion a solution $\psi$ of $D_A^*\psi=0$ out of
$\theta$, and within the same cohomology class. In other words, by
virtue of the extensibility hypothesis, we want to find a section
$s\in L^2(\Lambda^0E)$ such that $D_A^*(\theta+\del_As)=0$. Since
\linebreak $D_A^*=\del_A^*-\del_A$, this is the same as solving the
equation:
$$ \del_A^*\del_As \ = \ \Delta_As \ = \ -\del_A^*\theta $$
for a compactly supported $\theta$.

In the Fredholm theory for the Dirac operator developed above, we
constructed the Green's operator $G_A$ of the Dirac laplacian
$\Delta_A$. Thus, we can write $s=-G_A\del_A^*\theta$ and
$\psi=\theta-\del_AG_A\del_A^*\theta=P\theta$, where $P$ denotes
the $L^2$ projection $L^2(E\otimes S^-)\stackrel{P}{\seta}{\rm ker}D^*_A$.

We must verify that $\psi\in L^2(E\otimes S^-)$; it is enough to
show that $\del_AG_A\del_A^*\theta$ is square-integrable for any
compactly supported (0,1)-form $\theta$. First note that
$\gamma=\del_A^*\theta$ also has compact support, so that
$s=G_A\gamma\in L^2(\Lambda^0E)$. Therefore, we have:
\begin{eqnarray*}
||\del_A s||_{L^2}^2 & = & \langle \del_A s,\del_A s \rangle \
= \ \langle \del_A s,(\del_AG_A)\gamma \rangle \ = \\
& = & \langle (\del_AG_A)^*\del_A s , \gamma \rangle
\end{eqnarray*}
which is finite, since $\gamma$ is compactly supported. Note the the integration
by parts made from the first to the second line is justified by the same fact.
Therefore, $\psi$ is indeed a square-integrable solution of $D_A^*\psi=0$.

Finally, to see that the map defined above is injective (hence an isomorphism),
let $\theta'$ be another $(0,1)$-form supported around $B$ and within
the same cohomology class as $\theta$, so that $\theta-\theta'=\del_A\alpha$.
Thus:
\begin{equation} \begin{array}{rcl}
\psi-\psi' & = & (\theta-\del_AG_A\del_A^*\theta)-(\theta'-\del_AG_A\del_A^*\theta') = \\
& = & (\theta-\theta')-\del_AG_A\del_A^*(\theta-\theta') = \\
& = & \del_A\alpha-\del_AG_A\del_A^*\del_A\alpha = \del_A\alpha-\del_A\alpha \ = \ 0
\end{array} \end{equation}

This completes the proof. \pfend

\chapter{Nahm transform for instantons over $\torus$} \label{nahm}
\chaptermark{Nahm transform} 

We are finally ready to present the construction of the Nahm
transform for an instanton over $\torus$, proving theorem
\ref{nahmthm}. In the first section, we outline a purely
differential geometric approach to this construction. As we have
mentioned in the introduction, such approach is not powerful enough
due to the non-compactness of $\torus$, but has the virtue of being
very clear and explicit.

Inspired by this gauge-theoretical approach, we bring forth the
powerful tools of algebraic geometry to probe the singularity of
the Higgs field. The compactification results established in the
previous chapter puts us in position to approach the problem in a
holomorphic fashion, completing the proof of theorem \ref{nahmthm}
in chapter \ref{conc}.

\section{Gauge-theoretical construction} \label{diffgeo}

Recall that our starting point is a rank 2 vector bundle \linebreak
$E\seta\torus$ provided with an instanton connection
$A\in\cala_{(k,\ksi_0)}$, where the instanton number $k$ and the
asymptotic state $\ksi_0$ are from now on fixed.

Over the dual torus, consider the trivial Hilbert bundle
$\hat{H}\seta\dual\setminus\{\as\}$ whose fibres are
$\hat{H}_{\ksi}=L^2_1(E(\ksi)\otimes S^-)$. Taking the $L^2_1$-norm
on the fibres, $\hat{H}$ becomes an hermitian bundle. Moreover,
call $\hat{d}$ the trivial connection on $\hat{H}$; such connection
is clearly unitary, hence one can define a holomorphic structure
over $\hat{H}$.

Now, consider the finite-dimensional sub-bundle $V\hookrightarrow\hat{H}$
whose fibres are given by $V_{\ksi}={\rm ker}D_{A_{\ksi}}^*$. We
shall call $V\seta\dual\setminus\{\as\}$ the {\em dual bundle} of
$E$;  remark that this is actually the {\em index bundle} (see
\cite{BB} or \cite{DK}) for the family of Dirac operators
$D_{A_{\ksi}}$. Let $i:V\seta\hat{H}$ be the natural inclusion and
$P:\hat{H}\seta V$ the fibrewise orthogonal $L^2$ projection; more
precisely, $P_\ksi=I-D_{A_\ksi}G_{A_\ksi}D^*_{A_\ksi}$ for each
$\ksi\in\dual\setminus\{\as\}$, where $G_{A_\ksi}$ denotes the
Green's operator for (\ref{green}), $I$ is the identity operator.
We can define a hermitian connection on $V$ via the projection formula:
\eq \label{ft}
\nabla_{B} \ = \ P\circ\hat{d}\circ i
\end{equation}
where $B$ is the associated connection form.

Clearly, $V$ inherits the hermitian metric $h$ from $\hat{H}$, and
$B$ is also unitary with respect to this induced metric. Hence, we
can provide $V$ with the holomorphic structure coming from the
unitary connection $B$.

Alternatively, $V$ also admits an interpretation in terms of {\em monads},
see \cite{DK}. The Dirac operator can be unfolded into a family of elliptic
complexes parametrised by $\dual\setminus\{\as\}$, namely:
\begin{equation} \label{monad}
0 \seta L^2_2(\Lambda^0E(\ksi)) \stackrel{\del_{A_\ksi}}{\longrightarrow}
   L^2_1(\Lambda^{0,1}E(\ksi)) \stackrel{-\del_{A_\ksi}}{\longrightarrow}
   L^2(\Lambda^{0,2}E(\ksi)) \seta 0
\end{equation}
which, of course, are also Fredholm. Moreover, the cohomologies of
order 0 and 2 must vanish, by proposition \ref{vanish}. As in
\cite{DK}, such holomorphic family defines a holomorphic vector
bundle $V\seta(\dual\setminus\{\as\})$, with fibres
$V_\ksi=H^1(\ksi)={\rm ker}D_{A_\ksi}^*$, plus an unitary
connection, induced by orthogonal projection, which is compatible
with the given holomorphic structure. Such connection will be
denoted by $B$. We will invoke this construction repeatedly
throughout this work.

The curvature $F_B$ of $B$ is simply given by:
$$ F_B \ = \ \nabla_B\nabla_B \ = \ P\hat{d}(P\hat{d}) $$
Explicit formulas for the matrix elements on an arbitrary local
trivialisation of $V\seta(\dual\setminus\{\as\})$ will be useful
later on. For instance, pick up an orthonormal frame
$\{\psi_i\}_{n=1}^{k}$ over an open set
$U\subset\dual\setminus\{\as\}$. Then, we have that:
\begin{eqnarray}
(B)_{ij} & = & \langle \psi_j,\nabla_B\psi_i \rangle \ = \
               \langle \psi_j,\hat{d}\psi_j \rangle \nonumber \\
(F_B)_{ij} & = & \langle \psi_j,F_B\psi_i \rangle \ = \
                 \langle \psi_j,P\hat{d}(P\hat{d}\psi_i) \rangle \ = \
                 \langle \psi_j,\hat{d}(P\hat{d}\psi_i) \rangle \label{concurv}
\end{eqnarray}


\paragraph{Higgs field.}
We now define the Higgs field $\Phi\in{\rm End}(V)\otimes
K_{\dual}$. Recall that $w$ is the complex coordinate of the plane.
Let $\psi\in\Gamma(V)$, i.e. for each
$\ksi\in\dual\setminus\{\as\}$, $\psi[\ksi]\in{\rm
ker}D^*_{A\ksi}$. For a fixed $\ksi'$, the Higgs field will act on
$\psi[\ksi']$ by multiplying this section by the plane coordinate
$w$ and then projecting it back to ${\rm ker}D^*_{A_\ksi}$:
\begin{equation} \label{higgs.gt}
(\Phi(\psi))[\ksi'] \ = \ \frac{1}{\sqrt{2}}P_{\ksi'}(w\psi[\ksi'])d\ksi
\end{equation}
Its conjugate is clearly given by
$(\Phi^*(\psi))[\ksi'] \ = \ \frac{1}{\sqrt{2}}P_{\ksi'}
(\overline{w}\psi[\ksi'])d\overline{\ksi}$

Again, there is a subtle analytical point here. The spinors $\psi$
belong to $L^2(E(\ksi)\otimes S^-)$ but is not necessarily the case
that $w\psi$ also belong to $L^2(E(\ksi)\otimes S^-)$. To show this
is indeed the case, we have the following lemma:
\begin{lem}
If $\psi\in{\rm ker}D_A^*$ and $A$ has nontrivial asymptotic state, then
\linebreak $w\psi\in L^2(E\otimes S^-)$.
\end{lem}
\pf
The key result here is proposition \ref{est2}, and the observation
that follows it, in particular the invertibility of the operator
(\ref{df2}).

Let $K\subset\torus$ be a compact subset such that $D_A^*$ is
sufficiently close to the flat Dirac operator $D_{\as}^*$ outside
$K$. Thus, restricted to the complement of $K$, $D_A^*$ is
invertible acting from $\tilde{L^2}\seta L^2$.

Now if $\psi\in{\rm ker}D_A^*$, then $D_A^*(w\psi)=dw\cdot\psi\in
L^2(E(\ksi)\otimes S^+|_{\torus\setminus K})$ and the proposition
follows.
\pfend

Note that the dependence of $(B,\Phi)$ on the original instanton
$A$ is contained on the $L^2$-projection operator $P$, i.e. on the
$k$ solutions of $D_{A_\ksi}^*\psi=0$. It is easy to see that the
finite dimensional space spanned by these $\psi$ is gauge
invariant; moreover, the multiplication by $w$ also commutes with
gauge transformations $\hat{g}\in{\rm Aut}(V)$. Therefore, we have
that:

\begin{prop} \label{eqv2}
If $A$ and $A'$ are gauge equivalent irreducible instantons, then
the corresponding pairs $(B,\Phi)$ and $(B^\prime,\Phi^\prime)$ are
also gauge equivalent.
\end{prop}

A pair $(B,\Phi)$ is called a {\em Higgs pair} on the bundle
$V\seta\dual\setminus\{\as\}$ if it satisfies Hitchin's
self-duality equations:
\eq \label{hiteq2}
\left\{ \begin{array}{l}
{\rm (i)}\ F_{B}+[\Phi,\Phi^*]=0 \\
{\rm (ii)}\ \overline{\partial}_{B}\Phi=0
\end{array} \right. \end{equation}

Recall from section \ref{poin} that the unitary connection, and its
corresponding curvature, of the Poincar\'e line bundle
${\bf P}\seta T\times\dual$ are given by:
\begin{eqnarray*}
\omega(z,\ksi)=i\sum_{\mu=1}^{2}\ksi_\mu dz_\mu & \ \ &
\Omega(z,\ksi)=i\sum_{\mu=1}^{2}d\ksi_\mu\wedge dz_\mu
\end{eqnarray*}
From Braam \& Baal \cite{BVB}, we know that if $s\in\Gamma(E(\ksi)\otimes S^-)$,
then:
\begin{equation} \label{commut}
D^*_{A_\ksi}(\hat{d}s) \ = \ [D^*_{A_\ksi},\hat{d}]s \ = \ -\Omega\cdot s
\end{equation}
where $\cdot$ means Clifford multiplication. The local formula for the
curvature (\ref{concurv}) may now be cast on a more convenient form:
\begin{eqnarray*}
(F_B)_{ij} & = & \langle \psi_j,\hat{d}(P\hat{d}\psi_i) \rangle \ = \
\langle \psi_j,\hat{d}(D_{A_\ksi}G_{A_\ksi}D^*_{A_\ksi}\hat{d}\psi_i) \rangle \ = \\
& = &  \langle -D^*_{A_\ksi}\hat{d}\psi_j,G_{A_\ksi}(D^*_{A_\ksi}\hat{d}\psi_i) \rangle \ = \
\langle \Omega\cdot\psi_j,G_{A_\ksi}(\Omega\cdot\psi_i) \rangle
\end{eqnarray*}
Since the Clifford multiplication commutes with the Green's operator, we end up with:
\begin{eqnarray}
(F_B)_{ij}& = & - \langle (\Omega\wedge\Omega)\cdot\psi_i,G_{A_\ksi}\psi_i \rangle
                \ = \nonumber \\
& = & 2 \langle (dz_1\wedge dz_2)\cdot\psi_j,G_{A_\ksi}\psi_i \rangle
      d\ksi_1\wedge d\ksi_2 \ = \label{curv} \\
& = & -i \langle (dz_1\wedge dz_2)\cdot\psi_j,G_{A_\ksi}\psi_i \rangle
      d\ksi\wedge d\overline{\ksi} \nonumber
\end{eqnarray}
Note moreover that the inner product is taken in $L^2(E(\ksi)\otimes S^-)$,
integrating out the $(z,w)$ coordinates.


\paragraph{Hitchin's pairs from instantons.}
Our first step towards the proof of theorem \ref{nahmthm} is the following result:

\begin{thm} \label{soln}
If $A$ is an irreducible, extensible instanton connection on
$E\seta\torus$, then the associated pair $(B,\Phi)$ on the dual
bundle $V\seta\dual\setminus\{\as\}$ constructed above satisfies
the Hitchin's equations (\ref{hiteq2}).
\end{thm}

\pf Choose a point $\xi$ and an open neighbourhood
$\xi\in U\subset\dual\setminus\{\as\}$ and pick up a local
orthonormal trivialisation of $V\seta\dual\setminus\{\as\}$ over
$U$, such that the corresponding local frame $\{\psi_i\}_{n=1}^{k}$
is parallel at $\ksi$. Recall that $\psi_i(\ksi)\in
kerD^*_{A_\ksi}$.

First, we shall look at the second equation of (\ref{hiteq2}), and recall that \linebreak
$\dual\setminus\{\as\}$ was given the flat Euclidean metric induced from the quotient.
Once a local trivialisation is chosen, the endomorphism $\Phi$ can then be put
in matrix form, with matrix elements given by:
$$ a_{ij}(\ksi) \ = \ \langle \psi_j(\ksi),\Phi[\psi_i](\ksi) \rangle $$
where $\langle,\rangle$ is the inner product on $L^2(E(\ksi)\otimes S^-)$,
integrating out the $(z,w)$ coordinates. Clearly, $\Phi$ is a
holomorphic endomorphism if its matrix elements in holomorphic
trivialisation are holomorphic functions. However:
$$ \Phi[\psi_i](\ksi) \ = \ P_\ksi(w\psi_i(\ksi))d\overline{\ksi} \ = \
(I-D_{A_\ksi}G_{A_\ksi}D^*_{A_\ksi})(w\psi_i(\ksi))d\overline{\ksi} $$
so that:
\begin{eqnarray*}
a_{ij}(\ksi) & = & \frac{1}{\sqrt{2}}\left\{\langle \psi_j(\ksi),w\psi_i(\ksi) \rangle-
\langle \psi_j(\ksi),D_{A_\ksi}G_{A_\ksi}D^*_{A_\ksi}(w\psi_i(\ksi)) \rangle \right\} \ = \\
& = & \frac{1}{\sqrt{2}}\left\{\langle\psi_j(\ksi),w\psi_i(\ksi)\rangle-
\langle D^*_{A_\ksi}\psi_j(\ksi),G_{A_\ksi}D^*_{A_\ksi}(w\psi_i(\ksi))\rangle\right\} \ = \\
& = & \frac{1}{\sqrt{2}}\langle\psi_j(\ksi),w\psi_i(\ksi)\rangle
\end{eqnarray*}
Therefore:
\begin{eqnarray*}
\frac{\partial a_{ij}}{\partial\overline{\ksi}}(\ksi) &=&
\frac{1}{\sqrt{2}}\left\{ \langle \partial_{B}\psi_j,w\psi_i\rangle +
\langle\psi_j,\del_{B}(w\psi_i) \rangle \right\} \ =
\\ & = & \frac{1}{\sqrt{2}}\langle \psi_j,
         \left(\frac{\partial w}{\partial\overline{\ksi}}\right)\psi_i+
         \del_{B}\psi_i \rangle \ = \ 0
\end{eqnarray*}
as $\psi_i$ is parallel at $\ksi$. Since this can be done for all
$\ksi\in\dual\setminus\{\as\}$, the second equation is satisfied.

Now, we move back to (\ref{hiteq2}(i)). Let us first
compute the matrix elements $([\Phi,\Phi^*])_{ij}$. Note that:
\eq \label{ids} \left\{ \begin{array}{l}
(i)\ [D_{A_\ksi}^*,\overline{w}]\psi_i(\ksi) \ = \ D_{A_\ksi}^*(\overline{w}\psi_i(\ksi)) \ =
-d\overline{w}\cdot\psi_i(\ksi)\\
(ii)\ [D_{A_\ksi}^*,w]\psi_i(\ksi) \ = \ D_{A_\ksi}^*(w\psi_i(\ksi)) \ = \ 0
\end{array} \right. \end{equation}
where we used the fact that $D_{A_\ksi}=\del^*_{A_\ksi}-\del_{A_\ksi}$.

Recall that for 1-forms $[\Phi,\Phi^*]=\Phi\Phi^*+\Phi^*\Phi$.
We compute each term separately:
\begin{eqnarray*}
\Phi^*\Phi(\psi_i) & = & \frac{1}{2}P[\overline{w}P(w\psi_i)]d\ksi\wedge d\overline{\ksi} \ =\\
&=&\frac{1}{2}\left\{\overline{w}P(w\psi_i)-D_{A_\ksi}G_{A_\ksi}
D^*_{A_\ksi}\overline{w}P(w\psi_i)\right\}d\ksi\wedge d\overline{\ksi} \ =\\
&=& \frac{1}{2}\left\{\overline{w}w\psi_i-
\overline{w}D_{A_\ksi}G_{A_\ksi}D^*_{A_\ksi}(w\psi_i)-\right.\\
& &\left.-D_{A_\ksi}G_{A_\ksi}D^*_{A_\ksi}\overline{w}P(w\psi_i)\right\}
d\ksi\wedge d\overline{\ksi} \\
\Phi\Phi^*(\psi_i) & = & \frac{1}{2}P[wP(\overline{w}\psi_i)]d\overline{\ksi}\wedge d\ksi \ =\\
&=& \frac{1}{2}\left\{w\overline{w}\psi_i-wD_{A_\ksi}G_{A_\ksi}
D^*_{A_\ksi}(\overline{w}\psi_i)-\right.\\
& &\left.-D_{A_\ksi}G_{A_\ksi}D^*_{A_\ksi}wP(\overline{w}\psi_i)\right\}
d\overline{\ksi}\wedge d\ksi
\end{eqnarray*}

The two first terms of $\Phi\Phi^*$ and $\Phi^*\Phi$ cancel each other
and the third terms will cancel out when we take the inner product
with $\psi_j$. Moreover, the second term of $\Phi^*\Phi$ is zero by
(\ref{ids}(ii)). So we are left with:
\begin{eqnarray*}
([\Phi,\Phi^*])_{ij} &=&\! \frac{1}{2}\langle \psi_j,[\Phi,\Phi^*]\psi_i \rangle \ =
\frac{1}{2} \langle \psi_j,wD_{A_\ksi}G_\ksi
D^*_{A_\ksi}(\overline{w}\psi_i) \rangle\ d\ksi\wedge d\overline{\ksi} = \\
&=& \frac{1}{2}\langle D^*_{A_\ksi}(\overline{w}\psi_j),G_\ksi
D^*_{A_\ksi}(\overline{w}\psi_i) \rangle\ d\ksi\wedge d\overline{\ksi} = \\
&=& -\frac{1}{2} \langle (dw\wedge d\overline{w})\cdot\psi_j,G_\ksi\psi_i \rangle\
    d\ksi\wedge d\overline{\ksi} = \\
&=&\! -i \langle (dw_1\wedge dw_2)\cdot\psi_j,G_\ksi\psi_i \rangle d\ksi\wedge d\overline{\ksi}
\end{eqnarray*}
where we have once more used the fact that the Clifford multiplication commutes with the
Green's operator. Summing the final expression above with ({\ref{curv}), one gets:
$$ (F_B)_{ij}+([\Phi,\Phi^*])_{ij} =
-i \langle (dz_1\wedge dz_2+dw_1\wedge dw_2)\cdot\psi_j,G_\ksi\psi_i \rangle
d\ksi\wedge d\overline{\ksi} = 0 $$
for the first term of the inner product is zero since it consists of a
self-dual form (the K\"ahler form $\kappa$) acting on a negative spinor.
\pfend

\vskip8pt

Clearly, the above result has two weak points: it tells nothing
about the behaviour of the Higgs field around the singular points
$\as$; and it fails to show that the Higgs pairs so obtained are
admissible. In fact, establishing the first point requires the use
of algebraic-geometric methods, and will be taken up in section
\ref{holo} below. The second point will be clarified in section
\ref{inv} when we give the inverse construction, obtaining
instantons from singular Higgs pairs.



\section{Holomorphic approach} \label{holo}

The vanishing results of section \ref{spch} put us in position to
define the transformed bundle $\vv\seta\dual$. Indeed, consider the
following elliptic complex:
\begin{equation} \label{cpx1}
0 \seta L^2_2(\Lambda^0\ee(\ksi)) \stackrel{\del_{A_\ksi}}{\seta}
L^2_1(\Lambda^{0,1}\ee(\ksi)) \stackrel{-\del_{A_\ksi}}{\seta}
L^2(\Lambda^{0,2}\ee(\ksi)) \seta 0
\end{equation}
According to proposition \ref{vanish}, $H^1(\tproj,\ee(\ksi))$ is
the only nontrivial cohomology of this complex. It then follows
that the family of vector spaces given by
$\vv_\ksi=H^1(\tproj,\ee(\ksi))$ forms a holomorphic vector bundle
of rank $k$ over $\dual$; denote such holomorphic structure by
$\del_{\vv}$. Note that $\vv_\ksi$ is defined even if $\ksi=\as$.
Furthermore, by proposition \ref{spin/coho}, $\vv|_{\dual\setminus\as}$
coincides holomorphically with the dual bundle $V$ defined on the
previous section, i.e.:
$$ (\vv,\del_{\vv})|_{\dual\setminus\{\as\}}\simeq(V,\del_B) $$

Moreover, $\vv$ comes equipped with a hermitian metric $h'$, which we
want to compare with $h$, the hermitian metric on $V$ induced from the
monad (\ref{monad}). The key point is a fact we noted before in lemma
\ref{model}: given an 1-form $a$ on $\tproj$, its $L^2$-norm with respect
to the round metric is always larger than its $L^2$-norm with respect
to the flat metric on $T\times(\proj\setminus\{\infty\})$:
$$ ||a||_{L^2_R}>||a||_{L^2_F} $$
Thus, comparing the monads (\ref{monad}) and (\ref{cpx1}), one sees that
$h$ is bounded above by $h'$. In particular, the metric $h$ is bounded
at $\as$.

We can regard $\vv$ as an {\em index bundle} for the family of
Dirac operators over $\tproj$ parametrised by $\ksi\in\dual$.
Hence, its degree can be computed by the Atiyah-Singer index
theorem for families. Consider now the bundle
${\bf G}=p_{12}^*\ee\otimes p_{13}^*{\bf P}$ over $\tproj\times\dual$,
and note that ${\bf G}|_{\tproj\times\{\ksi\}}=\ee(\ksi)$. Then we
have:
\begin{eqnarray*}
ch(\vv) & = & - ch({\bf G})\cdot td(\tproj)/[\tproj] \ = \\
& = & - \left( 2+2c_1({\bf P})+c_1({\bf P})^2-c_2(\ee) \right)
\left( 1+\frac{1}{2}c_1(\proj) \right)/[\tproj] \ = \\
& = & k-\frac{1}{2}c_1({\bf P})^2c_1(\proj)/[\tproj] \ = \
k-2\hat{t}
\end{eqnarray*}
where the ``$-$" sign in the first line is needed since $\vv$ is
formed by the null spaces of the adjoint Dirac operator.

Summing up:
\begin{lem} \label{degree}
The dual bundle $(V,\del_B)\seta\dual\setminus\{\as\}$ admits a
holomorphic extension $\vv\seta\dual$ of degree $-2$. Moreover, its
hermitian metric $h$ is bounded above at the punctures $\as$.
\end{lem}

The determinant line bundle of $\vv$ is not fixed, however. In fact,
let $t_x:\tproj\seta\tproj$ be the translation of the torus by $x\in T$,
acting trivially on $\proj$, and let $\ee'=t_x^*\ee$. If $\vv'$ is
the dual bundle associated with $\ee'$ then $\vv'=\vv\otimes L_x$.
Indeed:
\begin{eqnarray*}
\vv_\ksi' = H^1(\tproj,\ee'(\ksi)) & = &
H^1\left( \tproj,p_{12}^*(t_x^*\ee)\otimes p_{13}^*{\bf P}|_{\tproj\times \{ \ksi \} } \right)= \\
& = & H^1 \left( \tproj,t_x^*(p_{12}^*\ee\otimes p_{13}^*{\bf P})\otimes p_3^*L_x
      |_{\tproj\times \{ \ksi \} } \right) = \\
& = & H^1 \left( \tproj,p_{12}^*\ee\otimes p_{13}^*{\bf P}|_{\tproj\times \{ \ksi \} } \right)
      \otimes (L_x)_\ksi \\
\Rightarrow \ \ \vv_\ksi' & = & \vv_\ksi\otimes(L_x)_\ksi
\end{eqnarray*}
as a canonical isomorphism for each $\ksi\in\dual$. Thus
$\vv'=\vv\otimes L_x$.

Note also that if $B$ is an admissible connection, $\vv$ admits no
splitting $\vv=\vv_0\oplus L$ compatible with $B$ for any flat line
bundle $L$.


\paragraph{Defining the Higgs field.}
The next step is to give a holomorphic description of the Higgs field $\Phi$.

Recall that $h^0(\tproj,p_2^*\oo_{\proj}(1))=2$, and
regarding $\proj=\cpx\cup\{\infty\}$, we can fix two holomorphic
sections $s_0,s_\infty\in H^0(\proj,\oo_{\proj}(1))$ such that
$s_0$ vanishes at $0\in\cpx$ and $s_\infty$ vanishes at the point
added at infinity. In homogeneous coordinates
$\{(w_1,w_2)\in\cpx^2|w_2\neq0\}$ and $\{(w_1,w_2)\in\cpx^2|w_1\neq0\}$,
we have that, respectively ($w=w_1/w_2$):
\begin{eqnarray*}
s_0(w) = w & \ \ \ & s_0(w) = 1 \\
s_\infty(w) = 1 & \ \ \ & s_\infty(w) = \frac{1}{w}
\end{eqnarray*}

Let us first consider an alternative definition of the transformed Higgs
field. For each $\ksi\in\dual$, we define the map:
\begin{equation} \begin{array}{rcl}
H^1(\tproj,\ee(\ksi))\times H^1(\tproj,\ee(\ksi)) &
\stackrel{\Psi_\ksi}{\longrightarrow} & H^1(\tproj,\tilde{\ee}(\ksi)) \\
(\alpha,\beta) & \mapsto & \alpha\otimes s_0-\beta\otimes s_\infty
\end{array} \end{equation}
If $(\alpha,\beta)\in{\rm ker}\Psi_\ksi$, we define an endomorphism
$\varphi$ of $H^1(\tproj,\ee(\ksi))$ at the point $\ksi\in\dual$ as
follows:
\begin{equation}
\varphi_\ksi(\alpha)=\beta
\end{equation}

We check that $\varphi$ actually coincides with the
Higgs field $\Phi$ we defined in the previous section, which is
part of the transformed Higgs pair. Note that:
$$ \alpha\otimes s_0-\beta\otimes s_\infty=0 \ \ \Leftrightarrow \ \
   \beta = \alpha(\otimes s_0)(\otimes s_\infty)^{-1} $$
Moreover, recall that, for any trivialisation of $\oo_{\proj}(1)$ with
local coordinate $w$ on $\proj$, the quotient $s_0(w)/s_\infty(w)=w$.
The claim now follows from the proof of proposition \ref{spin/coho};
we denote $\Phi_\ksi=\varphi_\ksi$.


\begin{prop} \label{higgs}
The eigenvalues of the Higgs field $\Phi$ have at most simple poles
at $\as$. Moreover, the residues of $\Phi$ are semi-simple and have
rank $\leq2$  if $\ksi_0$ is an element of order 2 in the Jacobian
of $T$, and rank $\leq1$ otherwise.
\end{prop}

\pf Suppose $\alpha(\ksi)$ is an eigenvector of $\Phi_\ksi$ with
eigenvalue $\epsilon'(\ksi)=1/\epsilon(\ksi)$, i.e.
$\Phi_\ksi(\alpha(\ksi)) = \epsilon'(\ksi)\cdot \alpha(\ksi)$.
Thus,
$$ \alpha(\ksi)\otimes s_0 - \epsilon'(\ksi)\cdot \alpha(\ksi)\otimes s_\infty =0
   \ \ \Rightarrow \ \
   \alpha(\ksi)\otimes(\epsilon(\ksi)\cdot s_0-s_\infty) =0$$
Therefore, denoting
$s_\epsilon(\ksi)=\epsilon(\ksi)\cdot s_0-s_\infty$,
we have that $\alpha(\ksi)\in{\rm ker}(\otimes s_\epsilon(\ksi))$.

On the other hand, consider the sheaf sequence:
$$ 0 \seta \ee(\ksi) \stackrel{\otimes s_\epsilon(\ksi)}{\seta}
     \widetilde{\ee}(\ksi) \seta \widetilde{\ee}(\ksi)|_{T_{\epsilon'(\ksi)}} \seta 0 $$
since the section $s_\epsilon(\ksi)$ vanishes at $\epsilon'(\ksi)$. It
induces the cohomology sequence:
\begin{equation} \label{sqc2}
0 \seta H^0(T_{\epsilon'(\ksi)},\tilde{\ee}(\ksi)|_{T_{\epsilon'(\ksi)}})
  \seta H^1(\tproj,\ee(\ksi)) \stackrel{\otimes s_{\epsilon}(\ksi)}{\seta} ...
\end{equation}
so that ${\rm ker}(\otimes s_\epsilon(\ksi))=
H^0(T_{\epsilon'(\ksi)},\tilde{\ee}(\ksi)|_{T_{\epsilon'(\ksi)}})$
which is non-empty if and only if
$\ee(\ksi)|_{T_{\epsilon'(\ksi)}}=L_\ksi\oplus L_{-\ksi}$ or ${\bf F}_2\otimes L_\ksi$.

Hence, as $\ksi$ approaches $\as$, we must have that one of the eigenvalues
of $\Phi$, say $\epsilon'(\xi)$ approaches $\infty$,
since $\ee|_{T_\infty}=L_{\xi_0}\oplus L_{-\xi_0}$. Moreover,
$s_\epsilon(\ksi)\seta s_\infty$, so that:
$$ \lim_{\ksi\seta\as}\alpha(\ksi) \in {\rm ker}(\otimes s_\infty)=
   H^0(T_\infty,\ee(\ksi)|_{T_\infty}) $$

Therefore, we conclude that, if $\xi_0\neq-\xi_0$, then one of the eigenvalues of
$\Phi$ has a simple pole at $\as$ since $h^0(T_\infty,\ee(\as)|_{T_\infty})=1$;
similarly, if $\xi_0=-\xi_0$, then two of the eigenvalues of
$\Phi$ have a simple poles at $\ksi_0$.

Note in particular that the images of the residues
of $\Phi$ at $\as$ are precisely given by:
$$ H^0(T_\infty,\tilde{\ee}(\as)|_{T_\infty}) \subset H^1(\tproj,\ee(\as)) $$
\pfend

This proposition almost concludes one way of the correspondence in the
statement of our main theorem; it only remains to be
shown that the Nahm transformed Higgs pair is admissible. We must then
show how to obtain an instanton connection $\check{A}$ on a bundle
$\check{E}\seta\torus$ from a singular Higgs pair, and match these
with the original objects $A$ and $E\seta\torus$. These tasks are
taken up in the following chapter.

\paragraph{A conjecture regarding the hermitian metric on ${\mathbf V}$.}
So far, we only know that the hermitian metric $h$ on the Nahm
transformed bundle is bounded above. Unfortunately, this is not
enough for the construction of the inverse transform in the next
chapter, where we shall need a precise knowledge of the behaviour
of $h$ at the punctures $\as$. More precisely, we must assume that:

\begin{quote}
The hermitian metric $h$ is non-degenerate along the
kernel of the residues of $\Phi$. Furthermore, in a holomorphic
trivialisation of $V$ over a sufficiently small neighbourhood around
$\as$, $h\sim O(r^{1\pm\alpha})$ along the image of the residues
of $\Phi$, for some alpha $0\leq\alpha<1/2$.
\end{quote}

In fact, we expect that $h$ indeed satisfy this assumption. However, further
technical work is necessary to establish this claim.

\paragraph{Final remarks.}
Before we proceed, let us make a few remarks about the proposition
\ref{higgs} above. In \cite{S}, a Higgs field is said to be {\em tame}
if its eigenvalues have at most simple poles. Kovalev has shown that,
if $(B,\Phi)$ is a Higgs pair on the punctured surface, this condition
is equivalent to the following regularity condition \cite{K}:
\begin{equation} \label{kov.reg}
\int_{D_0} \left( |\ksi|^2|F_B|^2 + |\nabla_B\hat{\Phi}|^2 \right)
d\ksi d\overline{\ksi}<\infty
\end{equation}
where $D_0$ is a punctured disc centred at $\as$ with complex
coordinate $\ksi$, and:
$$ \hat{\Phi}=\ksi\frac{\partial}{\partial\ksi}\llcorner\Phi $$

In other words, proposition \ref{higgs} shows that the transformed
Higgs pair $(B,\Phi)$ is {\em regular} in the sense of Kovalev,
i.e. it satisfies condition (\ref{kov.reg}) above. A direct proof of
the regularity condition (\ref{kov.reg}) within the gauge-theoretical
framework of section \ref{diffgeo} is possible; it involves an estimate
of the operator norm $||G_{A_\ksi}||$ as $\ksi\seta\as$, as in
(\ref{gr.est}). However, such approach  would not give the precise
form of the residue obtained in proposition \ref{higgs}.

Finally, we would like to emphasise that the transformed Higgs data
$(V,B,\Phi)$ depend on the original instanton connection only through
the induced holomorphic structure $\del_A$. Indeed, $(V,B,\Phi)$ arise
by looking at the kernel of the adjoint Dirac operator, which depend
only on the holomorphic structure on $E\seta\torus$ (which in turn depend on
the choice of complex structure on $\torus$) and on the choice of metric
on the base. Note also that the holomorphic structure $\del_A$ is entirely
encoded on the extended bundle $\ee\seta\tproj$. That is why we were able
to give a completely holomorphic description of the transform despite the
fact that, in principle, the extended holomorphic bundle $\ee$ contains
{\em less} information than $(E,A)$.


\section{A $\torus\times S^1$ action on the moduli space of instantons}
\label{actsec}

Seen as an abelian group, $\torus\times S^1$ acts on $\torus$ as
follows:
\begin{eqnarray}
(\torus\times S^1)\times(\torus) & \seta & \torus \nonumber \\
(x,y,\gamma)\cdot(z,w) & \seta & (z+x,e^{i\gamma}\cdot w+y) \label{act}
\end{eqnarray}
Clearly, this action lifts to an action of $\torus\times S^1$ on the
moduli space of extensible instantons. We are interested in
understanding the effect of this action on the Nahm transformed Higgs
pairs.

So let $t_{(x,y,\gamma)}(z,w)=(z+x,e^{i\gamma}\cdot w+y)$ and denote
$E'=t_{(x,y,\gamma)}^*E$, \linebreak $A'=t_{(x,y,\gamma)}^*A$ and
$\ee'=t_{(x,y,\gamma)}^*\ee$. Let $(V',\vv',B',\Phi')$ and
$(V,\vv,B,\Phi)$ be the corresponding objects obtained via Nahm
transform on $(E,\ee,A)$ and $(E',\ee',A')$.

Setting $y=\gamma=0$, we have seen that the effect of translations
on the torus $t_x^*$ is simply to add a flat tensor factor, i.e.:
$$ \vv'=\vv\otimes L_x $$
Of course, bundle $V$ and the connection $B$ are similarly twisted.
It is easy to see from the definition that the Higgs field remains
unaltered: $\Phi'=\Phi$.

Now set $x=0$. One sees from the calculations following lemma
\ref{degree} that $t_{(y,\gamma)}^*$ has no effect on the dual
bundle $\vv$, i.e. $\vv'=\vv$. On the other hand, (\ref{higgs.gt})
tells us that the Higgs field varies in a particularly simple way:
$$ \Phi'=e^{i\gamma}\cdot\Phi+y\cdot I $$
Clearly, the action of $t_\gamma^*$ multiplies the residues of
$\Phi$ by $e^{i\gamma}$, while the action of $t_y^*$ leaves them
unchanged.

\chapter{Constructing instantons via the inverse transform}
\label{inv}  \chaptermark{Inverse transform} 

Our task now is to construct a holomorphic rank 2 vector bundle
over $\torus$, with an instanton connection on it, departing from a
suitable singular Higgs pair. We will later show that these coincide
with the original objects from which we started in section \ref{diffgeo}.

Let $V\seta\dual\setminus\{\as\}$ be a hermitian, holomorphic
vector bundle of rank $k$ with a Higgs pair $(B,\Phi)$, as
described in theorem \ref{nahmthm}. More precisely, the connection
$B$ defines a holomorphic structure $\del_B$ on the bundle $V$,
which is also compatible with the hermitian metric; and
$\Phi\in{\rm End}V\otimes K_{\dual}$ has simple poles at $\as$,
with semi-simple residues of rank$\leq2$. Recall also that a
$(B,\Phi)$ is said to be {\em admissible} if there are no
covariantly constant sections of $V$, in other words, if the
following holds for every section $s\in\Gamma(V)$ which is not
constant:
\begin{equation} \label{adm}
||\nabla_Bs||_{L^2}>0
\end{equation}

Motivated by lemma \ref{degree}, we assume also that there is a
hermitian, holomorphic vector bundle $\vv\seta\dual$ of degree $-2$
such that:
$$ (\vv,\del_{\vv})|_{\dual\setminus\{\as\}}\simeq(V,\del_B) $$
Moreover, the hermitian metric on $V$ is bounded above by the
hermitian metric $\vv$.

Of course, this rigid set-up is motivated by the Nahm transform
construction of the previous chapter.

Let $S^+=\Lambda^0\oplus\Lambda^{1,1}$ and
$S^-=\Lambda^{1,0}\oplus\Lambda^{0,1}$. The idea is to study the
following elliptic operators:
\begin{eqnarray}
{\cal D}:\Gamma(V\otimes S^+) \seta \Gamma(V\otimes S^-) & &
{\cal D}^*:\Gamma(V\otimes S^-) \seta \Gamma(V\otimes S^+) \nonumber \\
{\cal D}=(\del_B+\Phi)-(\del_B+\Phi)^* & &
{\cal D}^*=(\del_B+\Phi)^*-(\del_B+\Phi) \label{tmp.dirac}
\end{eqnarray}
where $(B,\Phi)$ is a Higgs pair. Note that the operators in
(\ref{tmp.dirac}) are just the Dirac operators coupled to the
connection $\widetilde{B}$, obtained by lifting the Higgs pair
$(B,\Phi)$ to an invariant ASD connection $(\real^4)^*$ as in
the introduction. In particular, $D_B=\del_B-\del_B^*$ is
the coupled Dirac operator acting on $V\otimes S^-$.

Due to the non-compactness of the base space, the choice of metric
in $\dual\setminus\{\as\}$ is a delicate issue. From the point of
view of the Nahm transform, it is important to consider the Euclidean,
incomplete metric on the punctured dual torus, as we explained in
the introduction. However, such a choice of metric is not a good
one from the analytical point of view. For instance, one cannot
expect on general grounds to have a finite dimensional moduli space
of Higgs pairs.

Fortunately, as we mentioned before, Hitchin's equations are
conformally invariant, so that we are allowed to make conformal
changes in the Euclidean metric localised around the punctures to
obtain a complete metric on \linebreak $\dual\setminus\{\as\}$.
Thus, our strategy is to obtain results concerning the Euclidean
metric from known statements about complete metrics.

In \cite{B3}, Biquard considered the so-called {\em Poincar\'e metric},
which is defined as follows. We perform a conformal change on the
incomplete metric over the punctured torus localised on small punctured
neighbourhoods $D_0$ of $\as$, so that if $\ksi=(r,\theta)$ is a local
coordinate on $D_0$, we have the metric:
\begin{equation} \label{pmetric}
ds^2_P \ = \ \frac{d\ksi d\overline{\ksi}}{|\ksi|^2\log^2|\ksi|^2} \
= \ \frac{dr^2}{r^2\log^2r}+\frac{d\theta^2}{4\log^2r}
\end{equation}
We denote the complete metric so obtained by $g_P$. The Euclidean
metric is denoted by $g_E$. Whenever necessary, we will denote by
$L^2_E$ and $L^2_P$ the Sobolev norms in $\Gamma(\Lambda^*V)$ with
respect to $g_E$ and $g_P$, respectively, together with the
hermitian metric in $V$.


\paragraph{Admissibility and vanishing theorem.}
The next step is to prove that the admissibility condition
(\ref{adm}) implies the vanishing of the $L^2$-kernel of
${\cal D}$:

\begin{prop} \label{adm2}
The Higgs pair $(B,\Phi)$ is admissible if and only if \linebreak
$L^2_E{\rm -ker}{\cal D}=\{0\}$.
\end{prop}
\pf
Given a section $s\in L^2_2(V\otimes S^+)$, the Weitzenb\"ock
formula with respect to the Euclidean metric on the punctured torus
is given by:
\begin{eqnarray*}
(\del_B^*\del_B+\del_B \del_B^*)s & = &\nabla^*_B\nabla_Bs+F_Bs \ = \
\nabla^*_B\nabla_Bs-[\Phi,\Phi^*]s \\
\Rightarrow \ \ \ \nabla^*_B\nabla_Bs & = &
(\del_B^*\del_B+\del_B\del_B^*+\Phi\Phi^*+\Phi^*\Phi)s \\
& = & \left\{ (\del_B+\Phi)(\del_B^*+\Phi^*)+(\del_B^*+\Phi^*)(\del_B+\Phi) \right\}s \\
& = &{\cal D}^*{\cal D}s
\end{eqnarray*}
and integrating by parts, we get:
$$ ||{\cal D}s||^2_{L^2_E} \ = \ ||\nabla_Bs||^2_{L^2_E} $$
Thus, if $B$ is admissible, then the $L^2_E$-kernel of $\cal D$ must
vanish. The converse statement is also clear.
\pfend

In other words, the above proposition implies that the $L^2_E$-cohomologies
of orders 0 and 2 of the complex:
\begin{equation} \label{tmp.cpx}
{\cal C} \ : \ 0 \seta L^2_{2,E}(\Lambda^0V) \stackrel{\Phi+\del_B}{\longrightarrow}
  L^2_{1,E}(\Lambda^{1,0}V\oplus\Lambda^{0,1}V) \stackrel{\del_B+\Phi}{\longrightarrow}
  L^2_E(\Lambda^{1,1}V) \seta 0
\end{equation}
must vanish. On the other hand, since the $L^2$-norm for 1-forms is
conformally invariant, so the $L^2$-cohomology $H^1({\cal C})$ does
not depend on the metric itself, only on its conformal class.

Motivated by a result of Biquard (theorem 12.1 in \cite{B3}) we
will see how one can identify $H^1({\cal C})$ in terms of certain
hypercohomology vector spaces which we now introduce.

Let $\vv\seta\dual$ be the extended holomorphic vector bundle mentioned
above. Recall that if $\ksi_0$ is not an element of order 2 then the
residue of the Higgs field $\Phi$ at $\as$ is a $k\times k$ matrix of
rank 1. Therefore, if $s$ is a local holomorphic section on a neighbourhood
of $\as$, $\Phi(s)$ has at most a simple pole at $\as$ and its residue
has the form $(*,0,\ldots,0)$ on some suitable trivialisation.

Similarly, if $\ksi_0$ is an element of order 2, $\Phi(s)$ has at most
a simple pole at $\as$ and its residue has the form $(*,*,0,\ldots,0)$
on some suitable trivialisation.

This local discussion motivates the definition of a sheaf $\pp_{\as}$
such that, given an open cover $\{U_\alpha\}$ of $\dual$:
\begin{itemize}
\item $\pp_{\as}(U_\alpha)=\oo_{\dual}(\vv)(U_\alpha)$, if $\as\notin U_\alpha$;
\item $\pp_{\as}(U_\alpha)=\{$meromorphic sections of $U_\alpha\seta U_\alpha\times\cpx^k$
which have at most a simple pole at $\as$ with residue lying either along a
2-dimensional subspace of $\cpx^k$ if $\ksi_0$ has order 2, or along a 1-dimensional
subspace of $\cpx^k$ otherwise$\}$, if $\as\in U_\alpha$.
\end{itemize}
It is easy to see that such $\pp_{\as}$ is a coherent sheaf. To simplify notation,
we drop the subscript $\as$ out.

Hence, $\Phi$ can be regarded as the map of sheaves:
\begin{equation} \label{hypcpx}
\Phi : \vv \seta \pp\otimes K_{\dual}
\end{equation}
Seen as a two-term complex of sheaves, the map (\ref{defhiggs}) induces
an exact sequences of hypercohomology vector spaces (see for
example \cite{Bo}, section 3.1) parametrised by $(z,w)\in\torus$, namely:
\begin{equation} \label{hcoho} \begin{array}{cccccccc}
0 & \seta & {\Bbb H}^0(\dual,\Phi) & \seta & H^0(\dual,\vv) &
    \stackrel{\Phi}{\seta} & H^0(\dual,\pp\otimes K_{\dual}) & \seta \\
  & \seta & {\Bbb H}^1(\dual,\Phi) & \seta & H^1(\dual,\vv) &
    \stackrel{\Phi}{\seta} & H^1(\dual,\pp\otimes K_{\dual}) & \seta \\
  & \seta & {\Bbb H}^2(\dual,\Phi) & \seta & 0 & & &
\end{array} \end{equation}

It is easy to see that:
\begin{eqnarray*}
{\Bbb H}^0(\dual,\Phi) & = & {\rm ker}\left\{ H^0(\dual,\vv)
\stackrel{\Phi}{\seta} H^0(\dual,\pp\otimes K_{\dual}) \right\} \\
{\Bbb H}^2(\dual,\Phi) & = & {\rm coker}\left\{ H^1(\dual,\vv)
\stackrel{\Phi}{\seta} H^1(\dual,\pp\otimes K_{\dual}) \right\}
\end{eqnarray*}
and admissibility implies that the right-hand side must vanish:
restricted to $\dual\setminus\{\as\}$, a section there would
give a section in the kernel of ${\cal D}$.
Therefore, the dimension of ${\Bbb H}^1$ is equal to
$\chi(\pp\otimes K_{\dual})-\chi(\vv)$.

To compute this number, note that there is also a natural map \linebreak
$\vv\stackrel{\iota}{\seta}\pp$ defined as the {\em local inclusion}
of holomorphic local sections (elements of $\oo_{\dual}(\vv)(U_\alpha)$),
into the meromorphic ones (elements of $\pp(U_\alpha)$). It fits
into the following sequence of sheaves:
\begin{equation} \label{id} \begin{array}{cl}
0\seta \vv \stackrel{\iota}{\seta} \pp
\stackrel{res_{\ksi_0}}{\longrightarrow} \rr_{\ksi_0}\seta 0
& {\rm if}\ \ksi_0\ {\rm has\ order\ 2,} \\
0\seta \vv \stackrel{\iota}{\seta} \pp
\stackrel{res_{\as}}{\longrightarrow} \rr_{\as}\seta 0 & {\rm otherwise}
\end{array} \end{equation}
where $\rr_{\ksi_0}$ is the skyscraper sheaf supported at $\ksi_0$
and stalk isomorphic to $\cpx^2$ and $\rr_{\as}$ is the skyscraper
sheaf supported at $\as$ and stalks isomorphic to $\cpx$. Since
$\chi(\rr_{\as})=2$, we conclude that ${\Bbb H}^1$ is a 2-dimensional
complex vector space.

\begin{prop} \label{ker/hcoh}
The hypercohomology induced by the map of sheaves (\ref{hypcpx})
coincides with the $L^2_P$-cohomology of the complex (\ref{tmp.cpx}).
\end{prop}

In particular, we have identifications:
$$ {\Bbb H}^1(\dual,\Phi) \ \equiv \ L^2_P{\rm -cohomology}\ H^1({\cal C})
   \ \equiv \ L^2_E{\rm -cohomology}\ H^1({\cal C}) $$
Furthermore, note also that the $L^2$-cohomology of 1-forms with
respect to the Euclidean metric is a 2 dimensional complex vector
spaces.

\pf The hypercohomology defined by the map (\ref{hypcpx}) is given by the
total cohomology of the double complex:
\begin{eqnarray*}
\Lambda^0\vv & \stackrel{\Phi}{\seta} & \Lambda^{1,0}\pp \\
\del\ \downarrow & & \downarrow\ \del \\
\Lambda^{0,1}\vv & \stackrel{\Phi}{\seta} & \Lambda^{1,0}\pp
\end{eqnarray*}
which in turns is just the cohomology of the complex:
$$ 0\seta \Lambda^0\vv \stackrel{\Phi+\del}{\seta} \Lambda^{1,0}\pp \oplus
\Lambda^{0,1}\vv \stackrel{\del+\Phi}{\seta} \Lambda^{1,0}\pp \seta 0 $$
Now restricting the complex above to the punctured torus
$\dual\setminus\{\as\}$, we get:
$$ 0 \seta \Lambda^0V \stackrel{\Phi+\del_B}{\seta} \Lambda^1V
\stackrel{\del_B+\Phi}{\seta} \Lambda^2V \seta0 $$
which is, of course, the complex ${\cal C}$.

So, let $s$ be a section of $\Lambda^{1,0}\pp\oplus\Lambda^{0,1}\vv$
defining a class in ${\Bbb H}^1(\dual,\Phi)$. Thus, restricting $s$
to $\dual\setminus\{\as\}$ yields a section $s_r$ of $L^2(\Lambda^1V)$
defining a class in $H^1({\cal C})$.

Such {\em restriction map} is clearly a well-defined map:
\begin{eqnarray*}
R:{\Bbb H}^1(\dual,\Phi) & \seta & H^1({\cal C}) \\
<s> & \seta & <s_r>
\end{eqnarray*}
We claim that it is also injective. Indeed, suppose that $s_r$
represents the zero class, i.e. there is $t\in L^2_2(\Lambda^0V)$
such that $s_r=(\del_B+\Phi)t$. However, $L^2_2\hookrightarrow C^0$
is a bounded inclusion in real dimension 2. Thus, $h(t,t)$ must be
bounded at the punctures $\as$, and $t$ must be itself bounded
along the kernel of the residues of $\Phi$. On the other hand, the
hermitian metric degenerates along the image of the residues of
$\Phi$, so $t$ might be singular on this direction. However, $h\sim
O(r^{1\pm\alpha})$ is a holomrophic trivialisation, so that $t\sim
O(r^{-\frac{1}{2}(1\pm\alpha)})$. But then the derivatives of $t$
will not be square integrable, contradicting our hypothesis that
$t$ belongs to $L^2_2$. So $t$ must be bounded at $\as$.

This implies that $t\in L^2_2(\Lambda^0\vv)$ also with respect to
the $h'$ metric, so that $s_r$ is indeed the restriction of a
section representing the zero class in ${\Bbb H}^1(\dual,\Phi)$.

Finally, to show that $R$ is an isomorphism, it is enough by admissibility
to argue that the $L^2$ index of the complex ${\cal C}$ is $-2$.

It was shown by Biquard (theorem 5.1 in \cite{B3}) the laplacian
associated to the complex ${\cal C}$ is Fredholm when acting
between $L^2_P$ sections. This implies that ${\cal D}$ is also
Fredholm. Its index can be computed via Gromov-Lawson's relative
index theorem, and it coincides with the index of the Dirac
operator on $\vv$:
$$ \ind({\cal D})=\ind(\del_B-\del_B^*)={\rm deg}\vv=-2 $$
as desired
\pfend


\paragraph{Constructing the transformed bundle.}
We are finally in a position to construct a vector bundle with
connection over $\torus$ out of a Higgs pair $(B,\Phi)$. Let
$L_z\seta\dual\setminus\{\as\}$ be a flat line bundle as in
section \ref{poin}, with its natural connection $\omega_z$, and
form the tensor product $V(z)=V\otimes L_z$. The connection $B$ can
be tensored with $\omega_z$ to obtain another connection that we
denote by $B_z$.

Let $i:V(z)\seta V(z)$ be the identity bundle automorphism and define
$\Phi_{w}=\Phi-w\cdot i$, where $w$ is a complex number. It is easy
to see that $(B_z,\Phi_w)$ is still an admissible Higgs pair, for all
$(z,w)\in\torus$.

Now, consider the following continuous family of Dirac-type operators:
\begin{equation} \label{prt}
{\cal D}_{(z,w)}=(\del_{B_z}+\Phi_w)-(\del_{B_z}+\Phi_w)^*
\end{equation}
From proposition \ref{adm2}, we have that ${\rm ker}{\cal D}_{(z,w)}=\{0\}$
for all $(z,w)\in\torus$, and since its index remains invariant under this
continuous deformation, we conclude that ${\rm ker}{\cal D}_{(z,w)}^*$ has
constant dimension equal to 2.

Consider now the trivial Hilbert bundle $\check{H}\seta\torus$
with fibres given by $L^2(V(z)\otimes S^-)$. It follows that
$\check{E}_{(z,w)}={\rm ker}{\cal D}_{(z,w)}^*$ forms a vector
sub-bundle $\check{E}\stackrel{i}{\hookrightarrow}\check{H}$
of rank 2. Furthermore \cite{DK}, $\check{E}$ is also equipped
with an hermitian metric, induced from $\check{H}$, which we denote
by  $H$; and an unitary connection $\check{A}$, so-called {\em inverse
transformed connection}, defined as follows:
\begin{equation} \label{invconn}
\nabla_{\check{A}}=P \circ d \circ i
\end{equation}
where $d$ means differentiation with respect to $(z,w)$ on the
trivial Hilbert bundle and $P$ is the fibrewise orthogonal
projection $P:L^2(V(z)\otimes S^-)\seta {\rm ker}{\cal D}_{(z,w)}^*$,
with respect to the natural hermitian metric on the Hilbert bundle.
Clearly, $\check{A}$ defined on (\ref{invconn}) is unitary.

Note also that the hermitian metric in $\check{H}$ is actually
conformally invariant with respect to the choice of metric in
$\dual\setminus\{\as\}$, since the inner product in
$L^2(V(z)\otimes S^-)$ is. Therefore, the induced hermitian
metric $H$ in $\check{E}$ depends only on the conformal class
of the metric on the punctured dual torus.

Finally, it is not difficult to see that gauge equivalent Higgs
pairs $(B,\Phi)$ and $(B',\Phi')$ will produce gauge equivalent
instantons $\check{A}$ and $\check{A}'$. The dependence of
$\check{A}$ on the Higgs pair $(B,\Phi)$ is contained on the
$L^2$-projection operator $P$, i.e. on the 2 linearly independent
solutions of ${\cal D}^*_{(z,w)}\psi=0$. Gauge equivalence of $(B,\Phi)$
and $(B',\Phi')$ gives an automorphism of the transformed bundle
$\check{E}$, in other words, a gauge equivalence between $\check{A}$
and $\check{A}'$.


\paragraph{Anti-self-duality.}
In order to complete the inverse transform we must check if the connection
$\check{A}$ is anti-self-dual and if it is extensible. We now consider
the first problem; the second will be treated in the following section.

\begin{prop} \label{asd.inv}
$\check{A}$ is irreducible and anti-self-dual.
\end{prop}

\pf Irreducibility follows from proposition \ref{chch}. Since $\check{A}$
is an unitary connection, we only have to verify that the component
of $F_{\check{A}}$ along the K\"ahler class $\kappa$ of $\torus$
vanishes. Calculations are similar to those in the proof of theorem
\ref{soln}. Let $\{\psi_1,\psi_2\}$ be a local orthonormal frame
for $\check{E}$, with respect to the hermitian metric induced from
$\check{H}$. Fix some $(z,w)\in\torus$ so that, as a section of
$V(z)\otimes S^-\seta\dual$, we have
$\psi_i=\psi_i(\ksi;z,w)\in{\rm ker}{\cal D}_{(z,w)}^*$.

In this trivialisation, the matrix elements of the
curvature $F_{\check{A}}$ can then be written as follows:
\begin{eqnarray*}
(F_{\check{A}})_{ij} & = & \langle \psi_j,\nabla_{\check{A}}\nabla_{\check{A}}\psi_i \rangle
\ = \ \langle \psi_j,d\circ P\circ d\psi_i \rangle \ = \\
& = & \langle {\cal D}_{(z,w)}^*(d\psi_j),G_{(z,w)}{\cal D}_{(z,w)}^*(d\psi_j) \rangle
\end{eqnarray*}
where the inner product is taken in $L^2(V(z)\otimes S^-)$, integrating
out the $\ksi$ coordinate. Recall also that this is conformally invariant
with respect to the choice of metric on $\dual\setminus\{\as\}$.

Moreover, $G_{(z,w)}$ is the Green's operator for
${\cal D}_{(z,w)}^*{\cal D}_{(z,w)}$. Note that:
$$ [{\cal D}_{(z,w)}^*,d]\psi_i=\Omega^\prime\cdot \psi_i $$
where $\Omega^\prime=(idz_1+dw_1)\wedge d\ksi_1+(idz_2+dw_2)\wedge
d\ksi_2$ and ``$\cdot$'' denotes Clifford multiplication; compare with
(\ref{commut}). So,
$$ \kappa\llcorner(F_{\check{A}})_{ij} \ = \
   \langle \psi_j,\kappa\llcorner(\Omega'\wedge\Omega')\cdot G_{(z,w)}\psi_i \rangle \ = \ 0 $$
\pfend


\paragraph{Asymptotic estimate of the curvature.}

We must now work towards establishing that the inverse transformed
instanton connection $\check{A}$ satisfies the extensibility
conditions described in the introduction. We start with the
following result:

\begin{prop} $|F_A|\sim O(r^{-2})$. \end{prop}

\pf As in proposition \ref{asd.inv}, the matrix elements of the
curvature, in the local frame $\{\psi_i\}$, are given by:
$$ (F_{\check{A}})_{ij} = \langle (\Omega'\wedge\Omega')\cdot\psi_j,G_{(z,w)}\psi_i \rangle $$
Therefore, it is easy to see that the asymptotic behaviour of
$|(F_{\check{A}})_{ij}|$ depends only on the behaviour of the operator
norm $||G_{(z,w)}||$ for large $|w|$.

We can estimate $||G_{(z,w)}||$ by looking for a lower bound for the
eigenvalues of the associated laplacian acting on $V(z)\otimes S^-$:
\begin{equation} \label{zwlap}
{\cal D}_{(z,w)}{\cal D}_{(z,w)}^* \ = \
{\cal D}_{z}{\cal D}_{z}^* - w\phi^* - \overline{w}\phi+|w|^2
\end{equation}
where ${\cal D}_{z}={\cal D}_{(z,w=0)}$ and $\Phi=\phi d\ksi$, with
$\phi\in{\rm End}V$; $\phi^*$ denotes the adjoint (conjugate transpose)
endomorphism.

In other words, we want to find a lower bound for the following expression:
\begin{eqnarray}
& \left| \langle ({\cal D}_{z}{\cal D}_{z}^* + |w|^2)s,s \rangle -
         \langle(w\phi^* + \overline{w}\phi)s,s \rangle \right| \geq & \nonumber \\
& \geq \left| \  \langle({\cal D}_{z}{\cal D}_{z}^* + |w|^2)s,s  \rangle \ - \
            | \langle (w\phi^* + \overline{w}\phi)s,s \rangle | \ \right| \label{fin} &
\end{eqnarray}
for $s\in L^2_1(V\otimes S^-)$ of unit norm.

For the first term in the second line, it is easy to see that:
\begin{equation} \label{w2term}
|  \langle ({\cal D}_{z}{\cal D}_{z}^*+|w|^2)s,s  \rangle | \ = \
||{\cal D}_{z}^*s||^2 + |w|^2\cdot||s||^2 \ = \ c_1 + |w|^2
\end{equation}
for some non-zero constant $c_1=||{\cal D}_{z}^*||^2$ depending only
on $z\in T$.

The second term in (\ref{fin}) is more problematic; first note
that:
$$ |  \langle (w\phi^* + \overline{w}\phi)s,s  \rangle | \ \leq \
   |w| \cdot \left( | \langle \phi(s),s \rangle | + | \langle \phi^*(s),s \rangle | \right) $$
In a small neighbourhood $D_0$ of each singularity $\as$, we have:
\begin{eqnarray*}
\langle \phi(s),s \rangle_{L^2(D_0)} & = &
\int_{D_0} \langle \frac{\phi_0(s)}{\ksi},s \rangle rdrd\theta+
\left( \begin{array}{c} {\rm regular} \\ {\rm terms} \end{array} \right) \\
& \sim & \int_{D_0} \frac{|\phi_0|}{r}\cdot|s|^2 rdrd\theta
+ \left(\begin{array}{c} {\rm regular} \\ {\rm terms} \end{array}\right)
\end{eqnarray*}
Let $1<p<2$; using H\"older inequality, we obtain:
\begin{eqnarray*}
\int_{D_0} \frac{|\phi_0|}{\ksi}\cdot|s|^2 & \leq &
\left\{ \int_{D_0} \left( \frac{|\phi_0|}{r} \right)^p rdrd\theta \right\}^{1/p}
\left\{ \int_{D_0} |s|^{2q} \right\}^{1/q} \\
& \leq & c \cdot ||s||^2_{L^{2q}}
\end{eqnarray*}
where $q=\frac{p}{p-1}$, and for some real constant $c$.

Since $2q>4$, the Sobolev embedding theorem tells us that
$L^2_1\hookrightarrow L^{2q}$ is a bounded inclusion (in real
dimension 2). In other words, there is a constant $C$ depending
only on $q$ such that $||s||_{L^{2q}} \leq C \cdot ||s||_{L^2_1}$.
Thus, arguing similarly for the $\langle \phi^*(s),s \rangle$ term,
we conclude that:
$$ |  \langle (w\phi^* + \overline{w}\phi)s,s  \rangle | \ \leq \ c_2 \cdot |w| $$
where $c_2$ is a real constant depending neither on $z$ nor on $w$, but only on
the Higgs field itself and on the choice of $p$.

Putting everything together, we have:
$$ \left| \langle ({\cal D}_{z}{\cal D}_{z}^* - w\phi^* - \overline{w}\phi + |w|^2)s,
   s \rangle \right| \geq \left| |w|^2-c_2|w|+c_1 \right| $$
so that
$$ \lim_{|w|\seta\infty}|w|^2\cdot||G_{(z,w)}||<1 $$
and the statement follows. \pfend


\paragraph{Monad description.}
As in the definition of the dual bundle, $\check{E}$ also admits
a monad type description. More precisely, once a metric is chosen,
the family of Dirac operators (\ref{prt}) can be unfolded into the
following family of elliptic complexes ${\cal C}(z,w)$:
\footnotesize \baselineskip18pt
\begin{equation} \label{mnd}
0 \seta L^2_{2,E}(\Lambda^0V(z)) \stackrel{\Phi_w+\del_{B_z}}{\longrightarrow}
L^2_{1,E}(\Lambda^{1,0}V(z)\oplus\Lambda^{0,1}V(z))
\stackrel{\del_{B_z}+\Phi_w}{\longrightarrow} L^2_E(\Lambda^{1,1}V(z)) \seta 0
\end{equation}
\normalsize \baselineskip18pt

Admissibility implies that $H^0({\cal C}(z,w))$ and $H^2({\cal C}(z,w))$ must
vanish, and $H^1({\cal C}(z,w))$ coincides with $L^2_E{\rm -ker}{\cal D}^*_{(z,w)}$. As
$(z,w)$ sweeps out $\torus$, $H^1({\cal C}(z,w))$ forms a rank 2 holomorphic vector
bundle with a natural hermitian metric and a compatible unitary connection
$A$, equivalent to the ones defined as above; see \cite{DK}.

We now pass to the holomorphic description of the inverse transform.
It will allow us to compute the instanton number and the asymptotic
state of inverse transformed connection $\check{A}$.


\section{Holomorphic description} \label{invholo}

Motivated by section 2.1, one can expect to find a holomorphic
vector bundle $\check{\ee}\seta\tproj$ which extends
$(\check{E},\del_{\check{A}})$. The idea is to find a suitable
perturbation of the Higgs field $\Phi$ for which $w=\infty$ makes
sense.

As above, the torus parameter $z\in T$ simply twists the holomorphic
bundle $\vv\seta\dual$. We denote:
\begin{equation} \begin{array}{ccc}
\vv(z)=\vv\otimes L_z & \ \ \ & \pp(z)=\pp\otimes L_z
\end{array} \end{equation}
Since $\Phi\in H^0(\dual,{\rm Hom}(\vv,\pp)\otimes K_{\dual})$,
tensoring both sides of (\ref{hypcpx}) by the line bundle $L_z$
does not alter the sheaf homomorphism $\Phi$, so we have the
family of maps:
$$ \Phi : \vv(z) \seta \pp(z)\otimes K_{\dual} $$
parametrised by $z\in T$.

To define the perturbation $\Phi_w$, recall that, regarding
$\proj=\cpx\cup\{\infty\}$, we can fix two holomorphic sections
$s_0,s_\infty\in H^0(\proj,\oo_{\proj}(1))$ such that $s_0$ vanishes
at $0\in\cpx$ and $s_\infty$ vanishes at the point added at infinity.
In homogeneous coordinates $\{(w_1,w_2)\in\cpx^2|w_2\neq0\}$ and
$\{(w_1,w_2)\in\cpx^2|w_1\neq0\}$, we have that, respectively
($w=w_1/w_2$):
\begin{eqnarray*}
s_0(w)=w & \ \ \ & s_0(w)=1 \\
s_\infty(w)=1 & \ \ \ & s_\infty(w)=\frac{1}{w}
\end{eqnarray*}

Consider now the map of sheaves parametrised by pairs $(z,w)\in\tproj$:
\begin{eqnarray}
& \Phi_w : \vv(z) \seta \pp(z)\otimes K_{\dual} & \nonumber \\
& \Phi_w = s_\infty(w)\cdot\Phi-s_0(w)\cdot\iota\cdot d\ksi & \label{defhiggs}
\end{eqnarray}
Clearly, on $\proj\setminus\{\infty\}=\cpx$ this is just
$\Phi_w=\Phi-w\cdot\iota$, the same perturbation we defined before.
On the other hand, if $w=\infty$, then $\Phi_\infty = - \iota\cdot d\ksi$

The hypercohomology vector spaces ${\Bbb H}^0(\dual,\Phi_w)$ and
${\Bbb H}^2(\dual,\Phi_w)$ of the two-term complex (\ref{defhiggs})
must vanish by admissibility. On the other hand,
${\Bbb H}^1(\dual,\Phi_w)$ also makes sense for $\infty\in\proj$,
the inverse transformed bundle with connection $(\check{E},\check{A})$
admits a compatible holomorphic extension to a bundle $\check{\ee}\seta\tproj$
(in the sense of section 2.1.2), with fibres given by
$\check{\ee}_{(z,w)}={\Bbb H}^1(\dual,\Phi_w)$, as desired.

Equivalently, $\check{\ee}$ can be seen as the hermitian holomorphic vector
bundle induced by the monad
\begin{equation} \label{op.db.cpx}
0\seta \Lambda^0\vv \stackrel{\Phi+\del}{\seta} \Lambda^{1,0}\pp \oplus
\Lambda^{0,1}\vv \stackrel{\del+\Phi}{\seta} \Lambda^{1,0}\pp \seta 0
\end{equation}

Consider the metric $H'$ induced from the monad (\ref{op.db.cpx}) above,
while $H$ is induced from the monad (\ref{mnd}). Now, $H$ is bounded above
by $H'$ because the hermitian  metric $h$ on the bundle $V$ in (\ref{mnd})
is bounded above by the metric $h'$ on the bundle $\vv$ in (\ref{op.db.cpx}).

Let us now compute the Chern character of $\check{\ee}$.

\begin{lem} \label{chch}
Using the notation of section \ref{poin},
$ch(\check{\ee}) = 2 - k\cdot t\wedge p$.
\end{lem}

\pf
The exact sequence:
\begin{equation} \label{hcoho2} \begin{array}{ccl}
0 & \seta & H^0(\dual,\vv(z)) \stackrel{\Phi_w}{\seta}
            H^0(\dual,\pp(z)\otimes K_{\dual}) \seta {\Bbb H}^1(\dual,(z,w)) \seta \\
  & \seta & H^1(\dual,\vv(z)) \stackrel{\Phi_w}{\seta}
            H^1(\dual,\pp(z)\otimes K_{\dual}) \seta 0
\end{array} \end{equation}
induces a sequence of coherent sheaves over $\torus$, with stalks over $(z,w)$
given by the above cohomology groups:
\begin{equation} \label{hcoho3} \begin{array}{ccl}
0 & \seta & {\cal H}^0(\dual,\vv(z)) \stackrel{\Phi_w}{\seta}
            {\cal H}^0(\dual,\pp(z)\otimes K_{\dual}) \seta \check{\ee} \seta \\
  & \seta & {\cal H}^1(\dual,\vv(z)) \stackrel{\Phi_w}{\seta}
            {\cal H}^1(\dual,\pp(z)\otimes K_{\dual}) \seta 0
\end{array} \end{equation}
In this way, the Chern character of $\check{\ee}$ will then be given
by the alternating sum of the Chern characters of these sheaves, which can
be computed via the usual Grothendieck-Riemann-Roch for families.

Consider the bundle ${\bf G}_1\seta\tproj\times\dual$ given by
${\bf G}_1=p_3^*\vv\otimes p_{13}^*{\bf P}$. Clearly,
${\bf G}_1|_{(z,w)\times\dual}=\vv(z)$, so that:
\begin{equation} \label{yyy}
ch({\cal H}^0(\dual,\vv(z)))-ch({\cal H}^1(\dual,\vv(z)))=
ch({\bf G}_1)td(\dual)/[\dual]
\end{equation}
where $t$ is the generator of $H^2(T)$, as in section \ref{poin}.

Now consider the sheaf
${\bf G}_2=p_3^*{\cal P}\otimes p_{13}^*{\bf P}\otimes p_2^*\oo_{\proj}(1)$.
The twisting by $\oo_{\proj}(1)$ accounts for the multiplication by the
section $s_0\in H^0(\proj,\oo_{\proj}(1))$ contained in $\Phi_w$. As
above, ${\bf G}_1|_{(z,w)\times\dual}=\pp(z)$, and we have:
\begin{equation} \label{yyz}
ch({\cal H}^0(\dual,\pp(z)\otimes K_{\dual})) -
ch({\cal H}^1(\dual,\pp(z)\otimes K_{\dual})) = ch({\bf G}_2)td(\dual)/[\dual]
\end{equation}
where $p$ is the generator of $H^2(\proj)$, as in section \ref{poin}.

Therefore:
\begin{eqnarray*}
ch(\check{\ee}) & = & ( \ref{yyz}) - (\ref{yyy}) \ = \\
                & = & \left( c_1(\pp) - c_1(\vv) + c_1(\pp)\wedge p -
                             \frac{k}{2}c_1({\bf P})^2\wedge p \right) / [\dual] \ = \\
                & = & \chi(\pp) - {\rm deg}\vv + \chi(\pp)\cdot p - k\cdot t\wedge p \ = \
                      2 - k\cdot t\wedge p
\end{eqnarray*}
as desired.
\pfend

The next lemma determines the asymptotic state of the inverse
transformed connection.

\begin{lem} \label{infty}
$\check{\ee}|_{T_{\infty}}\equiv L_{\ksi_0}\oplus L_{-\ksi_0}$
\end{lem}

\pf Substituting $w=\infty\in\proj$, we get from (\ref{defhiggs})
that $\Phi_\infty=\iota\cdot d\ksi$. Therefore, the induced
hypercohomology sequence (\ref{hcoho2}) coincides with the long
exact sequence of cohomology induced by the sheaf sequence (\ref{id}),
which is given by:
\begin{equation} \label{hcoho4} \begin{array}{rcl}
0 & \seta & H^0(\dual,\vv(z)) \stackrel{\Phi_\infty}{\seta}
            H^0(\dual,\pp(z)\otimes K_{\dual}) \seta H^0(\dual,\rr_{\as}(z)) \seta \\
  & \seta & H^1(\dual,\vv(z)) \stackrel{\Phi_\infty}{\seta}
            H^1(\dual,\pp(z)\otimes K_{\dual}) \seta 0
\end{array} \end{equation}
Hence, ${\Bbb H}^1(\dual,(z,\infty))=H^0(\dual,\rr_{\as}(z))$.
The right hand side is canonically identified with
$(L_z)_{\ksi_0}\oplus (L_z)_{-\ksi_0}$, where by $(L_z)_{\ksi_0}$
we mean the fibre of $L_z\seta\dual$ over the point $\ksi_0\in\dual$.

On the other hand, $(L_z)_{\ksi_0}={\bf P}_{(z,\ksi_0)}=(L_{\ksi_0})_z$,
where ${\bf P}\seta T\times\dual$ is the Poincar\'e line bundle. Thus,
the bundle over $T_\infty$ with fibres given by $H^0(\dual,\rr_{\as}(z))$
is isomorphic to $L_{\ksi_0}\oplus L_{-\ksi_0}$, as we wished to prove.
\pfend

Finally, we argue that the determinant bundle of $\check{\ee}$ is
trivial, so that $\check{A}$ is indeed an $SU(2)$ instanton. Note
that ${\rm det}\check{\ee}$ is a line bundle with vanishing first
Chern class, so it must be the pull back of a flat line bundle
$L_\ksi\seta T$. But ${\rm det}\check{\ee}|_{T_\infty}=\underline{\cpx}$,
hence ${\rm det}\check{\ee}$ must be holomorphically trivial, as desired.

Thus, we conclude that $\check{A}\in\cala_{(k,\ksi_0)}$.

\paragraph{Final remark.}
Summing up the work done in this section, we established a map from
the set of equivalence classes of Higgs pairs $(B,\Phi)$ on a vector
bundle $V\seta\dual\setminus\{\as\}$ of rank $k$, such that $\Phi$
has simple poles at $\as$ with a residue of rank 1 or 2 (depending
on the order of $\ksi_0$), to the set of gauge equivalence classes
of unitary instanton connections $\check{A}\in\cala_{(k,\ksi_0)}$
on a rank 2 bundle $\check{E}\seta\torus$.

Note however that this procedure depends on the connection $B$ only
through the holomorphic structure it induces in $V$. Of course, this
piece of information is fully contained in the extended holomorphic
bundle $\vv\seta\dual$.

Finally, the abelian group $\torus\times S^1$ acts on the
set of Higgs bundles as follows:
\begin{equation} \label{inv.act.sec}
(x,y,\gamma)\cdot(\vv,\Phi) \  \mapsto \
(\vv\otimes L_x,e^{i\gamma}\cdot\Phi+y\cdot I)
\end{equation}
and this clearly corresponds to the action of $\torus\times S^1$ on
the set of extensible instanton connections via pullback, see section
\ref{actsec}.

Note also that $\dual$ does not act on the moduli of Higgs bundles via
pullback: since the singularities are fixed at $\as$, we are not
allowed to make translations on $\dual$.

\chapter{Completing the proof of theorem \ref{nahmthm}} \label{conc} 

We finally arrived to the final stage of the proof of the Nahm
transform theorem. More precisely, there are still two issues to be
addressed: first, we must show that the Higgs pairs initially
constructed from an instanton connection are indeed admissible;
second, we need to verify that $(\check{E},\check{A})$ is
equivalent to the original data $(E,A)$.

First, consider the six-dimensional manifold $\torus\times(\dual\setminus\{\as\})$.
To shorten notation, we denote $M_\ksi=\torus\times\{\ksi\}$ and
$\dual_{(z,w)}=\{z\}\times\{w\}\times(\dual\setminus\{\as\})$.

Now take the bundle $\calg=p_{12}^*E\otimes p_{13}^*\mathbf{P}$
over $\torus\times(\dual\setminus\{\as\})$; note that
$\calg|_{M_\ksi}\equiv E(\ksi)$ and $\calg|_{\dual_{(z,w)}} \equiv
\underline{E_{(z,w)}}\otimes L_z$, where $\underline{E_{(z,w)}}$
denotes a trivial rank 2 bundle over $\dual\setminus\{\as\}$ with
the fibres canonically identified with the vector space
$E_{(z,w)}$.

$\cal G$ is clearly holomorphic; we denote by $\del_M$ the action
of the associated Dolbeault operator along the $\tproj$ direction,
and by $\del_{\dual}$ its action along the $\dual$ direction. In
particular, $\del_M|_{M_\ksi} \equiv \del_{A_\ksi}$.

Let ${\bf C}^{p,q}=\Lambda^{0,p}_{\torus}({\cal G})\otimes\Lambda^q_{\dual}({\cal G})$;
in other words, ${\bf C}^{p,q}$ consists of the $(p+q)$-forms over
$\torus\times(\dual\setminus\{\as\})$ with values in $\cal G$ spanned by
forms of the shape:
\begin{equation} \label{cpq} \begin{array}{c}
s(z,w,\ksi)d\overline{z}_{i_1}d\overline{w}_{i_2}d\ksi_{j_1}d\overline{\ksi}_{j_2}, \\
i_1,i_2,j_1,j_2\in\{0,1\} \ \ \ {\rm and} \ \ \ i_1+i_2=p,\ j_1+j_2=q
\end{array} \end{equation}
Analytically, we want to regard ${\bf C}^{p,q}$ as the completion of the
set of smooth forms of the shape above with respect to a Sobolev norm
described as follows:
$$ \begin{array}{cl}
\left| s|_{\torus\times\{\ksi\}} \right| \in L^2_q(\Lambda^{2-q}E(\ksi)) &
{\rm for\ each}\ \ksi\in\dual\setminus\{\as\} \\
\left| s|_{\{(z,w)\}\times\dual\setminus\{\as\}} \right| \in L^2_q(\Lambda^{2-q}L_z) &
{\rm for\ each}\ (z,w)\in\torus
\end{array} $$

Now, define the maps:
\begin{equation} \label{this} \begin{array}{rcl}
\mathbf{C}^{p,0} \stackrel{\delta_1}{\seta} & {\bf C}^{p,1} &
\stackrel{\delta_2}{\seta} \mathbf{C}^{p,2} \\
\delta_1(s)=(\del_{\dual}s,-w\cdot s\wedge d\ksi) & &
\delta_2(s_1,s_2)=(\del_{\dual}s_2+w\cdot s_1\wedge d\ksi)
\end{array} \end{equation}
for $(s_1,s_2)\in\Lambda^{0,p}_{\tproj}({\cal G})\otimes
\left( \Lambda^{0,1}_{\dual}({\cal G})\oplus\Lambda^{1,0}_{\dual}({\cal G}) \right)
\equiv {\bf C}(p,1)$. Note that (\ref{this}) does define a complex.

The inversion result will follow from the analysis of the spectral
sequences associated to the following double complex
(for the general theory of spectral sequences and double complexes,
we refer to \cite{BT}):
\begin{equation} \label{dblcpx} \begin{array}{ccccc}
\mathbf{C}^{0,2} & \stackrel{\del_M}{\seta} & \mathbf{C}^{1,2} &
\stackrel{-\del_M}{\seta} & \mathbf{C}^{2,2} \\
\ \ \uparrow \delta_2 & & \ \ \ \uparrow -\delta_2 & & \ \ \uparrow \delta_2 \\
\mathbf{C}^{0,1}& \stackrel{\del_M}{\seta} &
\mathbf{C}^{1,1}& \stackrel{-\del_M}{\seta} & \mathbf{C}^{2,0} \\
\ \ \uparrow \delta_1 & & \ \ \ \uparrow -\delta_1 & & \ \ \uparrow \delta_1 \\
\mathbf{C}^{0,0} & \stackrel{\del_M}{\seta} & \mathbf{C}^{1,0} &
\stackrel{-\del_M}{\seta} & \mathbf{C}^{2,0}
\end{array} \end{equation}
The idea is to compute the total cohomology of the spectral sequence
in the two possible different ways and compare the filtrations of
the total cohomology.

\begin{lem} \label{rows}
By first taking the cohomology of the rows, we obtain
\begin{equation} \begin{array}{cccccc}
          & \ & & 0 & H^2({\cal C}(e,0)) & 0 \\
E^{p,q}_2 & \ & & 0 & H^1({\cal C}(e,0)) & 0 \\
& \ & q\uparrow & 0 & H^0({\cal C}(e,0)) & 0 \\
& \ & & \stackrel{\seta}{p} & &
\end{array} \end{equation}
where $H^i({\cal C}(e,0))$ are the cohomology groups of the complex
(\ref{tmp.cpx}).
\end{lem}

\pf First, note that the rows coincide with the complex (\ref{monad}) of
section \ref{diffgeo}.

Moreover, we can regard elements in $\mathbf{C}^{p,q}$ as $q$-forms over
$\dual$ with values in $L^2_{2-p}(\Lambda^{0,p}_{\torus}{\cal G})$.
To see this, fix some $\ksi'\in\dual$; by (\ref{cpq}),
$s(z,w,\ksi')\in\Lambda^{0,p}{\cal G}|_{M_{\ksi'}}$. So,
by varying $\ksi'$ we get the interpretation above.

This said, it is clear that the first and second columns of
$E^{p,q}_1$ must vanish, since $A$ is irreducible. In the middle column,
we get $q$-forms over $\dual$ with values in ${\rm ker}(\del_M^*-\del_M)$,
which for a fixed $\ksi'$ restricts to ${\rm ker}(D^*_{A_{\ksi'}})$.

Therefore, after taking the cohomologies of the rows, we are left with:
\begin{equation} \label{rows2} \begin{array}{cccccc}
                & \ & & 0 \ \ & L^2(\Lambda^{1,1}V) & \ \ 0 \\
                & \ & &  & \ \ \ \ \ \ \uparrow(\del_B+\Phi) &  \\
{\bf C}^{p,q}_1 & \ & & 0 \ \ & L^2_1(\Lambda^{1,0}V\oplus\Lambda^{0,1}V) &\ \ 0 \\
                & \ & &  & \ \ \ \ \ \ \uparrow(\Phi+\del_B) &  \\
                & \ & q\uparrow & 0 \ \ & L^2_2(\Lambda^0V) &\ \ 0 \\
                & \ & & \stackrel{\seta}{p} & &
\end{array} \end{equation}
But this is just the complex (\ref{tmp.cpx}). The lemma follows after
taking the cohomology of the remaining column.
\pfend


\paragraph{Total cohomology and admissibility.}
Note that, as we pointed out in the beginning of this section, we
still do not know if the Higgs pair $(B,\Phi)$ arising from the
instanton $(E,A)$ is admissible or not, so that $H^0$ and $H^2$
might be nontrivial. The next lemma deals with this problem.

\begin{lem} \label{t.coh}
The only nontrivial cohomology of the total complex is \linebreak
$H^2({\bf C}(p,q))$, which is naturally isomorphic to the fibre
$E_{(e,0)}$.
\end{lem}

In particular, this shows that the Higgs pairs $(B,\Phi)$ obtained
via Nahm transform on instanton connection $A\in\cfg$ are indeed
admissible, by proposition \ref{adm2}.

\pf First note that we can regard an element in ${\bf C}^{p,q}$ as a
$(0,p)$-form over $\torus$ with values in
$\Lambda^{q_1,q_2}_{\dual}(\calg)$. Since $\calg|_{\dual_{(z,w)}}
\equiv \underline{E_{(z,w)}}\otimes L_z$, ${\rm ker}\del_M$ and
${\rm ker}\del_M^*$ are nontrivial only if $z=e$, the identity
element in the group law of $T$. Hence, it is enough to work on a
tubular neighbourhood of $\{e\}\times\proj\times(\dual\setminus\{\as\})$.

More precisely, we define another double complex $({\rm germ}\ {\bf C})^{p,q}$,
consisting of forms defined on arbitrary neighbourhoods of
$\{e\}\times\proj(\dual\setminus\{\as\})$. Then we have a restriction
map ${\bf C}^{p,q}\seta({\rm germ}\ {\bf C})^{p,q}$ commuting with
$\del_M$, $\delta_1$ and $\delta_2$. Such map also induces an
isomorphism between the total cohomologies of ${\bf C}^{p,q}$ and
$({\rm germ}\ {\bf C})^{p,q}$. So we can work with $({\rm germ}\ {\bf C})^{p,q}$
to prove the lemma.

Let $V_e$ be some neighbourhood of $e\in T$. By the Poincar\'e lemma
applied to $\del_T$, we get:
\begin{equation} \label{t.coh2} \begin{array}{ccccccc}
& \ & & & \Lambda^2_{V_e}(\calg) & 0 & 0 \\
& \ & &  & \ \ \ \ \ \ \uparrow & & \\
({\rm germ}\ {\bf C})^{p,q}_1
& \ & & & \Lambda^1_{V_e}(\calg) & 0 & 0 \\
& \ & &  & \ \ \ \ \ \ \uparrow & &  \\
& \ & q\uparrow & & \Lambda^0_{V_e}(\calg) & 0 & 0 \\
& \ & & \stackrel{\seta}{p} & & &
\end{array} \end{equation}
where $V_e$ denotes a tubular neighbourhood of
$N_e=\{e\}\times\proj\times(\dual\setminus\{\as\})$

As in \cite{DK} (see pages 91-92), the complex in the first row is,
after restriction, mapped into a Koszul complex over $N_e$:
$$ \oo_{N_e}(\calg) \stackrel{(w\ \ksi)}{\longrightarrow}
   \oo_{N_e}(\calg)\oplus \oo_{N_e}(\calg)
   \stackrel{(-\ksi,z)}{\seta}\oo_{N_e}(\calg) $$
so that:
\begin{equation} \label{t.coh3} \begin{array}{cccccc}
& \ & & E_{(e,0)} & 0 & 0 \\
({\rm germ}\ {\bf C})^{p,q}_2
& \ & & 0 & 0 & 0 \\
& \ & q\uparrow & 0 & 0 & 0 \\
& \ & & \stackrel{\seta}{p} & &
\end{array} \end{equation}
\pfend

It then follows from lemmas \ref{rows} and \ref{t.coh} that there is a
natural isomorphism of vector spaces
${\cal I}_I:H^1({\cal C}(e,0))\equiv\check{E}_{(e,0)}\seta E_{(e,0)}$,
which in principle may depend on the choice of complex structure $I$ on $\torus$.

\paragraph{Matching $\mathbf{(\check{E},\check{A})}$ with the original data.}
Since the choice of identity element in $T$ and of origin in $\cpx$
is arbitrary, we can extend ${\cal I}_I$ to a bundle isomorphism
$E\seta\check{E}$. More precisely, let $t_{(u,v)}:\torus\seta\torus$
be the translation map $(z,w)\seta(z+u,w+v)$. Clearly, the connection
$t_{(u,v)}^*A$ on the pullback bundle $t_{(u,v)}^*E$ is also irreducible
and $t_{(u,v)}^*E_{(e,0)}\equiv E_{(u,v)}$. Computing the total
cohomology of the double complex (\ref{dblcpx}) associated to the
bundle $t_{(u,v)}^*{\cal G}$ (where $t_{(u,v)}^*$ acts trivially on
$\dual$ coordinate), lemmas \ref{rows} and \ref{t.coh} lead to an isomorphism
of vector spaces $H^1({\cal C}(u,v))\equiv\check{E}_{(u,v)}\seta E_{(u,v)}$.

It is clear from the naturality of the constructions that these
fibre isomorphisms fit together to define a holomorphic bundle
isomorphism \linebreak ${\cal I}_I:E\seta\check{E}$. In particular,
${\cal I}_I$ takes the Dolbeault operator $\del_A$ of the
holomorphic bundle $E\seta\torus$ to the Dolbeault operator
$\del_{\check{A}}$ of the holomorphic bundle
$\check{E}\seta\torus$. It also follows from this observation that
the holomorphic extensions $\ee$ and $\check{\ee}$ must be
isomorphic as holomorphic vector bundles.

However, such fact still does not guarantee that the connections
$A$ and $\check{A}$ are gauge-equivalent. This is accomplished if
we can show that ${\cal I}_I$ is actually {\em independent} of the
choice of complex structure in $\torus$. Therefore, the proof of
the main theorem \ref{nahmthm} is completed by the following
proposition:

\begin{prop}
The bundle map ${\cal I}_I:\check{E}\seta E$ is independent of the
choice of complex structure on $\torus$.
\end{prop}

\pf Again, it is sufficient to consider only the fibre over
$(e,0)$. As in \cite{DK} (p. 94-95), the idea is to present an
explicit description of ${\cal I}_I:\check{E}_{(e,0)}\seta E_{(e,0)}$,
and then show that it is Euclidean invariant.

Let $\alpha\in H^1({\cal C}(e,0))\subset{\bf C}^{1,1}$. To find
${\cal I}_I([\alpha])$ we have to find $\beta\in{\bf C}^{0,2}$ such that
$\del_M\beta=\delta_2\alpha$. A solution to this equation is provided
by the Hodge theory for the $\del_M$ operator:
$$ \beta=G_M(\del_M^*\delta_2\alpha) $$
where $G_M$ denotes the Green's operator for $\del_M^*\del_M$,
which can be regarded fibrewise as the family of Green's operators
$G_{A_\ksi}=G_M|_{M_\ksi}$ parametrised by $\ksi\in(\dual\setminus\{\as\})$.

In principle, $\beta$ depends on the complex structure $I$ via the
operators $\del_M$ and $G_M$. However, by the Weitzenb\"ock formula
applied to the bundle $\calg$, we have:
$$ \del_M^*\del_M=\nabla^*_M\nabla_M $$
Here, $\nabla_M$ is the covariant derivative in the $\torus$ direction
on $\calg$. With this interpretation, $G_M=(\nabla^*_M\nabla_M)^{-1}$
is seen to be independent of the complex structure $I$; in fact, it is
Euclidean invariant.

Now $\beta$ as an element of ${\bf C}^{1,1}$ has the form
$\beta(z,w;\ksi)d\ksi d\overline{\ksi}$, so that the restriction
$r_{(e,0)}(\beta)=\beta|_{\dual_{(e,0)}}$ is a $(1,1)$-form over
$\dual\setminus\{\as\}$ with values in $E_{(e,0)}$. Take its
cohomology class in $H^2(\dual\setminus\{\as\},\underline{\cpx}\otimes\underline{E_{(e,0)}})$,
so that:
$$ {\cal I}_I([\alpha])=\int_{\dual_{(e,0)}}r_{(e,0)}(\beta) $$
which is the desired explicit description. \pfend

This finally completes the proof of the main theorem \ref{nahmthm}.



\chapter{Further Remarks} \label{rks}

We now want to look more closely at a few consequences of the Nahm
transform theorem.

Our first remark concerns the non-emptiness of the moduli space of
doubly-periodic instantons. As we mentioned in the introduction,
singular solutions of Hitchin's equations are quite well studied,
being closely related to the so-called {\em parabolic Higgs
bundles}. In particular, existence of Higgs pairs of the type we want
is determined by some holomorphic data. Model solutions in a
neighbourhood of the singularities were described by Biquard \cite{B4}:
\begin{eqnarray*}
B & = & b\frac{d\xi}{\xi}+b^*\frac{d\overline{\xi}}{\overline{\xi}}\\
\Phi & = & \phi_0\frac{d\xi}{\xi}
\end{eqnarray*}
Every meromorphic Higgs pair with a simple pole approaches this
model solution close enough to the singularities.  These observations
together with our main theorem \ref{nahmthm} guarantees the
existence of doubly-periodic instantons of any given charge and
asymptotic state.

\paragraph{Holomorphic version.}
Now take the bundle $\calg=p_{12}^*\ee\otimes p_{13}^*\mathbf{P}$ over
\linebreak $\tproj\times\dual$ and consider the appropriate double
complex analogous to (\ref{dblcpx}). It is then easy to establish
results identical to lemmas \ref{rows} and \ref{t.coh}. This in turn
leads to a holomorphic bundle isomorphism between $\ee$ and
$\check{\ee}$, as above. Hence, as a by-product of the Nahm
transform theorem, we obtain the following result, which can be seen
as the holomorphic version of theorem \ref{nahmthm}:

\begin{thm} \label{holov}
There is a bijective correspondence between the following objects:
\begin{itemize}
\item holomorphic vector bundles $\ee\seta\tproj$ with
      ${\rm det}\ee=\underline{\cpx}$, $c_2(\ee)=k>0$ and such that
      $\ee|_{T_\infty}=L_{\ksi_0}\oplus L_{-\ksi_0}$;

\item Higgs bundles $(\vv,\Phi)$ consisting of a rank $k$ holomorphic vector
      bundle $\vv\seta\dual$ of degree $-2$ and a Higgs field
      $\Phi$, which is a meromorphic section of ${\rm End}\vv$
      having simple poles at $\as$ with semi-simple residues of rank $\leq2$,
      if $\ksi_0$ has order 2, and rank $\leq1$ otherwise.
\end{itemize}
\end{thm}

\paragraph{Generalisation to higher rank.}
The attentive reader might have noticed that there is nothing
really special about rank two bundles, and that the whole proof
could easily be generalised to higher rank. Indeed, the only point
in choosing the rank two case is to reduce the number of possible
vector bundles over an elliptic curve, and avoid a tedious
case-by-case study throughout the various stages of the proof.

Before we can state the generalisation of the main theorem \ref{nahmthm},
we must review our definitions of asymptotic state and irreducibility.

The restriction of the holomorphic extension $\ee\seta\tproj$ to
the added divisor $T_\infty$ is a flat $SU(n)$ bundle, i.e.
$$ \ee|_{T_\infty}=L_{\ksi_1}\oplus\dots\oplus L_{\ksi_k} $$
$$ {\rm such\ that}\ \ \ \ \ \bigotimes_{l=1}^k L_{\ksi_l}=\oo_{T} $$
In other words, $\ee|_{T_\infty}$ is determined by a set of points
$(\ksi_1,\dots,\ksi_j)\in{\cal J}(T)$ with multiplicities
$(m_1,\dots,m_j)$, and such that $\sum_{l=1}^j m_l\ksi_l=0$. We
call such data the {\em generalised asymptotic state}.

Moreover, we will say that $(E,A)$ is {\em 1-irreducible} if there is
no flat line bundle $E\seta\torus$ such that $E$  admits a splitting
$E'\oplus L$ which is compatible with the connection $A$.

\begin{thm} \label{nahmthmgen}
There is a bijective correspondence between the following objects:
\begin{itemize}
\item gauge equivalence classes of 1-irreducible, extensible $SU(n)$
instantons over $\torus$ with fixed instanton number $k>0$ and generalised
asymptotic state $(\ksi_1,\dots,\ksi_j)$ with multiplicities $(m_1,\dots,m_j)$  and
\item admissible $U(k)$ solutions of the Hitchin's equations over the dual torus
$\dual$, such that the Higgs field has at most simple poles at
$\{\ksi_1,\dots,\ksi_j\}$; moreover, its residue at $\ksi_j$ is
semi-simple and has rank $\leq m_j$.
\end{itemize}
\end{thm}

Of course, the holomorphic version \ref{holov} can be similarly
generalised. Also, the same remark about the possibility of
removing the technical hypothesis on the non-triviality of the
asymptotic states holds.

\paragraph{Extra parameters for Higgs bundles.}
On the Hitchin's equations side of our picture, there are two types
of parameters one generally fixes, namely the eigenvalues of the
residues of the Higgs field $\Phi$ and the limiting holonomy of the
connection $B$ around the singularities (or equivalently, the {\em
parabolic structure}; see also \cite{K} \cite{S}). In the
terminology of Kovalev, such parameters are called {\em commuting
triples}, for they are equivalent to specifying three mutually
commuting matrices in ${\frak u}(k)$.

In our situation however, only the rank of the residue of the Higgs
field is fixed, while its non-zero eigenvalues are free to vary.
However, Tr$(\Phi)$ is a meromorphic 1-form on $\dual$ with poles
at $\as$, and the sum of the residues must vanish. If $\ksi_0$ is
not of order 2, this implies that the unique non-zero eigenvalue of
the residue of $\Phi$ at $\ksi_0$ is minus the unique non-zero
eigenvalue of the residue of $\Phi$ at $-\ksi_0$. If $\ksi_0$ has
order 2, then the sum of the two non-zero eigenvalues of the
residue of $\Phi$ at $\ksi_0$ must vanish. Therefore, the
eigenvalues of the residues of $\Phi$ account for only one complex
degree of freedom, which we denote by $\epsilon$.

The parabolic structure consists of a filtration of $\vv_{\as}$,
the fibre of $\vv$ over the singularities $\as$, plus a choice of
weights $0\leq\alpha_i(\as)<1$. From proposition \ref{higgs}, a
natural choice of filtration would be, generically:
\begin{eqnarray*}
\begin{array}{rcl}
\vv_{\as} = F_1\vv_{\as} \supset & \underbrace{F_2\vv_{\as}} & \supset F_3\vv_{\as} =\{0\} \\
& {\rm dim}=1 &
\end{array}
& & {\rm order}(\ksi_0) \neq 2 \\ & & \\
\begin{array}{rcl}
\vv_{\as} = F_1\vv_{\as} \supset & \underbrace{F_2\vv_{\as}} & \supset F_3\vv_{\as} =\{0\} \\
& {\rm dim}=2 &
\end{array} &
& {\rm order}(\ksi_0) = 2
\end{eqnarray*}
More precisely, from (\ref{sqc2}) we have that in either case:
$$ F_2\vv_{\as} = H^0(T_\infty,\widetilde{\ee}(\as)|_{T_\infty})
   \hookrightarrow H^1(\tproj,\ee(\as)) = F_1\vv_{\as} $$

To complete the parabolic structure, we would have to choose four
weights (two for each parabolic point) in the first case and two
weights in the second case:
$$ 0 \leq \alpha_1(\as) < \alpha_2(\as) < 1 $$
From the point of view of the Higgs pair $(B,\Phi)$, these parameters
can also be interpreted as the rate of growth of local holomorphic sections of
$V\seta\dual\setminus\{\as\}$ near the singular points with respect
to the hermitian metric induced from the Hilbert bundle $\hat{H}$.

If $(\vv,\Phi)$ is $\alpha$-stable in the sense of parabolic Higgs bundles,
then the existence of a meromorphic Higgs pair as above is guaranteed
\cite{S}.

These are natural parameters in the theory of Higgs
bundles, and one would like to interpret them on the instanton side
of the correspondence. However, it is reassuring to know that if two
sets of parameters $(\alpha,\epsilon)$ and $(\alpha',\epsilon')$ are
chosen in generic position, then $\alpha$-stability and $\alpha'$-stability
are in fact equivalent conditions \cite{Nkj}.

\paragraph{Limiting holonomy.}
On the instanton side of the picture there is one further real
parameter that we have not discussed so far: the limiting holonomy
of the instanton connection $A$ around the added divisor
$T_\infty$.

More precisely, write the connection in {\em radial gauge} so that
$$ A = a_xdx + a_ydy + a_\theta d\theta $$
and look at the following initial value problem for a function
$f:S^1\seta SU(2)$:
$$ \frac{df_r}{d\theta}+a_\theta f_r=0  \ \ \  f_r(\theta=0)=I $$
where the other 3 variables are fixed. It admits
an unique solution $f_r(\theta)$, which we can consider as
parametrised by the $r$, the radial coordinate on $\cpx$. Set
$f_r(2\pi)=F_r$ and note that the conjugacy class $[F_r]\subset SU(2)$
is gauge-invariant (see \cite{SS}, lemma 3.2). We ask if the limit:
\begin{equation} \label{limhol}
\lim_{r\seta\infty}[F_r]=[F]
\end{equation}
is well-defined as a conjugacy class in $SU(2)$. Since conjugacy
classes in $SU(2)$ are parametrised by the half-open interval $[0,1)$,
the limiting holonomy $[F]$ is just a real number $0\leq c<1$.

Under suitable conditions (see appendix B), it is reasonable to expect
that (\ref{limhol}) will be indeed well-defined. One can then ask how it
behaves under Nahm transform, trying to see how it is translated into
the transformed Higgs pair.

\vskip8pt
The task of understanding how the limiting holonomy and the
parabolic weights behave under Nahm transform probably involves a
more detailed study of the asymptotic behaviour of the connections
$A$ on the bundle $E$ and $B$ on the bundle $V$ (or, equivalently,
of the corresponding hermitian metrics).

\chapter{Spectral data} \label{spec} 

In this chapter we close our circle of ideas by showing a
correspondence between instantons and what we call {\em spectral
data}, i.e. pairs consisting of a complex curve
$S\hookrightarrow\dual\times\proj$ and a line bundle over it
$\call\seta S$. In the light of main theorem \ref{nahmthm}, the
existence of such correspondence should not be surprising, for a
similar correspondence, between Higgs pairs and curves
with a line bundle over it, was shown by Hitchin in \cite{H3} in the
smooth case and by Bottacin \cite{Bo} and Markman \cite{Mk} in the
meromorphic case. These ideas are developed in the first two sections.

The main topic of this chapter is the proof of our third main
result. It is carried out in section \ref{match}.


\section{The instanton spectral data} \label{instspec}

Our first step towards the main theorem \ref{specthm} is to
construct a complex curve $S\hookrightarrow\dual\times\proj$
associated to a holomorphic vector bundle $\ee\seta\tproj$ as
defined in the beginning of chapter \ref{back}. To do this,
we follow Friedman, Morgan \& Witten \cite{FMW}.

Recall from section \ref{poin} that a semi-stable rank 2 holomorphic vector
bundle over an elliptic curve with trivial determinant either
splits as a sum of line bundles or is the unique non-trivial
extension ${\mathbf F}_2$ of $\underline{\cpx}$ by itself, tensored
with a line bundle of order two. From section \ref{spch}, we know
that that $\ee|_{T_w}$ splits as a sum of flat line bundles for all
but finitely many points $w\in\proj$.

We will assume that the restriction of $\ee\seta\tproj$
to the elliptic fibres is semi-stable for all $w\in\proj$.
Moreover, $\ee$ is defined to be {\em good} if there is no $w\in\proj$
such that $\ee|_{T_w}=L_\ksi\oplus L_{\ksi}$, for some $\ksi$ of order
two in $\dual$. In particular, we assume that the asymptotic
state $\ksi_0$ is not of order 2. From now on we restrict ourselves
to such bundles, unless otherwise stated.

The motivation for this definition will be made clear later on: the
spectral curves associated to {\em good} bundles are smooth. Note also
that {\em good} bundles form an open dense subset of the moduli space of
bundles $\ee$.

The {\em instanton spectral curve} $S\hookrightarrow\dual\times\proj$ is defined
as follows:
\begin{equation} \label{def.spec}
S=\{(\ksi,w)\in\dual\times\proj\ | \ {\rm either}\ \ee|_{T_w}=L_\ksi\oplus L_{-\ksi}\
   {\rm or}\ \ee|_{T_w}={\bf F}_2\otimes L_\ksi\}
\end{equation}

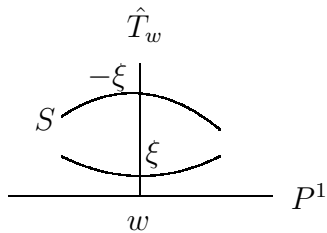
\begin{center} \begin{picture}(120,100)
            \put(0,20){\line(10,0){100}}
            \put(50,20){\line(0,10){50}}
            \put(45,7){$w$}
            \put(107,15){$P^1$}
            \put(45,80){$\hat{T}_w$}
            \put(52,33){$\xi$} \put(30,62){$-\xi$}
            \put(10,45){$S$}
            \qbezier(20,50)(50,70)(80,45)
            \qbezier(20,35)(50,20)(80,35)
\end{picture} \end{center}

Clearly, the natural projection $\pi_2:S\seta\proj$ is a branched double cover. More
precisely, for generic $w\in\proj$, $\pi^{-1}(w)=\{-\ksi,\ksi\}\in\dual\times\{w\}$.
There are then two types of branch points:
\begin{itemize}
\item those $w\in\proj$ for which $\ee|_{T_w}$ is indecomposable.
\item those $w\in\proj$ for which $\ee|_{T_w}$ splits as a sum of line bundles of
order two (i.e. $L_\ksi=L_{-\ksi})$;
\end{itemize}
Of course, the spectral curve associated to {\em good} bundles
$\ee$ only have branch points of the first type, since those of the
second type were excluded by definition.

Since $\ee$ is irreducible, there must be at least one branch point.
Its clear from the definition (\ref{def.spec}) that $S$ is a compact,
connected submanifold of $\dual\times\proj$ of complex dimension 1.
It inherits a complex structure from the chosen complex structure
on the ambient surface $\dual\times\proj$.

\begin{lem} \label{bp}
The map $\pi_2:S\seta\proj$ has $4k$ branch points, and the spectral
curve has genus $g(S)=2k-1$.
\end{lem}
\pf
This is an application of the Riemann-Roch theorem for the family
of Dolbeault operators $\del_w$ on $\ee|_{T_w}$, parametrised by
$\proj$. For generic $w\in\proj$, ${\rm dim}({\rm ker}\del_w)=0$;
this dimension jumps only when either $\ee|_{T_w}={\mathbf F}_2$ or
$\ee|_{T_w}=\underline{\cpx}\oplus\underline{\cpx}$ (again, this
second case is excluded from {\em good} bundles). From index
theory, we know that the number of jumping points is computed by
the first Chern class of the index bundle:
$$ c_1({\rm index}(\del_p))=\int_{\proj}
\left\{ch(\ee)td(p_1^*K_{T}^{-1})/[T]\right\}=-\int_{\tproj}c_2(\ee)=-k $$

This means that $\pi_1^{-1}(e)$ consists of $k$ points. Furthermore,
the points in the pre-image of $\pi_1$ of each element of order two
of $\dual$ are also branch points of $\pi_2$. As there are four such
points, we conclude that the covering map $S\seta\proj$ has $4k$ branch
points.

The second statement follows from the Riemann-Hurwitz formula.

Note however that branch points of the second type would count as a
double point, since the kernel of the Dolbeault operator of
$\underline{\cpx}\oplus\underline{\cpx}$ has dimension 2. For
instance, if there is exactly one point $p\in\proj$ such that
$\ee=\underline{\cpx}\oplus\underline{\cpx}$, then $\pi_1^{-1}(e)$
consists of $k-1$ points and there are $4k-4$ branch points
altogether. While this decreases the {\em real genus} of $S$, its
{\em virtual genus} is still $2k-1$.
\pfend

The curve $S$ admits an involution $\tau:S\seta S$ defined as
follows. Take $s\in S$ and let $w_s=\pi_2(s)$ and $\ksi_s=\pi_1(s)$
be its coordinates on $\dproj$; thus:
\begin{equation} \label{invol} \begin{array}{rcl}
\tau:S & \seta & S \\
(\ksi_s,w_s) & \mapsto & (-\ksi_s,w_s)
\end{array} \end{equation}
It is easy to see that the fixed points of $\tau$ are exactly the
branch points of the map $\pi_2:S\seta\proj$. Hence, $S/\tau$ is a
rational curve.

Once the topological type of $\ee$ is fixed, we show that, as we
vary the holomorphic structure on $\ee$, the respective spectral
curves lie within the same homology class in:
\begin{equation} \label{split.hom}
H_2(\dproj,\zed) = H_2(\dual,\zed) \oplus H_2(\proj,\zed)
\end{equation}

In fact, let $[p]$ be the generator of $H_2(\proj,\zed)$ and
$[\hat{t}]$ be the generator of $H_2(\dual,\zed)$. Regarding
$\dproj\stackrel{\pi_1}{\seta}\dual$ as a ruled surface, these can
be interpreted in $H_2(\dproj,\zed)$ as representing, respectively,
a fibre of $\pi_1$ and a constant section of $\pi_1$. They form a
basis for $H_2(\dproj,\zed)$, in which the intersection form looks
like:
$$ \left( \begin{array}{cc} 0 & 1 \\ 1 & 0 \end{array} \right) $$
Furthermore, the canonical divisor of $\dproj$ is given by
$K=-2[\hat{t}]$.

\begin{lem} \label{homclass}
As a homology class, $[S]=(k,2)$ in the (\ref{split.hom})
decomposition, and the map $\pi_1:S\seta\dual$ is a $k$-fold
branched covering map.
\end{lem}
\pf
S is a double cover of each fibre of the ruled surface, and we can
write the homology class of S as $[S]=2[p]+x[\hat{t}]$, for some
integer $x$. By the adjunction formula, we have:
\begin{eqnarray*}
g(S) & = & 1+\frac{1}{2}(K\cdot S+S^2) \\
2k-1 & = & 1+\frac{1}{2}(-4+4x) \\
\end{eqnarray*}
so $x=k$, as desired.

The second statement is now obvious. Note that the lemma could also
be proved by applying the proof of lemma \ref{bp} to the bundle
$\ee(\ksi)$ for each $\ksi\in\dual$.
\pfend

In other words, the topology of the bundle $\ee$ fixes the topology
of its spectral curve $S$. The holomorphic information is contained
on the choice of an embedding $S\hookrightarrow\dproj$ and on a line
bundle over $S$ that we now define.


\paragraph{Defining the line bundle over the spectral curve.}
The second part of our spectral data consists of a line bundle over
the spectral curve. Let $S$  be the spectral curve associated with
the holomorphic bundle $\ee\seta\tproj$, and consider the maps:
\begin{eqnarray}
T\times\dual \ \stackrel{\tau_1}{\longleftarrow} & T\times S &
\stackrel{\sigma}{\longrightarrow}\ S \nonumber \\
& \! \! \tau_2\downarrow & \label{maps} \\
& \tproj & \nonumber
\end{eqnarray}
where $\tau_1$ and $\tau_2$ are given by product of the identity on
the first factor and $\pi_1$ and $\pi_2$, respectively, on the
second factor. Clearly, $\tau_2$ is a double cover branched at 4k
elliptic curves $T_w\hookrightarrow\tproj$, where $w\in\proj$ are
the branch points of $\pi_2$. Furthermore, $\tau_1$ is also a
$k$-fold covering map.

We define a holomorphic line bundle $\call\seta S$ as follows:
\begin{equation} \label{lbdefn}
\call \ = \ \sigma_*(\tau_2^*\ee\otimes\tau_1^*{\bf P})
\end{equation}
where the subscript ``*'' denotes the direct image operation
on sheaves.

To identify the fibres of $\call$, denote $\ksi_s=\pi_1(s)$ and
$w_s=\pi_2(s)$, for $s\in S$. Relative Serre duality tells us that:
$$ \sigma_*(\tau_2^*\ee\otimes\tau_1^*{\bf P})^* \ = \
   R^1\sigma_*(\tau_2^*\ee\otimes\tau_1^*{\bf P}^*) $$
and this means that $\call^*=H^1(T_{w_s},\ee\otimes{\bf P}^*|_{T_{w_s}})$.
Thus, the fibre of $\call\seta S$ over $s\in S$ is given by:
\begin{equation} \label{lbdefn2}
\call_s=H^0(T_{w_s},\ee(\ksi_s)|_{T_{w_s}})
\end{equation}
If $\ee$ is {\em good}, it is easy to check that $\call_s$ is a
1-dimensional complex vector space for all $s\in S$, so that $\call$
is actually a line bundle. Otherwise, $\call$ is only a coherent
sheaf, since the dimension of (\ref{lbdefn2}) jumps at a finite
number of points; we will return to this point below.

\newpage

\begin{lem} \label{jac.spec}
The line bundle $\call$ has zero degree.
\end{lem}

\pf Look at the family of $\del$-operators on $T$ parametrised by $s\in S$:
$$ \del_s:\Lambda^0\ee(\ksi_s)|_{T_{w_s}}\seta\Lambda^{0,1}\ee(\ksi_s)|_{T_{w_s}} $$
and let ${\cal I}\in K(S)$ denote the corresponding {\em index bundle}.
Now, ${\rm det}\ {\cal I}$ is a genuine line bundle over $S$, with fibre
over $s\in S$ given by:
\begin{eqnarray*}
({\rm det}\ {\cal I})_s & = & \Lambda^{\rm max}({\rm ker}\del_s)\otimes
(\Lambda^{\rm max}({\rm coker}\del_s))^* \ = \\
& = & H^0(T_{w_s},\ee(\ksi_s)|_{T_{w_s}})\otimes H^1(T_{w_s},\ee(\ksi_s)|_{T_{w_s}})^* \ = \\
& = & H^0(T_{w_s},\ee(\ksi_s)|_{T_{w_s}})\otimes H^0(T_{w_s},\ee(-\ksi_s)|_{T_{w_s}})
\end{eqnarray*}
by Serre duality on $\ee(\ksi_s)|_{T_{w_s}}$. Thus
${\rm det}\ {\cal I}=\call\otimes(\tau^*\call)$, and
${\rm deg}\ {\cal I}=2{\rm deg}\ \call$.

Now, the degree of ${\cal I}$ can be computed via Riemann-Roch for
families, as follows:

\begin{eqnarray*}
{\rm deg}\ {\cal I} & = & ch\left( \tau_2^*\ee\otimes\tau_1^*{\bf P} \right)
                          td(T_FS)/[T\times S] \ = \\
                    & = & \left( 2 - k\cdot t\wedge(2s) \right) \cdot
                          \left( 1 + \tau_1^*c_1({\bf P}) + \frac{1}{2}(2t)\wedge(ks) \right)/
                          [T\times S] \ = \\
                    & = & 0
\end{eqnarray*}
as desired.
\pfend

\paragraph{Reconstructing the original bundle.}
We now want show how to reconstruct $\ee\seta\tproj$ from its
spectral pair $(S,\call)$ obtained as above, consisting of a
curve $S\hookrightarrow\dproj$ plus a line bundle $\call\seta S$
of degree 0. We define:
\begin{equation} \label{rconst}
\check{\ee}=\tau_{2*}(\tau_1^*{\bf P}\otimes\sigma^*\call^*)
\end{equation}
Clearly, $\check{\ee}$ is a locally free sheaf of rank 2.

\begin{prop} \label{rconstprop}
$\check{\ee}$ is holomorphically equivalent to $\ee$.
\end{prop}

\pf It is easy to see that $\check{\ee}$ and $\ee$ are topologically
equivalent, just by examining the effect of $\tau_{2*}$, $\tau_1^*$
and $\sigma^*$ on the Chern character of ${\bf P}$ and $\call$.

We want to show that there is a holomorphic bundle map
$\ee\stackrel{\varphi}{\seta}\check{\ee}$ whose determinant is
nowhere vanishing. In other words, $\varphi$ can be regarded as a
section in $H^0(\tproj,\ee\otimes\check{\ee})$, and $\det\varphi\in
H^0(\tproj,(\Lambda^2\ee)\otimes(\Lambda^2\check{\ee}))$. However,
$\Lambda^2\ee=\Lambda^2\check{\ee}=\underline{\cpx}$, so
$\det\varphi$ either vanishes identically or it is nowhere
vanishing. Thus, it is enough to verify that there is a section
\linebreak $\varphi\in H^0(\tproj,\ee\otimes\check{\ee})$ which is an
isomorphism at a single point \linebreak $(z,w)\in\tproj$.

The definition of $\call$ in (\ref{lbdefn}) gives us a canonical
identification:
\begin{equation} \label{canid}
\call\seta\sigma_*(\tau_2^*\ee\otimes\tau_1^*{\bf P})
\end{equation}
which can be interpreted as a canonical choice of section in
\linebreak $H^0(S,\call^*\otimes\sigma_*(\tau_2^*\ee\otimes\tau_1^*{\bf P}))$.
On the other hand, we have canonical identifications:
\begin{eqnarray*}
H^0(S,\call^*\otimes\sigma_*(\tau_2^*\ee\otimes\tau_1^*{\bf P})) & = &
H^0(T\times S,\sigma^*\call^*\otimes\tau_2^*\ee\otimes\tau_1^*{\bf P}) \ = \\
& = & H^0(\tproj,\tau_{2*}(\sigma^*\call^*\otimes\tau_1^*{\bf P})\otimes\ee)
\end{eqnarray*}
Thus, the identification (\ref{canid}) gives us a canonical choice
of a section \linebreak $\varphi\in
H^0(\tproj,\check{\ee}\otimes\ee)$ and according to the
observations made above is enough to check that this is an
isomorphism at one point.

Take $w\in\proj$ not a branch point of the spectral curve. Indeed,
it is then not difficult to see that $\varphi(z,w)$ is actually the
identity map on \linebreak $(L_{\ksi_w})_z\oplus(L_{-\ksi_w})_z$.
\pfend

\paragraph{Example: the Weierstrass $\wp$-function.}
The graph of the Weierstrass $\wp$-function:
\begin{eqnarray*}
& \wp:\dual\seta\proj & \\
& \Gamma_{\wp}=\{ (\ksi,w)\ | \ w=\wp(\ksi) \}
\end{eqnarray*}
is a curve of genus 1 inside $\dproj$. Clearly, projecting onto
each factor, $\Gamma_{\wp}$ is a 1-fold cover of $\dual$  and a
double cover of $\proj$, branched at 4 points. Together with any
line bundle of degree zero, $\Gamma_{\wp}$ can be used to construct
a {\em good} rank 2 holomorphic bundle $\ee\seta\tproj$, giving a
simple example of a charge 1 doubly-periodic instanton; the
asymptotic state can be chosen by changing the base point of $\wp$.

\paragraph{Relation with Fourier-Mukai transform.}
The spectral line bundle $\call$ can also be seen as a coherent sheaf
over $\dproj$ supported exactly over the spectral curve
$S\hookrightarrow\dproj$. Adopting this point of view, the appropriate
definition of $\call\seta\dproj$ is given by:
\begin{equation} \label{fmt}
\call^* = {\rm R}^1 p_{23*}(p_{12}^*\ee\otimes p_{13}^*{\bf P}^*)
\end{equation}
where $p_{ij}$ are the obvious projections of $\tproj\times\dual$
onto its factors.

The sheaf (\ref{fmt}) coincides with the so-called {\em
Fourier-Mukai transform} of the holomorphic vector bundle
$\ee\seta\tproj$ (see for instance \cite{Th} and the references
there). Proposition \ref{rconstprop} is then equivalent to the fact
that (\ref{rconst}) is inverse in a certain sense to (\ref{fmt}),
where these operations are regarded as functors acting between
certain derived categories over $\tproj$ and $\dproj$.

\paragraph{The geometry of the branch points.}
Let us now allow $S$ to have branch points of the second type. As
one approaches the branch points of $\pi_2:S\seta\proj$, the
behaviour of the spectral curve is roughly given by the pictures
below:

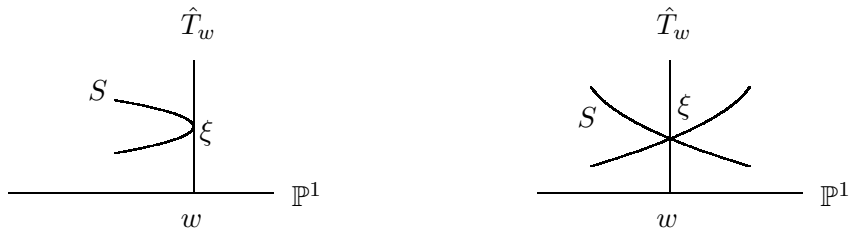
\begin{center} \small

\begin{picture}(400,100)

\put(40,0){\begin{picture}(120,100)
          \put(0,20){\line(10,0){100}}
          \put(70,20){\line(0,10){50}}
          \put(65,7){$w$}
          \put(107,15){$\proj$}
          \put(65,80){$\hat{T}_w$}
          \put(72,40){$\xi$}
          \put(30,55){$S$}
          \qbezier(40,35)(99,45)(40,55)
          \end{picture}}

\put(240,0){\begin{picture}(120,100)
            \put(0,20){\line(10,0){100}}
            \put(50,20){\line(0,10){50}}
            \put(45,7){$w$}
            \put(107,15){$\proj$}
            \put(45,80){$\hat{T}_w$}
            \put(53,50){$\xi$}
            \put(15,45){$S$}
            \qbezier(20,30)(70,45)(80,60)
            \qbezier(20,60)(30,45)(80,30)
            \end{picture}}

\end{picture} \\
Branch points of the spectral curve corresponding to \\
$E|_{T_w}={\bf F}_2\otimes L_\xi$ and $E|_{T_w}=L_\xi\oplus L_\xi$, respectively.
\end{center}

In other words, $S$ acquires a double-point over the points $w\in\proj$
for which $\ee|_{T_w}$ is a trivial extension of a line bundle of order
2 by itself. Moreover, $\call$ fails to be a genuine line bundle over $S$,
since the stalk over the double-point becomes 2-dimensional. Instead,
$\call$ is a coherent sheaf of degree 0 over the singular spectral curve.

Clearly, the presence of such points alters the genus of $S$, but not the
homology class within which $S$ lies. Furthermore, the bundle equivalence
established in proposition \ref{rconstprop} is still valid for bundles
$\ee$ which are not {\em good}.

We will show in the two following sections that the spectral curve
associated with a generic point in the moduli space of
doubly-periodic instantons must be smooth, i.e. there are no branch
points of the second type.


\section{Hitchin's spectral data} \label{hitspec}

We now look at the other side of the picture and study the spectral
curves coming from Higgs pairs. This time, our construction is
based on Hitchin's approach to non-singular Higgs pairs \cite{H2}.

Recall that $\vv\seta\dual$ is a holomorphic bundle of rank $k$,
and keeping in mind the holomorphic description of the Higgs field
discussed in section \ref{holo}, $\Phi$ is an endomorphism valued
$(1,0)$-form with simple poles at $\as$. Recall also that the
eigenvalues of the residues of $\Phi$ are non-vanishing. So, for
any fixed $\ksi\in\dual\setminus\{\as\}$, $\Phi(\ksi)$ is a
$k\times k$ matrix and one can compute its $k$-eigenvalues. As we
vary $\ksi$, we get  a $k$-fold covering, possibly branched, of
$\dual\setminus\{\as\}$ inside $\dual\times\cpx$. This {\em curve
of eigenvalues} is what we want to define as our spectral curve.

More precisely, we define the {\em Higgs spectral curve} to be the set:
\begin{equation} \label{h.spec}
C=\left\{ (\ksi,w)\in\dproj\ |\ {\rm det}(\Phi[\ksi]-w\cdot{\rm I}_k)=0 \right\}
\end{equation}
where $\proj=\cpx\cup\{\infty\}$. In other words, $C$ is the set of points
$(\ksi,w)\in\dproj$ such that $w$ is an eigenvalue of the endomorphism
$\Phi(\ksi):\vv_\ksi\seta\vv_\ksi$. Note in particular that the points
$(\as,\infty)$ belong to $C$ (with multiplicity one if $\ksi_0$ is not of
order 2).

\begin{prop} \label{smooth}
The spectral curve associated to a generic Higgs bundle $(\vv,\Phi)$ is smooth.
\end{prop}

\pf Let $Q\stackrel{p}{\seta}\dual$ be the line bundle with a section
$\sigma$ vanishing up to order 1 at $\as$. Thus, $\Psi=\Phi\otimes\sigma$
is a holomorphic section of ${\rm End}\vv\otimes Q\otimes K_{\dual}$.
Clearly, the value of $\Psi$ at $\as$ is a matrix of rank 1.

Usual Higgs bundle theory \cite{H2} yields a spectral curve
$C^\prime$ lying in the total space of the line bundle $Q$, which
we will denote by $X$. In other words, $C^\prime$ is the zero locus
of a section of $(\pi^*Q)^{\otimes k}$ given by the characteristic
polynomial of $\Psi$:
$$ \varphi \ = \ {\rm det}(\Psi-\lambda) \ = \
   \lambda^k+a_1\cdot\lambda^{k-1}+\dots+a_{k-1}\cdot\lambda+a_k $$
where $\lambda$ is a tautological section of the pull back of the line bundle
$Q\seta\dual$ to its total space, i.e. $p^*Q\seta X$. Since $\Psi(\as)$ is a
matrix of rank 1, the coefficients $a_2,\dots,a_k$ all have simple zeros
at $\as$. The coefficient $a_1(\as)$ is equal to the trace of $\Psi$ at these
points, which is simply given by its unique nonzero eigenvalue, i.e.
$a_1(\as)=\pm\epsilon$.

On the other hand, as the coefficients $a_1,\dots,a_k$ vary, the corresponding
zero locus $\{\varphi=0\}$ form a linear system of divisors on $X$, and hence
on its compactification $\overline{X}={\Bbb P}(Q\oplus\underline{\cpx})$. Since
$\lambda^k$ belongs to the system, any base point must lie in the 0-section of
$X$. So the base points of $|\{\varphi=0\}|$ are $\as$ in the 0-section of $X$,
since these are the only points where $a_k$ vanishes. Indeed, it is easy to see
that $a_k$ vanishes with order $k-1$ at $\as$.

Bertini's theorem guarantees that a generic element of the linear system is
smooth away from its base points, and it is singular there. In other words,
the spectral curve $C^\prime$ associated to a generic Higgs field $\Psi$ is
smooth away from $\as$ in the 0-section of $X$, which is a point of multiplicity
$k-1$.

We must now relate $C^\prime$ with our spectral curve $C$ defined in
(\ref{h.spec}). First note $\dproj$ can be obtained from $\overline{X}$
by performing {\em elementary transformations} based on $(\as,0)$
(see, for instance, \cite{L}). More precisely, we blow up
$(\as,0)\in\overline{X}$ and then blow down the proper transforms of
the fibres over $(\as,0)$. This gives a birational map
$\overline{X}\stackrel{\beta}{\dashrightarrow}\dproj$; we argue that
$C$ is the closure of $\beta(C^\prime)$, i.e. the {\em proper
transform} of $C^\prime$ under $\beta$.

Indeed, $\beta$ can also be represented as follows:
\begin{eqnarray*}
\overline{X} & \seta & \dproj \\
x & \seta & (p(x),(p^*\sigma)(x))=(p(x),\sigma(p(x)))
\end{eqnarray*}
Let $\dproj\stackrel{\pi}{\seta}\dual$ be the projection onto the first factor,
and denote by $\lambda'$ the tautological section of $\pi^*K_{\dual}$; clearly,
$\lambda=\lambda'\otimes\sigma$. If $x\in C^\prime$, then
$\left[{\rm det}(\Psi-p^*\lambda)\right](x)=0$, so that:
\begin{eqnarray*}
0 & = & {\rm det}(\Psi(p(x))-p^*\lambda(x)) \ = \\
  & = & {\rm det}\left( \Phi(p(x))\cdot\sigma(p(x))-p^*\lambda'(x)\cdot\sigma(p(x))\right) \ = \\
  & = & {\rm det}(\Psi(\pi(\beta(x)))-p^*\lambda'(x))\cdot\sigma(\pi(\beta(x)))^k \ = \\
  & = & \left[{\rm det}(\Phi-\pi^*\lambda')\right](\beta(x))\cdot\sigma(\pi(\beta(x)))^k
\end{eqnarray*}
$$ \Rightarrow\ \ \ \left[{\rm det}(\Phi-\pi^*\lambda')\right](\beta(x))=0 $$
Therefore, $\beta(x)\in C$ if $p(x)\neq\as$, since $\sigma(\pi(\beta(x)))$ vanishes
at these points.

The birational map $\beta$ is ill-defined on the fibres over $\as$; the situation
there is better understood by looking more closely at the elementary transformation.
Recall that $C^\prime$ has multiplicity $k-1$ at $(\as,0)$. After blowing up these
points, $\widetilde{C^\prime}$ (the proper transform of $C^\prime$) intersects the
exceptional divisor at $k-1$ generically distinct points. On the other hand,
$\widetilde{C^\prime}$ intersects $\widetilde{p^{-1}(\as)}$ (the proper transforms
of the fibres over $\as$) at a single point. Blowing down $\widetilde{p^{-1}(\as)}$
maps the exceptional divisors to the fibres of $\dproj$ over $\as$, so that
$C=\overline{\beta(C^\prime)}$ intersects $\pi^{-1}(\as)$ at generically $k$ distinct
points. This completes the proof, for $C$ is smooth elsewhere for generic Higgs
field $\Phi$. \pfend

In particular, it follows from the proof that all possible Higgs
spectral curves lie within the same linear system.

\paragraph{Defining the line bundle over the spectral curve.}
By definition, each point $c\in C$ corresponds to an eigenvalue of
$\Phi[\pi_1(c)]$. We define a line bundle ${\cal N}\seta C$ with
fibre over $c\in C$ given by the associated eigenspace. More precisely,
let $\ksi_c=\pi_1(c)$ and $w_c=\pi_2(c)$, and define:
$$ {\cal N}_c={\rm ker}\left\{ \Phi[\ksi_c]-w_c\cdot I_k \right\} $$
Generically, one expects the eigenvalues to be distinct, so that
${\cal N}$ is actually a line bundle.

\paragraph{Reconstructing the Higgs bundle.}
Conversely, the curve $C$ and the line bundle $\cal N$ determine $\vv$
and $\Phi$ over $\dual$. Indeed, Hitchin has shown that there is a torsion
sheaf ${\cal B}\seta\dual$ supported over the branch points of the $k$-fold
map $\pi_1:C\seta\dual$ such that:
$$ 0 \seta \oo_{\dual}(V)^* \seta \oo_{\dual}(\pi_{1*}{\cal N})^* \seta {\cal B} \seta 0 $$

Furthermore, the Higgs field $\Phi$ can be obtained as follows.
Pulling back $K_{\dual}$ to the spectral curve $C$ via the natural
$k$-fold covering map $\pi_1$ one obtains a tautological section
$\lambda$ in $H^0(C,\pi_1^*K_{\dual})$, the {\em section of
eigenvalues}. The operation of multiplication by $\lambda$ yields a
section of ${\rm End}(\pi_{1*}{\cal N})\otimes K_{\dual}$ which takes
$V$ to $V \otimes K_{\dual}$ and so defines the $\Phi\in{\rm End}V\otimes K_{\dual}$.

See \cite{H2} for more details.


\section{Matching the spectral data} \label{match}

So far we only know that our two spectral curves $S$ and $C$ lie
inside \linebreak $\dproj$ and that they have at least two points
in common, namely $(\as,\infty)$, since $\Phi$ has semi-simple
residues. We now show that if $(V,B,\Phi)$ is the Nahm transform of
$(E,A)$, then the instanton spectral curve $S$ associated to
$(E,A)$ actually coincides with the Higgs spectral curve $C$
associated to $(V,B,\Phi)$, thus proving our third main result.

Let us first consider an alternative definition of the transformed Higgs
field. Pick up the sections $s_0,s_\infty\in H^0(\proj,\oo_{\proj}(1))$,
as defined in section \ref{invholo}. For each $\ksi\in\dual$, we can
define the map:
\begin{equation} \label{alt.hig} \begin{array}{rcl}
H^1(\tproj,\ee(\ksi))\times H^1(\tproj,\ee(\ksi)) &
\stackrel{\Psi_\ksi}{\longrightarrow} & H^1(\tproj,\tilde{\ee}(\ksi)) \\
(\alpha,\beta) & \mapsto & \alpha\otimes s_0-\beta\otimes s_\infty
\end{array} \end{equation}

If $(\alpha,\beta)\in{\rm ker}\Psi_\ksi$, we define the
Higgs field $\Phi$ at the point $\ksi\in\dual$ as follows:
\begin{equation} \label{alt.hig2}
\Phi[\ksi](\alpha)=\beta
\end{equation}
It is easy to see that this is equivalent to our previous definition,
presented on section \ref{holo}.

Now suppose that $\alpha$ is an eigenvector of $\Phi[\ksi]$ with eigenvalue
$\epsilon$. In particular, the point $(\ksi,\epsilon)\in\dproj$ belongs
to the Higgs spectral curve $C$. By definition, we have:
$$ \Phi[\ksi](\alpha)=\epsilon\cdot\alpha \ \ \ \Rightarrow \ \ \
   \alpha\otimes(s_0-\epsilon\cdot s_\infty)=0 $$

Clearly, $s_\epsilon=s_0-\epsilon\cdot s_\infty$ is a holomorphic
section in $H^0(\proj,\oo_{\proj}(1))$ vanishing at $\epsilon\in\proj$.
So, it induces the following sheaf sequence:
$$ 0 \seta \ee(\ksi) \seta \widetilde{\ee(\ksi)} \seta
   \widetilde{\ee(\ksi)}|_{T_\epsilon} \seta 0 $$
which in turn induces the cohomology sequence:
\begin{equation} \label{sqc.xx} \begin{array}{cccccc}
0 & \seta & H^0(T_\epsilon,\widetilde{\ee(\ksi)}|_{T_\epsilon}) & \seta & & \\
  & \seta & H^1(\tproj,\ee(\ksi)) & \stackrel{\otimes s_\epsilon}{\seta} &
            H^1(\tproj,\widetilde{\ee(\ksi)}) & \stackrel{r}{\seta} \\
  & \stackrel{r}{\seta} & H^1(T_\epsilon,\widetilde{\ee(\ksi)}|_{T_\epsilon}) &  \seta & 0 &
\end{array} \end{equation}
Thus $\alpha\in{\rm ker}(\otimes s_\epsilon)=
H^0(T_\epsilon,\widetilde{\ee(\ksi)}|_{T_\epsilon})=
H^0(T_\epsilon,\ee(\ksi)|_{T_\epsilon})$.

In particular, $H^0(T_\epsilon,\ee(\ksi)|_{T_\epsilon})$ in non-empty, hence
either $\ee|_{T_\epsilon}=L_\ksi\oplus L_{-\ksi}$ or
$\ee|_{T_\epsilon}={\bf F}_2\otimes L_\ksi$. So, the point $(\ksi,\epsilon)\in\dproj$
also belongs to the instanton spectral curve $S$. Therefore, the two curves
$C$ and $S$ must coincide.

It also follows from the cohomology sequence (\ref{sqc.xx}) that the
$\epsilon$-eigenspace of $\Phi[\ksi]$ is exactly
$H^0(T_\epsilon,\widetilde{\ee(\ksi)}|_{T_\epsilon})=
H^0(T_\epsilon,\ee(\ksi)|_{T_\epsilon})$, i.e.
${\cal N}_{(\ksi,\epsilon)}=\call_{(\ksi,\epsilon)}$, and the spectral
bundles (or sheaves) also coincide.

This proves our main theorem \ref{specthm}. Note that the argument also
works if $\ee$ is not {\em good}.

In particular, we conclude that the instanton spectral curves lie within the same
linear system inside $\dproj$, and are smooth for a generic point in the moduli
space $\modspc^*$.


\section{The moduli space of spectral data}

Let ${\cal S}_{(k,\ksi_0)}$ denote the configuration space for the spectral
data $(S,\call)$. Let also $\Sigma_{(k,\ksi_0)}$ be the space spectral curves,
i.e. space of complex curves lying within the homology class
$(2,k)\in H_2(\dproj,\zed)$ and containing the points $(\as,\infty)\in\dproj$.
From section \ref{instspec}, it is easy to see that ${\cal S}_{(k,\ksi_0)}$
is the total space of a fibration over $\Sigma_{(k,\ksi_0)}$ whose fibres are
given by $\jj(S)$, the Jacobian of the curve $S\in\Sigma_{(k,\ksi_0)}$:
\begin{equation} \label{fib}
\jj \seta {\cal S}_{(k,\ksi_0)} \seta \Sigma_{(k,\ksi_0)}
\end{equation}

Let us compute the dimension of the space of spectral curves
$\Sigma_{(k,\ksi_0)}$. From Kodaira \cite{Kod}, we know that
deformations of a complex submanifold $S\hookrightarrow\dproj$ are
given by holomorphic sections of the normal line bundle $N_S$. On
the other hand, we want to keep the points $(\as,\infty)\in\dproj$
fixed. Thus, we are actually interested only on those elements of
$H^0(S,N_S)$ vanishing at these points. Hence:
\begin{equation} \label{dimbase}
{\rm dim}\ \Sigma_{(k,\ksi_0)} = {\rm dim}\ H^0(S,N_S)-1
\end{equation}

In order to compute the right hand side, we look at the following exact sequence:
$$ 0 \seta \oo_{\dproj} \seta \oo_{\dproj}(L_S) \seta \oo_S(N_S) \seta 0 $$
where by $L_S\seta\dproj$ we denoted the line bundle associated to the divisor
$S\hookrightarrow\dproj$. It induces the cohomology sequence ($M=\dproj$):
\small \baselineskip18pt
\begin{equation} \label{sqc.xy} \begin{array}{cccc}
0 & \seta H^0(M,\oo_M) \seta H^0(M,L_S) \seta H^0(S,N_S) \seta H^1(M,\oo_M) & \seta & \\
  & \seta H^1(M,L_S) \seta H^1(S,N_S) \seta H^2(M,\oo_M) \seta H^2(M,L_S) & \seta & 0
\end{array} \end{equation}
\normalsize \baselineskip18pt

By regarding $M=\dproj$ as a ruled surface over an elliptic curve,
we know that $H^2(M,\oo_M)=\{0\}$ (see \cite{Beau}, chapter 3).
Thus $H^2(M,L_S)$ must also vanish, and $h^0(L_S)-h^1(L_S)=2k+2$ by
Riemann-Roch for line bundles over surfaces.

On the other hand, we argue that $h^0(L_S)=2k+2$. Indeed, note that
$c_1(L_S)=2\cdot\hat{t}+k\cdot p$, so $L_S=p_1^*Q\otimes p_2^*\oo_{\proj}(k)$,
where $Q\seta\dual$ is line bundle of degree 2. Now, it follows from the Leray
spectral sequence that (see \cite{Beau}, chapter 3):
$$ \begin{array}{rccc}
H^0(\dproj,L_S)= & \underbrace{H^0(\dual,Q)} & \otimes &
\underbrace{H^0(\proj,\oo_{\proj}(k))} \\
& {\rm dim}=2 & & {\rm dim}=k+1
\end{array} $$
and the claim is now obvious.

Thus $h^1(L_S)=0$ and it follows from (\ref{sqc.xy}) that also $H^1(S,N_S)=\{0\}$.
In particular, one concludes that the deformation of spectral curves is
{\em unobstructed} \cite{Kod}. We are then left with:
\begin{equation} \label{sqc.zz} \begin{array}{rcccl}
0  \seta & \underbrace{H^0(M,\oo_M)} & \seta H^0(M,L_S) \seta H^0(S,N_S) \seta &
\underbrace{H^1(M,\oo_M)} & \seta 0 \\
& {\rm dim}=1 & & {\rm dim}=1 &
\end{array} \end{equation}
so that $h^0(M,L_S)=h^0(S,N_S)=2k+2$. It follows from (\ref{dimbase})
that ${\rm dim}\ \Sigma_{(k,\ksi_0)}=2k+1$. Thus,
$$ {\rm dim}\ {\cal S}_{(k,\ksi_0)}=
   {\rm dim}\ \Sigma_{(k,\ksi_0)}+{\rm dim}\ \jj(S)=4k $$
Furthermore, $\Sigma_{(k,\ksi_0)}$ is a smooth projective manifold,
since the deformation is unobstructed and all curves lie within the
same linear system. This implies that the whole moduli space of
spectral data ${\cal S}_{(k,\ksi_0)}$ is itself smooth and
projective.

Therefore, we conclude that the $\modspc^*$, the moduli space of
extensible instanton connections with fixed instanton number $k$
and asymptotic state $\as$, is a complex manifold of dimension $4k$,
containing ${\cal S}_{(k,\ksi_0)}$ as an open dense subset.

Finally, one would like to understand the action of $\torus\times S^1$
on $\modspc^*$ introduced in section \ref{actsec} in terms of the
fibration (\ref{fib}). We expect the torus translations $t_x^*$
to leave $\Sigma_{(k,\ksi_0)}$ invariant, acting only
on the jacobian fibres (by tensoring line bundles over $S$ with
$\pi_1^*L_z$). On the other hand, $\cpx\times S^1$ is
expected to preserve the fibres, acting only on the base space.

\paragraph{Conclusion.}
Summing up the work done so far, we note that the moduli spaces
of doubly-periodic instantons and the moduli space of singular Higgs
pairs are seen to be naturally identified via the construction of the
respective spectral data. The two moduli spaces are, in particular,
diffeomorphic. Since we know that the moduli of Higgs bundles is a
hyperk\"ahler manifold (once the parabolic structure and the residue
are fixed), one concludes that the moduli of instantons (with the
appropriate parameters fixed) is also hyperk\"ahler.


\section{Instantons and rational maps}

Donaldson has shown in \cite{D3} that monopoles are equivalent to
rational maps $\proj\seta\proj$. This was done via the equivalence
of monopoles and solutions of Nahm's equations obtained by Nahm
transform. It is reasonable to expect that a similar result should
hold for doubly-periodic instantons as well. As a by-product of the
spectral curve construction done above, we show that the space of
spectral curves $\Sigma_{(k,\ksi_0)}$ admits a parametrisation in
terms of rational maps.

First, recall that $\dual$ admits a $\zed_2$ action $\sigma$ (its
group involution), and that the quotient $\dual/\sigma$ is a
rational curve, which we denote by $\hproj$. Points in $\hproj$ can
be regarded as a pair of points $\{\pm\ksi\}\in\dual$. Moreover, it
is easy to see that the diagram:
\begin{eqnarray*}
\tau: S & \seta & S \\
\pi_1\downarrow \ & & \ \downarrow\pi_1 \\
\sigma: \dual & \seta & \dual
\end{eqnarray*}
commutes, where $\tau$ is the involution of the spectral curve defined
in (\ref{invol}).

So, let $\ee\seta\tproj$ be a {\em good} rank 2 holomorphic vector bundle
as above. We define a map $R:\proj\seta\hproj$ as follows.
Restricting $\ee$ to each elliptic fibre as in the construction of the
spectral curve, we get either $\ee|_{T_w}=L_\ksi\oplus L_{-\ksi}$ or
$\ee|_{T_w}={\bf F}_2\otimes L_\ksi$
We then define:
\begin{equation}
R(w)=[\pm\ksi_w]
\end{equation}
Lemma \ref{homclass} implies that $R$ has degree $k$. Fixing the
asymptotic state means fixing the image of $\infty$ under the map $R$.

The involution $\sigma:\dual\seta\dual$ can be seen as acting on the
product $\dproj$, with quotient $\hproj\times\proj$. Under this
quotient, the spectral curve is mapped to $S/\tau\hookrightarrow\hproj\times\proj$.
It is then easy to see that $\Gamma_R\hookrightarrow\proj\times\hproj$,
the graph of $R$, coincides with $S/\tau$. In particular, this implies
that $R$ is a rational map.

Recovering the spectral curve from the rational map $R$ is not
hard. Let $p_\sigma:\dproj\seta\hproj\times\proj$ be the projection
map naturally associated with the quotient $(\dual/\sigma)\times\proj$.
It is easy to see that $p_\sigma^{-1}(\Gamma_R)\hookrightarrow\hproj\times\proj$
coincides with the spectral curve $S$ associated with $\ee$.

In other words, we have shown that:

\begin{thm} \label{ratmap}
There is a bijective correspondence between $\Sigma_{(k,\xi_0)}$,
the space of instanton spectral curves, and rational maps
$R:\proj\seta\hproj$ of degree $k$ and such that $R(\infty)=[\as]$.
\end{thm}

It is easy to see that the set of rational maps as above is indeed
parametrised by $2k+1$ complex numbers. The map $R:\proj\seta\hproj$
has the form:
$$ \frac{a_kw^k + a_{k-1}w^{k-1} + \dots + a_0}{b_kw^k + b_{k-1}w^{k-1} + \dots + b_0} $$
which gives $2k+2$ parameters. Now fixing $R(\infty)=[\as]$ means
fixing the ratio $a_0/b_0$, killing the extra degree of freedom.

\appendix 

\chapter{Relative Index Theorem} \label{apprit}

Let $X$ be a connected, complete riemannian manifold, possibly
non-compact. Let $K\subset X$ be a compact subset and denote
$\Omega=X\setminus K$.

Consider complex vector bundles $E_0\seta X$ and $E_1\seta X$
and pick up two first-order, elliptic differential operators
$D_0:L^2_1(E_0)\seta L_2(E_0)$ and \linebreak $D_1:L^2(E_1)\seta L^2(E_1)$.
Suppose that there is a bundle isomorphism \linebreak
$F:E_0|_\Omega\seta E_1|_\Omega$.

We define the {\em relative topological index} of $D_0$ and $D_1$,
which we denote by $ind_t(D_1,D_0)$. First, if $X$ is a compact
manifolds, then we define $ind_t(D_1,D_0)=\ind(D_1)-\ind(D_0)$.
If not, we proceed as follows. Cut the set $\Omega$ out of $X$
along the hypersurface $M=\partial\Omega$ and compactify
$X$ by sewing on another compact manifold $\widetilde{\Omega}$
with boundary $M$; in particular, we can take $\widetilde{\Omega}$
to be the closure of $X\setminus K$. Extend $D_0$ and $D_1$ to
elliptic pseudo-differential operators $\widetilde{D_0}$ and
$\widetilde{D_1}$ over $\widetilde{X}$. Then, we define:
\eq \label{indt}
ind_t(D_1,D_0)=\ind(\widetilde{D_1})-\ind(\widetilde{D_0})
\end{equation}
a quantity that can be computed using the Atiyah-Singer index theorem.

It can be shown that the above expression is independent of the choice
of $\widetilde{\Omega}$ and of how the operators $D_0$ and $D_1$ are
extended to $\widetilde{D}_0$ and $\widetilde{D}_1$ (see lemma
\ref{lem} below). Note also that if $X$ is odd dimensional, then
$ind_t(D_1,D_0)=0$. Moreover, it is clear that perturbations of $D_0$
and $D_1$ supported at $\Omega$ leave $\ind(\widetilde{D}_1)$ and
$\ind(\widetilde{D}_0)$ unchanged.

Now suppose that $D_0$ and $D_1$ are Fredholm operators when acting between
the spaces considered above. We define the {\em relative analytical index}
as follows:
$$ \label{inda} ind_a(D_1,D_0)=\ind(D_1)-\ind(D_0) $$

We want to show that, under certain conditions, these relative indices coincide.
Let us start by reviewing some standard facts. Recall that if $D$ is a Fredholm
operator, there is a bounded, elliptic pseudo-differential operator $Q$, called
the {\em parametrix} of $D$, such that $DQ=I-S$ and $QD=I-S^\prime$, where $S$
and $S^\prime$ are compact smoothing operators, and $I$ is the
identity operator. Note that neither $Q$ nor $S$ and $S^\prime$ are unique.

In particular, there is a bounded operator $G$, called the Green's operator for $D$,
satisfying $DG=I-H$ and $GD=I-H^\prime$, where $H$ and $H^\prime$ are finite rank
projection operators $H:L^2_p(E)\seta{\rm ker}(D)$ and
$H^\prime:L^2(E)\seta{\rm coker}(D)$.

Let $K^H(x,y)$ be the Schwartzian kernel of the operator $H$. Its {\em trace function} is
defined by ${\rm tr}[H](x)=K^H(x,x)$; moreover, these are $C^\infty$ functions \cite{A}.
If $D$ is Fredholm, its index is given by:
\eq \label{ind} \ind(D)=\int_X\left({\rm tr}[H]-{\rm tr}[H^\prime]\right) \end{equation}
as it is well-known; recall that compact operators have smooth, square integrable kernels.
Furthermore, if $X$ is a closed manifold, we have \cite{A}:
\eq \label{ind2} \ind(D)=\int_X\left({\rm tr}[S]-{\rm tr}[S^\prime]\right) \end{equation}

Let us now return to the situation set up above. Consider the parametrices
and Green's operators ($j=0,1$):
\eq \label{ops} \begin{array}{ccc}
\left\{ \begin{array}{l} D_jQ_j=I-S_j \\ Q_jD_j=I-S^\prime_j \end{array} \right. &\ \ &
\left\{ \begin{array}{l} D_jG_j=I-H_j \\ G_jD_j=I-H^\prime_j \end{array} \right.
\end{array} \end{equation}

The two operators $D_0$ and $D_1$ are said to coincide at $\Omega$ if
\linebreak $D_0|_\Omega=F\circ(D_1|_\Omega)\circ F^{-1}$.
We are finally in position to state our relative index theorem:

\begin{thm} \label{rit}
Let $D_0$ and $D_1$ be first-order, elliptic pseudo-differential Fredholm operators
over a complete riemannian manifold $X$ as above and suppose that they
coincide at $\Omega$. Then $ind_a(D_1,D_0)=ind_t(D_1,D_0)$.
\end{thm}

The first step is to express the indices involved in terms of integral formulas. As in
(\ref{ind}), we have for the analytical index that:
\begin{equation} \label{inda2}
\ind_a(D_1,D_0)=\int_{X}\left({\rm tr}[H_1]-{\rm tr}[H_1^\prime]\right)-
\int_{X}\left({\rm tr}[H_0]-{\rm tr}[H_0^\prime]\right)
\end{equation}

For the relative topological index, we have the following lemma:
\begin{lem} \label{lem}
Under the hypothesis of the theorem, we have that:
\begin{equation} \label{mid}
ind_t(D_1,D_0)=\int_{X}\left({\rm tr}[S_1]-{\rm tr}[S_1^\prime]\right)-
\int_{X}\left({\rm tr}[S_0]-{\rm tr}[S_0^\prime]\right)
\end{equation} \end{lem}

\pf Compactify $X$ as explained above; one obtains the compact manifold
$\widetilde{X}$. Extend $D_0$ and $D_1$ to operators $\widetilde{D}_0$ and
$\widetilde{D}_1$, both defined over the whole $\widetilde{X}$. Let
$\widetilde{Q},\widetilde{Q}^\prime$ denote the extension of each $Q_j,Q_j^\prime$
from $\Omega$ to $\widetilde{\Omega}$, which are, by hypothesis, equal. Choose
cut-off functions $\beta_1,\beta_2:\widetilde{X}\seta\real$ such that:
\eq \label{betas} \begin{array}{cc}
(\beta_1)^2+(\beta_2)^2=1 & {\rm supp}\beta_1^j=K\
{\rm and\ supp}\beta_1= \widetilde{\Omega}
\end{array} \end{equation}
Suppose also that the differentials $d\beta_1,d\beta_2$ are supported in a small
neighbourhood of $M$. One can glue each $Q_j$ with $\widetilde{Q}$ using the cut-off
functions to obtain parametrix $\widetilde{Q}_j$ for $\widetilde{D}_j$ over the whole
$\widetilde{X}$. More precisely, let $s\in\Gamma(\widetilde{E})$:
\eq \label{para}
\widetilde{Q}_j(s)=\beta_1Q_j(\beta_1s)+\beta_2Q(\beta_2s)
\end{equation}
It is straightforward to verify that these are truly parametrix for
$\widetilde{D}_j$ and that:
\eq \label{comps}
\left\{ \begin{array}{l}
\widetilde{S}_j=\beta_1S_j(\beta_1s)+\beta_2S(\beta_2s)+
d\beta_1.Q_j(\beta_1s)+d\beta_2.Q(\beta_2s) \\
\widetilde{S}_j^\prime(s)=\beta_1S_j^\prime(\beta_1s)+\beta_2S^\prime(\beta_2s)
\end{array}
\right. \end{equation}
hence $\widetilde{S}_0,\widetilde{S}_1$ and $\widetilde{S}^\prime_0,\widetilde{S}^\prime_1$
coincide at $\widetilde{\Omega}=\widetilde{X}\setminus K$. Thus
\eq \label{xx1}
{\rm tr}[S_1]-{\rm tr}[S_1^\prime]-{\rm tr}[S_0]+{\rm tr}[S_0^\prime]=0
\end{equation}
at $\widetilde{\Omega}$. From (\ref{ind2}):
\eq \label{xx2}
\ind(\widetilde{D}_j)=\int_{\widetilde{X}}
\left({\rm tr}[\widetilde{S}_j]- {\rm tr}[\widetilde{S}_j^\prime]\right)
\end{equation}
and (\ref{mid}) follows immediately from the definition (\ref{indt}),
(\ref{xx1}) and (\ref{xx2}).
\pfend

As we noted before, this lemma shows also that the definition of relative topological
index is independent of the choice of extensions $\widetilde{D}_0$ and $\widetilde{D}_1$;
this is quite clear from (\ref{mid}).

Before we step into the proof of theorem \ref{rit} itself, we must introduce some notation.
Let $f:[0,1]\seta[0,1]$ be a smooth function such that $f=1$ on $[0,\frac{1}{3}]$, $f=0$ on
$[\frac{2}{3},1]$ and $f'\approx-1$ on $[\frac{1}{3},\frac{2}{3}]$. Pick up a point $x_0\in X$
and let $d(x)={\rm dist}(x,x_0)$. For each $m\in\zed^*$, consider the functions:
\eq \label{beta1}
f_m(x)=f\left(\frac{1}{m}e^{-d(x)}\right)
\end{equation}
Note that ${\rm supp}d(f_m)^{\frac{1}{2}}\subset B_{\log\frac{3}{4m}}-B_{\log\frac{3}{2m}}$ and
\eq \label{beta2}
||\nabla f_m||_{L^2}\leq\frac{C}{m}
\end{equation}
where $C=\left(\int_Xe^{-d(x)}\right)^{\frac{1}{2}}$. Here, $B_r=\{x\in X\ |\ d(x)\leq r\}$,
which is compact by the completeness of $X$.

{\em Proof of theorem \ref{rit}:} All we have to do is to show that
the right  hand sides of (\ref{inda2}) and (\ref{mid}) are equal.
In fact, let $V^*\subset V$ be small neighbourhoods of the diagonal
of $(X\times X)$ and choose $\psi\in C^\infty(X\times X)$
supported on $V$ and such that $\psi=1$ on $V^*$. Let $Q_j$ be
the operator whose Schwartzian kernel is
$K^{Q_j}(x,y)=\psi(x,y)K^{G_j}(x,y)$, where $G_j$ is the Green's
operator for $D_j$. Then $Q_j$ is a parametrix for $D_j$ with:
$$ \begin{array}{ccc} D_jQ_j=I-S_j & {\rm and} & Q_jD_j=I-S_j^\prime \end{array} $$
and clearly:
\eq \label{traces} \begin{array}{ccc}
{\rm tr}[S_j]={\rm tr}[H_j] & {\rm and} & {\rm tr}[S_j^\prime]={\rm tr}[H_j^\prime]
\end{array} \end{equation}

But is not necessarily the case that the two parametrix $Q_0$ and
$Q_1$ so obtained coincide at $\Omega$. We will glue them with $Q$,
the common parametrix of $D_0|_\Omega$ and $D_1|_\Omega$ using the
cut-off functions $f_m$ defined above (assume that the base points
are contained in the compact set $K$). More precisely:
\eq \label{para2}
Q_j^{(m)}(s)=(f_m)^{\frac{1}{2}}Q_j((f_m)^{\frac{1}{2}}s)+
(1-f_m)^{\frac{1}{2}}Q((1-f_m)^{\frac{1}{2}}s)
\end{equation}
which now coincide at $\Omega$. For the respective compact operators, we get:
$$ \left\{ \begin{array}{l}
S_j^{(m)}(s)=(f_m)^{\frac{1}{2}}S_j((f_m)^{\frac{1}{2}}s)+
(1-f_m)^{\frac{1}{2}}S((1-f_m)^{\frac{1}{2}}s)+\\
\hspace{2cm} +d(f_m)^{\frac{1}{2}}.(Q_j((f_m)^{\frac{1}{2}}s)-Q((1-f_m)^{\frac{1}{2}}s) \\
S_j^{(m)\prime}(s)=(f_m)^{\frac{1}{2}}S_j^\prime((f_m)^{\frac{1}{2}}s)+
(1-f_m)^{\frac{1}{2}}S^\prime((1-f_m)^{\frac{1}{2}}s) \end{array} \right. $$
therefore:
$$ \begin{array}{l} {\rm tr}[S_j^{(m)\prime}]-{\rm tr}[S_j^{(m)}]= \\
(f_m)^{\frac{1}{2}}({\rm tr}[S_j^\prime]-{\rm tr}[S_j])+
(1-f_m)^{\frac{1}{2}}({\rm tr}[S^\prime]-{\rm tr}[S])+{\rm tr}[d(f_m)^{\frac{1}{2}}.(Q_j-Q)]
\end{array} $$
and
$$ \begin{array}{c}{\rm tr}[S_1^{(m)\prime}]-{\rm tr}[S_1^{(m)}]-{\rm tr}[S_0^{(m)\prime}]+
{\rm tr}[S_0^{(m)}]= \\

=(f_m)^{\frac{1}{2}}\left({\rm tr}[S_1^\prime]-{\rm tr}[S_1]-{\rm tr}[S_0^\prime]+
{\rm tr}[S_0]\right)+ \\ +\underbrace{{\rm tr}[d(f_m)^{\frac{1}{2}}(Q_1-Q)]-
{\rm tr}[d(f_m)^{\frac{1}{2}}(Q_0-Q)]} \\
={\rm tr}[d(f_m)^{\frac{1}{2}}(Q_1-Q_0)] \end{array} $$

We must now integrate both sides of the expression above and take limits as
$m\seta\infty$. For $m$ sufficiently large, ${\rm supp}(1-f_m)\subset\Omega$ the LHS equals
$ind_t(D_1,D_0)$ by lemma \ref{lem}; on the other hand, the term inside the parenthesis on
the RHS equals $ind_a(D_1,D_0)$ by (\ref{traces}) and (\ref{inda2}). Thus, it is enough to show
that the last two terms on the RHS vanishes as $m\seta\infty$. Indeed, note that:
$$ {\rm tr}[d(f_m)^{\frac{1}{2}}(Q_1-Q_0)]=d(f_m)^{\frac{1}{2}}{\rm tr}[(Q_1-Q_0)] $$
hence, since ${\rm supp}(df_m)\subset\Omega$ for sufficiently large $m$ and using also
(\ref{beta2}), it follows that:
$$ \int_{\Omega}{\rm tr}[d(f_m)^{\frac{1}{2}}(Q_1-Q_0)]\leq
\frac{C}{m}\int_{\Omega}{\rm tr}[(G_1-G_0)]\seta0\ {\rm as}\ m\seta\infty $$
if the integral on the RHS is finite.

Indeed, let $D=D_i|_\Omega$; from the parametrix equation, we have:
$$ D(G_1|_\Omega-G_0|_\Omega)=H_1|_\Omega-H_0|_\Omega $$
Observe that $W={\rm ker}(H_1|_\Omega-H_0|_\Omega)$ is a closed
subspace of finite codimension in $L^2(\Omega)$. Moreover
$W\subseteq{\rm ker}(D)$; thus,
$(G_1|_\Omega-G_0|_\Omega)$ has finite dimensional range and hence
it is of trace class.

This concludes the proof. \pfend

\paragraph{Applications.}
In our applications, we have a Fredholm operator $D_1$ and an
invertible operator $D_0$. However, they do not exactly coincide
away from a compact set; instead, they are {\em asymptotically
equal}, i.e. given $\epsilon>0$, there is a compact set $K\subset
X$ such that:
$$ ||D_1-D_0||^2_{L^2(X\setminus K)}<\epsilon $$
In order to apply theorem \ref{rit}, we construct a new Fredholm operator
$D^\prime_1$ as follows. Let $\beta_1$ and $\beta_2$ be cut-off functions,
respectively supported over $K$ and $X\setminus K$ as before, and define:
$$ D^\prime_1=\beta_1 D_1 \beta_1 + \beta_0 D_0 \beta_0 $$
Now, it is clear that $D^\prime_1|_{X\setminus K}$ coincides with
$D_0|_{X\setminus K}$. Furthermore, since
$||D^\prime_1-D_1||_{L^2(X)}<\epsilon$ with $\epsilon$ arbitrarily
small, we know that $\ind(D_1^\prime)=\ind(D_1)$.

So, theorem \ref{rit} applies for the pair of operators $D_1^\prime$
and $D_0$. Since $\ind(D_0)=0$, one concludes
$\ind(D_1)=\ind(D_1^\prime)=ind_t(D_1^\prime,D_0)$. In this situation,
$D_0$ is often referred to as the {\em background operator}.

\paragraph{Final example.}
We conclude by treating one example particularly relevant to the
index problems we deal with in the bulk of the present work; see
also \cite{GL}. Suppose $X$ is a spin manifold and let $D$ be its
canonical Dirac operator acting on positive spinors over $X$.
Suppose that $E\seta X$ is a complex vector bundle of rank $n$
which is trivialised outside a compact subset of $X$. Let
$\underline{\cpx}^n$ denote the trivial complex bundle of rank $n$,
and consider the operators:
$$ \left\{ \begin{array}{l}
D_0:\Gamma(\underline{\cpx}^n\otimes S^+)\seta\Gamma(\underline{\cpx}^n\otimes S^-) \\
D_1:\Gamma(E\otimes S^+)\seta\Gamma(E\otimes S^-)
\end{array} \right. $$
Clearly, these operators coincide outside the support of $E$; thus:
\begin{eqnarray*}
ind_t(D_1,D_0) & = & ind_t(D_1,D_0)=\ind(\widetilde{D_1})-\ind(\widetilde{D_0})=  \\
& = & \left\{ ch(E)\cdot\hat{{\bf A}}(\widetilde{X}) \right\}[\widetilde{X}]-
\left\{ ch(\underline{\cpx}^n)\cdot\hat{{\bf A}}(\widetilde{X}) \right\}[\widetilde{X}]= \\
& = & \left\{ (ch(E)-n)\cdot\hat{{\bf A}}(\widetilde{X}) \right\}[\widetilde{X}]
\end{eqnarray*}


\chapter{On the asymptotic behaviour of extensible connections} \label{appextn}
\chaptermark{Extensibility}

Motivated by the properties of the inverse transformed bundle and
instanton connection, it seems fair to make the following conjecture:

\begin{conj} \label{dec.extn}
If $|F_A|\sim O(|w|^{-2})$ then there is a holomorphic
vector bundle $\ee\seta\tproj$ such that
$$ \ee|_{T\times(\proj\setminus\{\infty\})}\simeq(E,\del_A) $$
In other words, $A$ is extensible.
\end{conj}

Such conjecture motivates other questions, which we will not attempt
to address here:
\begin{itemize}
\item Do all anti-self-dual connections on $E\seta\torus$ with
      finite energy with respect to the Euclidean metric satisfy
      $|F_A|\sim O(|w|^{-2})$?
\item Does the converse holds, i.e. if $A$ is extensible then
      $|F_A|\sim O(|w|^{-2})$? If not, what are the necessary and
      sufficient analytical conditions for extensibility (in terms of
      the Euclidean metric)?
\item Given a holomorphic bundle $\ee\seta\tproj$, is there a
      connection $A$ on $\ee|_{T\times(\proj\setminus\{\infty\})}$
      such that $A$ is anti-self-dual and $|F_A|\sim O(|w|^{-2})$
      with respect to the Euclidean metric?
\end{itemize}
Note however that if the conjecture does hold, the Nahm transform
constructed in the bulk of the thesis would give a positive answer,
though a rather indirect one, to the last question. However, it
would be rather interesting to obtain a direct proof.

We would like to point out that the techniques applied to the solution
of this problem would probably extend to instanton connection on bundles
over surfaces of the form $\Sigma\times\cpx$, where $\Sigma$ is any
compact complex curve.


\paragraph{Ingredients for a proof.}
The key ingredient for a possible proof \ref{dec.extn} is the
following $L^p$ integrability result due to Buchdahl \cite{Bu}:

\begin{lem} \label{buch}
Let $\Delta$ be a unit polydisc in $\cpx^2$. Let $A$ be a matrix
valued $(0,1)$-form on $\Delta$ with coefficients in $L^p_j(\Delta)$,
where $p>2$ and $j\geq1$, such that $\del A+A\wedge A=0$. Then there
is a matrix-valued function $h\in L^p_{j+1}(\Delta)$, possibly defined
on a smaller polydisc, such that $\del h=-Ah$.
\end{lem}

The strategy is to use lemma \ref{buch} to construct local holomorphic
extensions of $E$, and then patch them together to give a global holomorphic
extension $\ee$.

More precisely, let $U\subset T$ be a small open set, with complex coordinate
$z$; and let $D_R\subset\cpx$ be the complement of a disc of large radius
$R\gg0$, with complex coordinate $w$. Define:
$$ \begin{array}{ccc}
\Delta_0=U\times(B_1(0)\setminus\{0\}) & \ \ {\rm and} \ \ &
\Delta=U\times B_1(0)
\end{array} $$
and consider the inversion map:
\begin{equation} \label{invmap} \begin{array}{ccc}
\iota:\Delta_0 & \seta &  U\times D_R \\
(z',w') & \seta & \left(z=z',w=\frac{R}{w'}\right)
\end{array} \end{equation}
It is also convenient to introduce polar coordinates for the above complex
coordinates:
\begin{eqnarray*}
w'=(\rho,\theta) & \stackrel{\iota}{\mapsto} & w=\left(r=\frac{R}{\rho},\theta\right)
\end{eqnarray*}
and this implies that:
\begin{eqnarray*}
dr=-\frac{d\rho}{\rho^2} & {\rm and} & d\rho=-\frac{dr}{r^2}
\end{eqnarray*}

In order to use Buchdahl's lemma (or some of its versions), one would have to
establish following gauge fixing lemma:

\begin{conj} \label{gaugefixlem}
If $|F_{\check{A}}|\sim O(|w|^{-2})$ then, for $R\gg0$, there is a gauge
$g:U\times D_R\seta SU(2)$ such that $\iota^*g(A)\in L^p_1(\Delta_0)$, $p>2$.
\end{conj}

This is a familiar problem in gauge theory, and there are various
results along these lines, see for instance \cite{DK}, \cite{U},
\cite{R}. The fact that we have a pointwise estimate on the curvature,
instead of some global $L^p_k$ bound, makes the conjecture possibly
easier to prove than the hard results mentioned above.

Now consider the local trivialisation of $E|_{U\times D_R}$ corresponding
to the gauge obtained in the above conjecture. Define
$F=\iota^*E|_{U\times D_R}\seta\Delta_0$ and $A'=\iota^*g(A)$.
Thus, by \ref{buch} and \ref{gaugefixlem}, we can find a gauge
$h\in L^p_2(\Delta_0)$ ($p>2$), possibly after shrinking
$\Delta^{(n)}_0$ if necessary, such that:
$$ h(A')=h^{-1}(A')h+h^{-1}\del h $$
is a (1,0)-form. Note that there are many functions satisfying the
above equation, for if $h$ is one, so is $hf$ for any holomorphic
matrix-valued function $f$ on $\Delta_0$. Since $\iota^*A$
vanishes at $\{w'=0\}$, we see that $h(z',0)$ is holomorphic in
$z'=z$. Thus, we suppose without loss of generality that $h$ is
the identity over $\{w'=0\}$, for we can always take
$h(z,w')\cdot h^{-1}(z,0)$ instead, if necessary.

Now let $g_2=(\iota^*g)h$. In this new gauge, the connection $\iota^*A$
is represented by a $(1,0)$-form. Thus, $g_2$ is a holomorphic basis for
$F$. We extend $F$ holomorphically over $\{w'=0\}$ by defining
$g_2$ as a holomorphic basis on $\overline{F}\seta\Delta$.

We must now show how to patch these local extensions together and produce
a global holomorphic extension of $E$ over $T_\infty$.

Let $U$ and $W$ be any two intersecting neighbourhoods in $T$. It
suffices to show that the transition function $\Psi$ for the gauges $g_2(U)$
and $g_2(W)$ on $\overline{F}_U\seta U\times B_1(0)$ and
$\overline{F}_W\seta W\times B_1(0)$, respectively, constructed as above does
extend to a holomorphic function on $(U\cap W)\times B_1(0)$. Let \linebreak
$g_2(U)=g_2(W)\Psi_{UW}$ be such a transition function; $\Psi_{UW}$ is defined
and holomorphic on $(U\cap W)\times(B_1(0)\setminus\{0\})$. If it can be extended
holomorphically over $\{w'=0\}$, the cocycle condition will follow from continuity
of the transition functions and the cocycle condition for $E$.

Let $\iota^*g(U)=\iota^*g(W)\Upsilon_{UW}$, where $\Upsilon_{UW}$
is a transition function for the original gauges. The gauges $\iota^*g(U)$
and $\iota^*g(W)$ are continuous, hence so is $\Upsilon_{UW}$.

On the other hand, we have:
\begin{equation} \label{trans}
g_2(U)=\iota^*g(U)\Psi_U=\iota^*g(W)\Upsilon_{UW}\Psi_U
=\iota^*g(W)\Psi_W^{-1}\Upsilon\Psi_U
\end{equation}
Since $\Psi_W$ and $\Psi_U$ are bounded and continuous, so is the matrix function
$\Psi_{UW}=\Psi_W^{-1}\Upsilon_{UW}\Psi_U$. But $\Psi_{UW}$ is holomorphic on
$(U\cap W)\times(B_1(0)\setminus\{0\})$, so it extends holomorphically over
$(U\cap W)\times B_1(0)$, as desired.

In other words, {\em quadratic curvature decay implies extensibility up to
the gauge fixing lemma \ref{gaugefixlem}}.


\section{Proof of the proposition \ref{goodgauge}}

Recall that we need to establish the following result
\footnote{I thank Olivier Biquard for showing me the arguments 
in this section}:

\begin{prop} \label{goodgauge2}
If $|F_A|\sim O(r^{-2})$, then, for $R$ sufficiently large, there
is a gauge over $T\times V_R$ and a constant flat connection $\Gamma$
on a topologically trivial rank two bundle over the elliptic curve
such that:
$$ |A-p^*\Gamma| \ = \ |\alpha| \ \sim O(r^{-1}\cdot\log r) $$
\end{prop}

First, we need the following lemma that we shall assume without
proof:

\begin{lem} \label{t3gauge}
Let $B$ be a connection on a rank two bundle over
$T^3=S^1\times S^1\times S^1$ satisfying  
$|F_B|\leq\epsilon$ for $\epsilon$ sufficiently small.
Choose $L_1$, $L_2$, $L_3$ such that $\exp(-2\pi L_k)$ 
is the monodromy of $B$ at the point $(0,0,0)$ around 
the $k^{\rm th}$-circle. Then there exists an unique gauge $g$ 
on $S^1\times S^1\times S^1$, such that:
\begin{enumerate}
\item $g(0,0,0)=I$;
\item $g(A)=M_1d\theta_1 + M_2d\theta_2 + M_3d\theta_3$, where:
	\begin{itemize}
	\item $M_1(\theta_1,0,0)=L_1$, $M_2(0,\theta_2,0)=L_2$, $M_3(0,0,\theta_3)=L_3$;
	\item $M_2(\theta_1,\theta_2,0)$ does not depend on $\theta_2$;
	\item $M_3(\theta_1,\theta_2,\theta_3)$ does not depend on $\theta_3$;
	\end{itemize}
\item in this gauge, one has the control:
$$ \sup\left\{|M_i-L_i|,\left|[M_i,M_j]\right|\right\} \leq c\cdot\epsilon $$
\end{enumerate} \end{lem}

Now, fix a ray $\{x_0\}\times\{y_0\}\times[R,\infty)\times\{\theta_0\}$
and trivialise the bundle $E\seta\torus$ on this ray by parallel
transport. Therefore we have fixed a gauge on this ray.

Using lemma \ref{t3gauge} on each 3-dimensional tori
$T\times\{r\}\times S^1$, where $r>R$, we extend the above gauge
to a global gauge on $T\times V_R$. This is the gauge we are
looking for.

Indeed, let $B_r=A|_{T\times S^1_r}$, then $|F_{B_r}|<C\cdot r^{-1}$
(we have to account for the fact that one circle is getting larger).
for some constant $C$. By lemma \ref{t3gauge}, for each $r$, we can
find a gauge on $T\times S^1_r$ and a constant connection
\footnote{Note that $a$, $b$ and $h$ are respectively $L_1$, $L_2$
and $L_3$ in the statement of lemma \ref{t3gauge}.}:
$$ \Gamma_r = a(r)dx+b(r)dy+h(r)d\theta $$
such that $|B_r-\Gamma_r|<C/r$.

Now it follows from the curvature bound that:
\begin{eqnarray} \label{decs}
a\sim O(r^{-1})+a_\infty\ \ & \ \ b\sim O(r^{-1})+b_\infty
\ \ & \ \ h\sim O(\log r)+c
\end{eqnarray}
Therefore, the torus components are well-defined limits
as $r\seta\infty$, which we denoted by $a_\infty$ and $b_\infty$, respectively.
Defining $\Gamma=a_\infty dx+b_\infty dy$, we have:
$$ \Gamma_r = a_\infty dx + b_\infty dy + \gamma(r), \ \ {\rm where}\ \
   \gamma(r)\sim O(r^{-1}\cdot\log r) $$
Thus:
\begin{equation} \label{cont1}
|B_r-\Gamma|<C\cdot\frac{\log r}{r}
\end{equation}
and note that $\Gamma$ is flat by the estimate in (3) of lemma \ref{t3gauge}.

On the other hand, the connection $A$ can now be written in the
global gauge as follows:
$$ A = a(x,y,r,\theta)dx + b(x,y,r,\theta)dy + f(x,y,r,\theta)dr + h(x,y,r)d\theta $$
$$ {\rm such\ that}\ \ f(0,0,r,0)=0 $$
A lemma due Biquard (lemma 1 in \cite{B5}) implies that
$\partial h/\partial r$ and:
$$ \frac{\partial a}{\partial r}(x,0,r,0) \ \ \ \ \ \ \ \ 
   \frac{\partial b}{\partial r}(x,y,r,0) $$
are controled by the curvature bound.
From this control and from the curvature bound, one can
deduce a control on the following terms (which can be regarded
as the curvature of the connection $A$ restricted to each of 
the three circles plus the radial derivatives):
$$ \frac{\partial f}{\partial x} + [a,f] \ \ \ \ \ \ \
   \frac{\partial f}{\partial y} + [b,f] \ \ \ \ \ \ \
   \frac{\partial f}{\partial \theta} + [f,h] $$
Now diagonalising $a$, $b$ and $h$ one at a time allows us to control
each summand of the three terms above separately, thus controlling $f$:
the third term gives an estimate to $\frac{\partial f}{\partial \theta}$, 
so it is enough to control $f(x,y,r,0)$; now the second term gives an 
estimate to $\frac{\partial f}{\partial y}$, so it is enough to control 
$f(x,0,r,0)$ and this is finally done using the first term.
In fact, $f\sim O(r^{-1})$.

Together with (\ref{cont1}), this concludes the proof. \pfend

Note that the gauge fixing result needed to prove extensibility from the
curvature bound would require much more delicate arguments in order to
give estimates on the derivatives of the connection components. 


\bibliographystyle{plain} \bibliography{trbib}

\begin{thebibliography}{10}

\bibitem{A}
G.~Arfken.
\newblock {\em Mathematical methods in physics}.
\newblock Academic Press, New York-London, 1966.

\bibitem{A3}
M.~Atiyah.
\newblock Vector bundles over an elliptic curve.
\newblock {\em Proc. London Math. Soc.}, 7:414--452, 1957.

\bibitem{ADHM}
M.~Atiyah, V.~Drinfel'd, N.~Hitchin, and Y.~Manin.
\newblock Construction of instantons.
\newblock {\em Phys. Lett. A}, 65:185--187, 1978.

\bibitem{Beau}
A.~Beauville.
\newblock {\em Complex algebraic surfaces}.
\newblock Cambridge University Press, Cambridge, 1996.

\bibitem{B5}
O.~Biquard.
\newblock Prolongement d'un fibr\'e holomorphe hermitien \`a courbure $l^p$ sur
  une courbe ouverte.
\newblock {\em Internat. J. Math.}, 3:441--453, 1992.

\bibitem{B4}
O.~Biquard.
\newblock Sur les \'equations de {N}ahm et la structure de {P}oisson des
  alg\'ebres de {L}ie semi-simple complexes.
\newblock {\em Math. Ann.}, 304:253--276, 1996.

\bibitem{B2}
O.~Biquard.
\newblock Sur les fibr\'es paraboliques sur une surface complexe.
\newblock {\em J. London Math. Soc. (2)}, 53:302--316, 1996.

\bibitem{B3}
O.~Biquard.
\newblock Fibr\'es de {H}iggs et connexions int\'egrables: le cas logarithmique
  (diviseur lisse).
\newblock {\em Ann. Scient. \'Ec. Norm. Sup. (4)}, 30:41--96, 1997.

\bibitem{BB}
B.~Booss and D.~Bleecker.
\newblock {\em Topology and analysis - The {A}tiyah-{S}inger index formula and
  gauge theoretic physics}.
\newblock Springer-Verlag, 1985.

\bibitem{BT}
R.~Bott and L.~Tu.
\newblock {\em Differential forms in algebraic topology}.
\newblock Springer-Verlag, New York, 1982.

\bibitem{Bo}
F.~Bottacin.
\newblock Symplectic geometry on the moduli space of stable pairs.
\newblock {\em Ann. \'Ecole Norm. Sup. (4)}, 28:391--433, 1995.

\bibitem{BVB}
P.~Braam and P.~van Baal.
\newblock Nahm's transform for instantons.
\newblock {\em Comm. Math. Phys.}, 122:267--280, 1989.

\bibitem{Bu}
N.~Buchdahl.
\newblock Hermitian-einstein connections.
\newblock {\em Math. Ann.}, 280:625--648, 1988.

\bibitem{D3}
S.~Donaldson.
\newblock Nahm's equations and classification of monopoles.
\newblock {\em Commun. Math Phys.}, 96:387--207, 1984.

\bibitem{DK}
S.~K. Donaldson and P.~Kronheimer.
\newblock {\em Geometry of four-manifolds}.
\newblock Clarendon Press, 1990.

\bibitem{FMW}
R.~Friedman, J.~Morgan, and E.~Witten.
\newblock Vector bundles and {F}-theory.
\newblock {\em Commun. Math. Phys.}, 187:679--743, 1997.

\bibitem{GM}
H.~Garland and M.~Murray.
\newblock Kac-{M}oody monopoles and periodic instantons.
\newblock {\em Commun. Math. Phys.}, 120:335--351, 1988.

\bibitem{GR}
I.~Gradshteyn, I.~Ryzhik, and A.~Jeffrey (editor).
\newblock {\em Table of integrals, products and series}.
\newblock Academic Press, 1994.

\bibitem{GL}
M.~Gromov and H.~B. Lawson.
\newblock Positive scalar curvature and the index of the {D}irac operator on
  complete riemannian manifolds.
\newblock {\em Inst. des Hautes \'Etudes Scientifiques Publ. Math.},
  58:295--408, 1983.

\bibitem{H3}
N.~Hitchin.
\newblock Construction of monopoles.
\newblock {\em Commun. Math. Phys.}, 89:145--190, 1983.

\bibitem{H}
N.~Hitchin.
\newblock The self-duality equations on a {R}iemann surface.
\newblock {\em Proc. London Math. Soc.}, 55:59--126, 1987.

\bibitem{H2}
N.~Hitchin.
\newblock Stable bundles and integrable systems.
\newblock {\em Duke Math. J.}, 54:91--114, 1987.

\bibitem{HM}
J.~Hurtubise and M.~Murray.
\newblock On the construction of monopoles for the classical groups.
\newblock {\em Commun. Math. Phys.}, 122:35--89, 1989.

\bibitem{Ka}
A.~Kapustin.
\newblock Solution of {N}=2 gauge theories via compactification to three
  dimensions.
\newblock {\em Preprint hep-th/9804069}, 1998.

\bibitem{KS}
A.~Kapustin and S.~Sethi.
\newblock Higgs branch of impurity theories.
\newblock {\em Adv. Theor. Math. Phys.}, 2:571--592, 1998.

\bibitem{Kod}
K.~Kodaira.
\newblock A theorem of completeness of characteristic systems for analytic
  families of compact submanifolds of complex manifolds.
\newblock {\em Ann. of Math.}, 75:146--162, 1962.

\bibitem{Ko}
H.~Konno.
\newblock Construction of the moduli space of stable parabolic higgs bundles on
  a {R}iemann surface.
\newblock {\em J. Math. Soc. Japan}, 45:253--276, 1993.

\bibitem{K}
A.~Kovalev.
\newblock {\em The geometry of dimensionally reduced anti-self-duality
  equations}.
\newblock PhD thesis, Oxford, 1995.

\bibitem{L}
H.~Lange.
\newblock On elementary transformations of ruled surfaces.
\newblock {\em J Reine Angew. Math.}, 346:32--35, 1984.

\bibitem{Mk}
E.~Markman.
\newblock Spectral curves and integrable systems.
\newblock {\em Compositio Math.}, 93:255--290, 1994.

\bibitem{MS}
V.~Mehta and C.~Seshadri.
\newblock Moduli of vector bundles on curves with parabolic structure.
\newblock {\em Math. Ann.}, 248:205--239, 1980.

\bibitem{N}
W.~Nahm.
\newblock Self-dual monopoles and calorons.
\newblock In Denardo et~al., editor, {\em Trieste Conference on magnetic
  monopoles}, Berlin, 1980. Springer-Verlag.

\bibitem{Nkj}
H.~Nakajima.
\newblock Hyperk\"ahler structures on the moduli spaces of parabolic higgs
  bundles on riemann surfaces.
\newblock In {\em Lecture Notes in Pure and Applied Mathematics, 179}, pages
  199--208, New York, 1996. Dekker.

\bibitem{R}
J.~Rade.
\newblock On singular {Y}ang-{M}ills fields {III}: global theory.
\newblock {\em Internat. J. Math.}, 5:491--521, 1994.

\bibitem{SS}
L.~Sibner and R.~Sibner.
\newblock Classification of singular {S}obolev connections by their holonomy.
\newblock {\em Commun. Math. Phys.}, 114:337--371, 1992.

\bibitem{S}
C.~Simpson.
\newblock Harmonic bundles on noncompact curves.
\newblock {\em J. of Am. Math. Soc.}, 3:713--770, 1990.

\bibitem{Th}
R.~P. Thomas.
\newblock Mirror symmetry and actions of braid groups on derived categories.
\newblock {\em Preprint}, 1999.

\bibitem{U}
K.~Uhlenbeck.
\newblock Removable singularities in {Y}ang-{M}ills fields.
\newblock {\em Comm. Math. Phys.}, 83:11--29, 1982.

\bibitem{Wa}
G.~N. Watson.
\newblock {\em A treatise on the theory of Bessel functions}.
\newblock Cambridge University Press, 1995.

\end{thebibliography}

\end{document}